\def\DateTime{22/April/2005, 17:00(JP)}
\def\Version{Version $1.01$}

\def\no{\if01}
\def\iftwelvept{\no}
\def\ifdedicatory{\no}
\def\ifusepdf{\no}
\def\ifpsfont{\no}

\iftwelvept
\documentclass[leqno,12pt]{amsart}
\else
\documentclass[leqno]{amsart}
\fi
\usepackage{amssymb}
\usepackage{amscd}
\usepackage{latexsym}
\usepackage{verbatim}
\usepackage{graphics}
\usepackage[all]{xy}

\ifusepdf
\usepackage{hyperref}
\else\fi
\ifpsfont
\usepackage[T1]{fontenc}
\usepackage{times}
\else\fi

\iftwelvept
\else\fi

\theoremstyle{plain}
\newtheorem{Theorem}{Theorem}[section]
\newtheorem{Proposition}[Theorem]{Proposition}
\newtheorem{Lemma}[Theorem]{Lemma}
\newtheorem{Corollary}[Theorem]{Corollary}

\newtheorem{Claim}{Claim}[Theorem]

\theoremstyle{definition}
\newtheorem{Definition}[Theorem]{Definition}
\newtheorem{Remark}[Theorem]{Remark}

\renewcommand{\theTheorem}{\arabic{section}.\arabic{Theorem}}
\renewcommand{\theClaim}{\arabic{section}.\arabic{Theorem}.\arabic{Claim}}
\renewcommand{\theequation}{\arabic{section}.\arabic{Theorem}.\arabic{Claim}}

\def\rom{\textup}
\newcommand{\ZZ}{{\mathbb{Z}}}

\newcommand{\NN}{{\mathbb{N}}}
\newcommand{\AAA}{{\mathbb{A}}}
\newcommand{\OO}{{\mathcal{O}}}
\newcommand{\BB}{{\mathcal{B}}}

\newcommand{\Proj}{\operatorname{Proj}}

\newcommand{\Hom}{\operatorname{Hom}}
\newcommand{\Ext}{\operatorname{Ext}}
\newcommand{\Sing}{\operatorname{Sing}}

\newcommand{\Ker}{\operatorname{Ker}}
\newcommand{\Coker}{\operatorname{Coker}}
\newcommand{\Spec}{\operatorname{Spec}}

\newcommand{\Supp}{\operatorname{Supp}}
\newcommand{\Map}{\operatorname{Map}}
\newcommand{\codim}{\operatorname{codim}}
\newcommand{\mult}{\operatorname{mult}}
\newcommand{\ord}{\operatorname{ord}}

\newcommand{\lformal}{[\![}
\newcommand{\rformal}{]\!]}

\newcommand{\Sym}{\operatorname{Sym}}
\newcommand{\Irr}{\operatorname{Irr}}
\newcommand{\pushout}{\boxplus}
\newcommand{\marking}{\operatorname{mark}}

\newcommand{\Proof}{{\sl Proof.}\quad}
\newcommand{\QED}{{\unskip\nobreak\hfil\penalty50\quad\null\nobreak\hfil
{$\Box$}\parfillskip0pt\finalhyphendemerits0\par\medskip}}
\newcommand{\rest}[2]{\left.{#1}\right\vert_{{#2}}}
\newcommand{\labelrightrightarrow}[2]%
{\raisebox{1.0ex}{$\mathrel{\mathop{\stackrel{{\scriptstyle #1}}{{\scriptstyle%
\longrightarrow}}}\limits_{\mathrel{\mathop{\longrightarrow}\limits_{{\scriptstyle #2}}}}}$}}

\begin{document}

\title[Rigidity of morphisms for log schemes]%
{Rigidity of morphisms for log schemes}
\address{Department of Mathematics, Faculty of Science,
Kyoto University, Kyoto, 606-8502, Japan}
\author{Atsushi Moriwaki}
\email{moriwaki@math.kyoto-u.ac.jp}
\date{\DateTime, (\Version)}
\ifdedicatory
\dedicatory{Dedicated to Professor Masaki Maruyama on his 60th birthday}
\else\fi
\begin{abstract}
In this paper, we give the rigidity theorem for a log morphism
as an extension of a fixed scheme morphism.
We also give several applications of the rigidity theorem.
\end{abstract}


\maketitle

\setcounter{tocdepth}{1}
\tableofcontents


\section*{Introduction}
\renewcommand{\theTheorem}{\Alph{Theorem}}

In the paper \cite{IwaMw}, we proved Kato's conjecture, that is,
the finiteness of dominant rational maps in the category of log schemes
as a generalization of Kobayashi-Ochiai theorem.
It guarantees the finiteness of $K$-rational points of
a certain kind of log smooth schemes for a big function field $K$,
which gives rise to an evidence for Lang's conjecture.
In the proof of the above theorem, the most essential part is
the rigidity theorem of log morphisms.
In this paper, we would like to generalize it to a semistable scheme over
an arbitrary noetherian scheme.

Let $f : X \to S$ be a scheme of finite type over a locally noetherian
scheme $S$. 
We assume that $f : X \to S$ is a semistable scheme over $S$,
that is, $f$ is flat and, for any morphism $\Spec(\Omega) \to S$ with
$\Omega$ an algebraic closed field,
the completion of the local ring 
of $X \times_{S} \Spec(\Omega)$ at every closed point
is isomorphic to a ring of the type
\[
\Omega \lformal X_1, \ldots, X_n \rformal/(X_1 \cdots X_l).
\]
Let $g : Y \to S$ be another semistable scheme over $S$, and
let $\phi : X \to Y$ be a morphism over $S$.
Let $M_X$, $M_Y$ and $M_S$ be fine log structures on $X$, $Y$ and $S$
respectively.
We assume that $(X, M_X)$ and $(Y, M_Y)$ are log smooth and integral over
$(S, M_S)$ and $\phi$ is admissible with respect to $M_Y/M_S$,
that is, for all $s \in S$ and any irreducible components $V$ of the geometric
fiber $X_{\bar{s}}$ over $s$,
\[
\phi_{\bar{s}}(V) \not\subseteq \Supp(M_{Y_{\bar{s}}}/M_{\bar{s}}),
\]
where $M_{Y_{\bar{s}}} = \rest{M_Y}{Y_{\bar{s}}}$,
$M_{\bar{s}} = \rest{M_S}{\bar{s}}$ and
\[
\Supp(M_{Y_{\bar{s}}}/M_{\bar{s}}) =
\{ y \in Y_{\bar{s}} \mid \text{$M_{\bar{s}} \times \OO_{Y,\bar{y}}^{\times}
\to M_{Y_{\bar{s}}, \bar{y}}$ is not surjective} \}.
\]
Then, the following theorem is one of the main results of this paper.

\begin{Theorem}[Rigidity theorem]
\label{thm:rigid:log:morphism:intro}
If we have log morphisms
\[
(\phi, h) : (X, M_X) \to (Y, M_Y)
\quad\text{and}\quad
(\phi, h') : (X, M_X) \to (Y, M_Y)
\]
over $(S, M_S)$ as extensions of $\phi : X \to Y$,
then $h = h'$.
\end{Theorem}

As corollary of the rigidity theorem, we have the following descent theorem
for log morphisms:

\begin{Theorem}
\label{thm:descent:log:smooth:intro}
\begin{list}{}{\itemindent=\parindent \leftmargin=0ex \parsep=1ex}
\item[\rom{(1)}]\rom{(Descent for a faithfully flat base change)}
Let $\pi : S' \to S$ be a faithfully flat and quasi-compact morphism
of locally noetherian schemes.
Let $X' = X \times_S S'$, $Y' = Y \times_S  S'$ and 
$\phi' = \phi \times_S \operatorname{id}_{S'}$, and
let us set the induced morphisms as follows:
\[
\begin{CD}
X @<{\pi_X}<< X' \\
@V{f}VV @VV{f'}V \\
S @<{\pi}<< S'
\end{CD}
\qquad\qquad
\begin{CD}
Y @<{\pi_Y}<< Y' \\
@V{g}VV @VV{g'}V \\
S @<{\pi}<< S'
\end{CD}
\]
If there is a log morphism
\[
(\phi', h') : (X', \pi_X^*(M_X)) \to (Y', \pi_Y^*(M_Y))
\]
over $(S', \pi^*(M_S))$, then it descends to
a log morphism
\[
(\phi, h) : (X, M_X) \to (Y, M_Y)
\]
over $(S, M_S)$.

\item[\rom{(2)}]\rom{(Descent for the smooth topology)}
Let 
\[
\{ \pi_i : X_i \to X \}_{i \in I}
\quad\text{and}\quad
\{ \mu_j : Y_j \to Y \}_{j \in J}
\]
be families of smooth morphisms such that
$X = \bigcup_{i \in I} \pi_i(X_i)$ and $Y = \bigcup_{j \in J} \mu_j(Y_j)$.
We assume that, for each $i \in I$, there are $j \in J$ and
a morphism $\phi_i : X_i \to Y_j$ such that
the diagram
\[
\begin{CD}
X_i @>{\phi_i}>> Y_j \\
@V{\pi_i}VV @VV{\mu_j}V \\
X @>{\phi}>> Y
\end{CD}
\]
is  commutative and $\phi_i$ extends to
a log morphism 
\[
(\phi_i, h_i) : (X_i, \pi_i^*(M_X))  \to (Y_j, \mu_j^*(M_Y)).
\]
Then, there is a log morphism
\[
(\phi, h) : (X, M_X) \to (Y, M_Y)
\]
as an extension of $\phi : X \to Y$ with
the following diagram commutative:
\[
\begin{CD}
(X_i,  \pi_i^*(M_X))@>{(\phi_i, h_i)}>> (Y_j, \mu_j^*(M_Y))\\
@V{(\pi_i, \operatorname{nat})}VV @VV{(\mu_j,  \operatorname{nat})}V \\
(X, M_X) @>{(\phi, h)}>> (Y, M_Y),
\end{CD}
\]
where $\operatorname{nat}$ is the natural homomorphism.
\end{list}
\end{Theorem}

For the proofs of Theorem~\ref{thm:rigid:log:morphism:intro} and 
Theorem~\ref{thm:descent:log:smooth:intro},
our starting point is the following local structure theorem,
which asserts the local description of integral and smooth
log morphisms of semistable schemes.

\begin{Theorem}[Local structure theorem]
\label{thm:local:structure:theorem:intro}
Let $(f, h) : (X, M_X) \to (S, M_S)$ be a smooth and  integral morphism of fine log schemes.
Let $x$ be an element of $X$ and $s = f(x)$.
We assume that $f : X \to S$ is semistable at $x$, i.e.,
$f$ is flat at $x$ and, for any morphism $\eta : \Spec(\Omega) \to S$ with
$\Omega$ an algebraic closed field and $\eta(0) = s$,
the completion of the local ring 
of $X \times_{S} \Spec(\Omega)$ at every closed point lying over $x$
is isomorphic to a ring of the type
\[
\Omega \lformal X_1, \ldots, X_n \rformal/(X_1 \cdots X_l).
\]
Then, we have the following for each case:
\begin{enumerate}
\renewcommand{\labelenumi}{\Roman{enumi}.}
\item
If $f$ is smooth at $x$, then $\overline{M}_{X,\bar{x}} \simeq
\overline{M}_{S,\bar{s}} \times \NN^a$
for some non-negative integer $a$.

\item
If $f$ is not smooth at $x$ and 
$\bar{h}_{\bar{x}} : \overline{M}_{S,\bar{s}} \to
\overline{M}_{X,\bar{x}}$ splits, then
$\overline{M}_{X,\bar{x}} \simeq
\overline{M}_{S,\bar{s}} \times N$, where
$N$ is the monoid arising from monomials of 
\[
\ZZ[U_1,U_2, \ldots, U_a]/(U_1^2 - U_2^2)
\]
for some $a \geq 2$.

\item
If $f$ is not smooth at $x$ and
$\bar{h}_{\bar{x}} : \overline{M}_{S,\bar{s}} \to
\overline{M}_{X,\bar{x}}$ does not split, then
$\overline{M}_{X,\bar{x}}$ has the unique semistable structure 
$(\sigma, q_0, \Delta, B)$ over $\overline{M}_{S,\bar{s}}$
for some $\sigma \subseteq \overline{M}_{X,\bar{x}}$ with $\#(\sigma) \geq 2$,
$q_0 \in \overline{M}_{S,\bar{s}}$ and $\Delta, B \in \NN^{\sigma}$.
\end{enumerate}
\end{Theorem}

Based on the local structure theorem,
the proof of the rigidity theorem is carried out as follows:
Clearly we may assume that $S = \Spec(A)$ for some
noetherian local ring $(A, m)$.
First, we establish the theorem in the case where $A$ is an algebraically closed field.
This was proved actually in the previous paper \cite{IwaMw}. 
Next, by induction on $n$, we see that the assertion holds for the case $S = \Spec(A/m^n)$.
Finally, using the Krull intersection theorem, we can conclude its proof.
As an application, we give a criterion for the extension of a scheme
morphism to the log morphism
(cf. Theorem~\ref{thm:extension:log:mor}).
Besides them, we show the following extension problem of automorphisms to
the central fiber, which is a generalization of the result in
\cite{DM}.

\begin{Theorem}
\label{thm:log:canonical:isom;intro}
Let $(A, tA)$ be a discrete valuation ring. Let
$f : X \to \Spec(A)$ and $f' : X' \to \Spec(A)$ be proper and generically smooth semistable 
schemes over $\Spec(A)$, and let $B$ and $B'$ be horizontal effective Cartier divisors
on $X$ and $X'$ respectively. 
Let $\phi : X \dasharrow X'$ be a birational map over $\Spec(A)$ such that
it is an isomorphism on the generic fibers $X_{\eta}$ and $X'_{\eta}$ and that
$\phi_{\eta}^*(B'_{\eta}) = B_{\eta}$.
We assume that
$(X, B)$ and $(X', B')$ are of log canonical type 
\rom{(}for the definition of log canonical type, see \S\S\rom{\ref{subsec:log:canonical:auto}}\rom{)}
and that $\omega_{X/A} + B$ and $\omega_{X'/A} + B'$ are ample
with respect to $f$ and $f'$.
Then $\phi$ extends to an isomorphism.
\end{Theorem}

In \S1, we give the definition of semistable schemes and show their several properties.
In \S2, we recall several facts concerning log schemes and prove
the local structure theorem.
\S3 contains the proof of the rigidity theorem.
In \S4, we consider a descent problem for log schemes.
\S5 contains the explicit form of the determinant of the sheaf of log 1-forms.
In \S6, we consider a criterion for the extension of a scheme
morphism to the log morphism. Finally, in \S7, the extension problem of automorphisms
to the central fiber are considered.

\renewcommand{\thesubsubsection}{\arabic{subsubsection}}

\bigskip
\subsection*{Conventions and terminology}
Here we will fix several conventions and terminology for this paper.

\subsubsection{}
\label{CT:commutative:ring}
Throughout this paper, a ring means a commutative ring with the unity.

\subsubsection{}
\label{CT:log:structure}
In this paper, the logarithmic
structures of schemes means the sense of J.-M Fontaine, L. Illusie, 
and K. Kato. For the details, we refer to \cite{KatoLog}.
For a log structure $M_X$ on a scheme $X$,
we denote the quotient $M_{X}/\mathcal{O}_{X}^{\times}$
by $\overline{M}_{X}$.

\subsubsection{}
\label{CT:germ:etale:zariski}
Let $X$ be a scheme and $F$ a sheaf in the \'{e}tale topology.
For a point $x \in X$,
the germ of $F$ at $x$ with respect to the Zariski topology 
(resp. the \'{e}tale topology)
is denoted by $F_x$ (resp. $F_{\bar{x}}$).

\subsubsection{}
\label{CT:commutative:monoid}
Throughout this paper, a monoid is a commutative monoid with
the unity.
The binary operation of a monoid is often written additively.
We say a monoid $P$ is {\em finitely generated} if there are
$p_1, \ldots, p_n$ such that 
$P = \NN p_1 + \cdots + \NN p_r$.
Moreover, $P$ is said to be {\em integral} if
$x + z = y + z$ for $x, y, z \in P$, then $x = y$.
An integral and finitely generated monoid is said to be {\em fine}.
We say $P$ is {\em sharp} if $x + y = 0$ for $x, y \in P$,
then $x = y = 0$. For a sharp monoid $P$, an element $x$ of $P$
is said to be {\em irreducible} if
$x = y + z$ for $y, z \in P$, then either $y = 0$ or $z = 0$.
A homomorphism $f : Q \to P$ of monoids is said to be {\em integral} if
$f(q) + p = f(q') + p'$ for $p, p' \in P$ and $q, q' \in Q$,
then there are $q_1, q_2 \in Q$ and $p'' \in P$
such that $q + q_1 = q' + q_2$, $p = f(q_1) + p''$
and $p' = f(q_2) + p''$.
Note that an integral homomorphism of sharp monoids is injective.
Moreover, we say an injective homomorphism
$f : Q \to P$ {\em splits} if there is a submonoid $N$ of $P$ with
$P = f(Q) \times N$.

\subsubsection{}
\label{CT:map:N:support}
Let $X$ be a set. We denote the set of all maps $X \to \NN$ by
$\NN^X$.  
For $T \in \NN^X$,
$\Supp(T)$ is given by $\{ x \in X \mid T(x) > 0 \}$.
Moreover, for $T, T' \in \NN^X$,
\[
T \leq T'\ \overset{\operatorname{def}}{\Longleftrightarrow}\ T(x) \leq T'(x)\ \forall x \in X.
\]
In the case where $X = \{ 1, \ldots, n \}$, 
$\NN^X$ is sometimes denoted by $\NN^n$.

\subsubsection{}
\label{CT:monomial:X:T}
Let $M$ be a monoid, $X$ a finite subset of $M$ and
$T \in \NN^X$.
For simplicity,
$\sum_{x \in X} T(x) x$ is
often denoted by $T \cdot X$.
If we use the product  symbol for
the binary operation of the monoid $M$, then
$\prod_{x \in X} x^{T(x)}$ is written by
$X^T$.
Moreover, if $X = \{ X_1, \ldots, X_n \}$ and $I \in \NN^n$,
then $I \cdot X$ and $X^I$ means 
$\sum_{i=1}^n I(i) X_i$ and $\prod_{i=1}^n X_i^{I(i)}$
respectively according to a way of the binary operator
of $M$.

\subsubsection{}
\label{CT:monomial}
Let $A$ be a ring and $R$ be either the ring of polynomials of $n$-variables
over $A$, or the ring of formal power series of $n$-variables over $A$,
that is, $R = A[X_1, \ldots, X_n]$ or $A\lformal X_1, \ldots, X_n\rformal $.
For $I \in \NN^n$, we denote the monomial $X_1^{I(1)} \cdots X_n^{I(n)}$
by $X^I$.

\subsubsection{}
\label{CT:henselization:polynomial}
Let $(A, m)$ be a local ring.
The henselization of $A$ and the completion of $A$ with
respect to $m$ are denoted by
$A^h$ and $\widehat{A}$ respectively.

\renewcommand{\theTheorem}{\arabic{section}.\arabic{Theorem}}
\renewcommand{\thesubsubsection}{\arabic{section}.\arabic{subsection}.\arabic{subsubsection}}


\renewcommand{\theTheorem}{\arabic{section}.\arabic{subsection}.\arabic{Theorem}}
\renewcommand{\theClaim}{\arabic{section}.\arabic{subsection}.\arabic{Theorem}.\arabic{Claim}}
\renewcommand{\theequation}{\arabic{section}.\arabic{subsection}.\arabic{Theorem}.\arabic{Claim}}

\section{Semistable schemes over a scheme}

\subsection{Algebraic preliminaries}
\setcounter{Theorem}{0}

Here we consider several lemmas which will be used later.
Let us begin with the following lemma.

\begin{Lemma}
\label{lem:local:hom:same:residue}
Let $f : (A, m_A) \to (B, m_B)$ be a local homomorphism of
noetherian local rings
such that  $f$ induces an isomorphism 
$A/m_A \overset{\sim}{\longrightarrow} B/m_B$.
\begin{enumerate}
\renewcommand{\labelenumi}{(\arabic{enumi})}
\item
Let $x_1, \ldots, x_n$ be generator of $m_B$, that is,
$m_B = B x_1 + \cdots + B x_n$.
If $(B, m_B)$ is complete, then
\[
 B = 
 \sum_{a_1 \geq 0, \ldots, a_n \geq 0}
  f(A)x_1^{a_1} \cdots x_n^{a_n}.
\]

\item
Let $x_1, \ldots, x_n$ be elements of $m_B$ with
$m_B = Bx_1 + \cdots + Bx_n + m_A B$.
If  $(A, m_A)$ and $(B, m_B)$ are complete,
then
\[
 B = 
 \sum_{a_1 \geq 0, \ldots, a_n \geq 0}
  f(A)x_1^{a_1} \cdots x_n^{a_n}.
\]
\end{enumerate}
\end{Lemma}

\Proof
(1)
First, we claim the following:

\begin{Claim}
\label{claim:lem:local:hom:same:residue:1}
\[
m_B^d \subseteq \sum_{\substack{a_1 \geq 0, \ldots, a_n \geq 0 \\a_1 + \cdots + a_n = d}}
f(A) x_1^{a_1} \cdots x_n^{a_n} + m_B^{d+1}
\]
for all $d \geq 0$. 
\end{Claim}

We prove this claim by induction on $d$.
Since $A/m_A \simeq B/m_B$, we have $B = f(A) + m_B$, which means that
the assertion holds for $d=0$.
Thus, 
\begin{align*}
m_B & =  (f(A) + m_B) x_1 + \cdots + (f(A) + m_B) x_n  \\
& \subseteq  f(A) x_1 + \cdots + f(A) x_n + m_B^2,
\end{align*}
which show that the assertion holds for $d =1$, so that we assume $d \geq 2$.
By the hypothesis of induction,
\begin{align*}
m_B^{d} & = m_B \cdot m_B^{d-1} \\
& \subseteq  \left( f(A) x_1 + \cdots + f(A) x_n + m_B^2 \right) \cdot
 \left(\sum_{\substack{a'_1 \geq 0, \ldots, a'_n \geq 0 \\a'_1 + \cdots + a'_n = d-1}}
f(A) x_1^{a'_1} \cdots x_n^{a'_n} + m_B^{d}\right) \\
& \subseteq
\sum_{\substack{a_1 \geq 0, \ldots, a_n \geq 0 \\a_1 + \cdots + a_n = d}}
f(A) x_1^{a_1} \cdots x_n^{a_n} + m_B^{d+1}
\end{align*}
Hence, we get the claim.

\medskip
In order to complete the proof of the lemma,
it is sufficient to see the following claim:

\begin{Claim}
\label{claim:lem:local:hom:same:residue:2}
For all $b \in B$, there is a sequence 
$\{ b_d \}_{d=0}^{\infty}$ of $B$
such that
\[
b_d \in \sum_{\substack{a_1 \geq 0, \ldots, a_n \geq 0 \\a_1 + \cdots + a_n = d}}
f(A) x_1^{a_1} \cdots x_n^{a_n}
\]
and
\[
b - (b_0 + \cdots + b_d) \in m_{B}^{d+1}
\]
for all $d \geq 0$.
\end{Claim}

Since $B = f(A) + m_B$, we can set $b = b_0 + c$ with $b_0 \in f(A)$ and $c \in m_B$.
We assume that $b_0, \ldots, b_{d-1}$ are given.
Then, by  Claim~\ref{claim:lem:local:hom:same:residue:1} again,
\[
b - (b_0 + \cdots + b_{d-1}) = b_d + c',
\]
where $b_d \in \sum_{\substack{a_1 \geq 0, \ldots, a_n \geq 0 \\a_1 + \cdots + a_n = d}}
f(A) x_1^{a_1} \cdots x_n^{a_n}$ and $c' \in m_B^{d+1}$.
Therefore, we get the second claim.

\medskip
(2) Let us choose $y_1, \ldots, y_r \in A$ with $m_A = y_{1} A + \cdots + y_r A$.
Then, 
\[
m_B = x_1 B + \cdots + x_n B + f(y_{1}) B + \cdots + f(y_{r}) B.
\]
Note that
\[
x_1^{a_1} \cdots x_n^{a_n} f(y_{1})^{b_{1}} \cdots f(y_{d})^{b_r} =
f(y_{1}^{b_{1}} \cdots y_{y}^{b_r})x_1^{a_1} \cdots x_n^{a_n}.
\]
Therefore, since $(A, m_A)$ is complete,
\begin{multline*}
\sum_{\substack{a_1 \geq 0, \ldots, a_n \geq 0 \\ b_1 \geq 0, \ldots, b_r \geq 0}}
 f(A)x_1^{a_1} \cdots x_n^{a_n} f(y_{1})^{b_{1}} \cdots f(y_{d})^{b_r} \\
= \sum_{a_1 \geq 0, \ldots, a_n \geq 0} 
\left(\sum_{b_1 \geq 0, \ldots, b_r \geq 0} f(A y_{1}^{b_{1}} \cdots y_{y}^{b_r}) \right)
x_1^{a_1} \cdots x_n^{a_n} \\
= \sum_{a_1 \geq 0, \ldots, a_n \geq 0}  f(A) x_1^{a_1} \cdots x_n^{a_n}
\end{multline*}
Thus, (2) is a consequence of (1).
\QED

Next let us consider the following lemma.

\begin{Lemma}
\label{lem:flatness:x:1:x:l}
Let $(A, m)$ be a noetherian local ring and
$T \in \NN^n \setminus \{ (0, \ldots, 0)\}$.
Let $G \in m \lformal X_1, \ldots, X_n \rformal$,
$R = A \lformal X_1, \ldots, X_n \rformal/(X^T - G)$ and
$\pi : A \lformal X_1, \ldots, X_n \rformal \to R$ the canonical homomorphism.
Then, we have the following:
\begin{enumerate}
\renewcommand{\labelenumi}{(\arabic{enumi})}
\item
Let $M$ be an $A$-submodule of
$A \lformal X_1, \ldots, X_n \rformal$ given by
\[
M = \left\{ \sum_{T \not\leq I} a_I X^I \mid a_I \in A \right\}.
\]
If $(A, m)$ is complete, then $\rest{\pi}{M} : M \to R$
is bijective.

\item
$A \lformal X_1, \ldots, X_n \rformal/(X^T - G)$
is flat over $A$.
\end{enumerate}
\end{Lemma}

\Proof
(1) 
We denote $\pi(X_i)$ by $x_i$.
First, we claim the following:

\begin{Claim}
For $f \in R$, there is a sequence $\{ F_n \}_{n=0}^{\infty}$ in $M$ such that
$F_{n+1} - F_{n} \in m^n \lformal X_1, \ldots, X_n \rformal$ and
$f - \pi(F_n) \in m^n R$
for all $n$.
\end{Claim}

We will construct a sequence $\{ F_n\}_{n=0}^{\infty}$ inductively.
Clearly, we may set $F_0 = 0$.
We assume that $F_0, F_1, \ldots, F_n$ have been constructed.
Then, we can set
$f - \pi(F_n) = \pi(H) + x^{T} \pi(H')$ for some $H, H' \in 
m^n \lformal X_1, \ldots, X_n \rformal$ with
$H \in M$.
Here, $x^{T} \pi(H')= \pi(G) \pi(H') \in m^{n+1} R$.
Thus, if we set $F_{n+1} = F_n + H$, then we get our desired $F_{n+1}$.

\medskip
The above claim shows that $\rest{\pi}{M}$ is surjective.
Next, let us consider the injectivity of $\rest{\pi}{M}$.
We assume
\[
\pi \left( \sum_{T \not\leq I} a_I X^I  \right) = 0.
\]
Then, there is $H \in A\lformal X_1, \ldots, X_n \rformal$
with
\[
\sum_{T \not\leq I} a_I X^I  = (X^{T} - G) H.
\]
Here we set
\[
G = \sum_{I \in \NN^n} g_I X^I
\quad\text{and}\quad
H = \sum_{I \in \NN^n} h_I X^I.
\]
Then, $g_I \in m$ for all $I$ and
\[
\sum_{T \not\leq I} a_I X^I  =
\sum_{I \in \NN^n} h_I X^{T + I} - \sum_{I \in \NN^n} 
\left( \sum_{J + J' = I} g_J h_{J'} \right) X^I.
\]
On the left hand side of the above equation,
there is no term of a form $X^{I + T}$.
Thus,
\[
h_I = \sum_{J + J' = I + T} g_J h_{J'}
\]
for all $I$.
Here we claim that $h_I \in m^{n}$ for all $n$ and all $I$.
We see this fact by induction $n$.
First of all, since $g_I \in m$ for all $I$, we have $h_I \in m$ for all $I$.
We assume that $h_I \in m^n$ for all $I$.
Then, by the above equation, we can see that $h_I \in m^{n+1}$.
By this claim, $h_I$ must be zero for all $I$
because $\bigcap_{n \geq 0} m^n = 0$. Therefore, 
$a_I = 0$ for all $I$.

\bigskip
(2)  If $(A, m)$ is complete, then
the assertion follows from (1) by Chase's theorem \cite{Chase}.
In general, let $\widehat{A}$ be the completion of $A$ and
\[
R' = \widehat{A} \lformal X_1, \ldots, X_n \rformal/(X^T - G).
\]
Then, we have the following commutative diagram:
\[
\begin{CD}
R @>{h'}>> R' \\
@A{f}AA @AA{f'}A \\
A @>{h}>> \widehat{A}.
\end{CD}
\]
Note that $f'$, $h$ and $h'$ are faithfully flat.
Thus, so is $f$.
\QED

\begin{Remark}
\label{rem:flatness:x:1:x:l}
Let $A$ be a ring and $A[X_1, \ldots, X_n]$ the polynomial ring of
$n$-variables over $A$. Then, in the same way as in Lemma~\ref{lem:flatness:x:1:x:l},
we can see that $R = A[X_1, \ldots, X_n]/(X^T - a)$
is flat over $A$ for $T \in \NN^n \setminus \{ (0, \ldots, 0) \}$ and $a \in A$.
Indeed, if we set
\[
M = \left\{ \sum_{T \not\leq I} a_I X^I 
\in A[X_1, \ldots, X_n]  \right\},
\]
then the natural homomorphism
$\pi: M \to R$ is bijective.
The surjectivity of $\pi$ is obvious. 
We assume that an element $\sum a_I X^I$ of $M$ is zero
in $R$.
Then,
\[
\sum a_I X^I = (X^T - a) \sum_{J} b_J X^J
\]
for some $\sum_{J} b_J X^J \in A[X_1, \ldots, X_n]$.
Thus, $b_J = a b_{J + T}$ for all $J$.
Therefore,
$b_J = a^m b_{J + m T}$ for all $m > 0$.
On the other hand, since $\sum_{J} b_J X^J \in A[X_1, \ldots, X_n]$,
$b_{J + m T} = 0$ if $m \gg 0$. Hence $b_J = 0$ for all $J$.
\end{Remark}

Here we consider an approximation by an \'{e}tale neighborhood.

\begin{Proposition}
\label{prop:approximation:etale}
Let $(A, m_A)$ be a noetherian local ring essentially of finite type
over an excellent discrete valuation ring or a field.
Let $f : X \to \Spec(A)$ be a scheme of finite type over $A$.
Let $x$ be a point of $X$ such that
$f(x) = m_A$ and $A/m_A$ is isomorphic to $\OO_{X,x}/m_{X,x}$.
We assume that there are $F_1, \ldots, F_r \in A[X_1, \ldots, X_n]$
\rom{(}the polynomial ring of $n$-variables over $A$\rom{)} and
an isomorphism
\[
\phi : \widehat{A} \lformal X_1, \ldots, X_n \rformal/(F_1, \ldots, F_r)
\overset{\sim}{\longrightarrow} \widehat{\OO}_{X, x}
\]
over $\widehat{A}$ with
$\phi(\bar{X}_i) \in \widehat{m}_{X,x}$ for all $i$,
where $\bar{X}_i = X_i \mod (F_1, \ldots, F_r)$.
Then, there is an \'{e}tale neighborhood $(U, x')$ of $X$ at $x$ together
with an \'{e}tale morphism
\[
\rho : U \to \Spec(A[T_1, \ldots, T_n]/(F_1(T), \ldots, F_r(T)))
\]
such that $\rho(x') = (m_A, \bar{T}_1, \ldots, \bar{T}_n)$, where
$\bar{T}_i = T_i \mod (F_1(T), \ldots, F_r(T))$.
\end{Proposition}

\Proof
First note that
\[
F_1(\phi(\bar{X}_1), \ldots, \phi(\bar{X}_n)) = \cdots =
F_r(\phi(\bar{X}_1), \ldots, \phi(\bar{X}_n)) = 0.
\]
Thus, by Artin's approximation theorem \cite{Artin},
there are $t_1, \ldots, t_n \in \OO_{X,x}^h$
such that
\[
F_1(t_1, \ldots, t_n) = \cdots =
F_r(t_1, \ldots, t_n) = 0
\]
and $t_i - \phi(\bar{X}_i) \in \widehat{m}^2_{X,x}$ for all $i$.
Here we claim the following:

\begin{Claim}
\label{claim:prop:approximation:etale:1}
\begin{align*}
\widehat{m}_{X,x} & = \phi(\bar{X}_1) \widehat{\OO}_{X,x} +
\cdots + \phi(\bar{X}_n) \widehat{\OO}_{X,x} + m_A \widehat{\OO}_{X,x} \\
& = t_1 \widehat{\OO}_{X,x} +
\cdots + t_n \widehat{\OO}_{X,x} + m_A \widehat{\OO}_{X,x}.
\end{align*}
\end{Claim}

Clearly,
\[
\widehat{m}_{X,x} \supseteq \phi(\bar{X}_1) \widehat{\OO}_{X,x} +
\cdots + \phi(\bar{X}_n) \widehat{\OO}_{X,x} + m_A \widehat{\OO}_{X,x}.
\]
Conversely, let us pick up
$f \in \widehat{m}_{X,x}$. Then, we can write
$f = \sum_{I} a_I \bar{X}^{I}$.
If $a_{(0,\ldots, 0)} \in \widehat{A}^{\times}$, then
$f$ must be a unit because $f \in a_{(0, \ldots, 0)} + \widehat{m}_{X,x}$.
This is a contradiction. Thus, $a_{(0, \ldots, 0)} \in m_A$, which means that
\[
f \in \phi(\bar{X}_1) \widehat{\OO}_{X,x} +
\cdots + \phi(\bar{X}_n) \widehat{\OO}_{X,x} + m_A \widehat{\OO}_{X,x}.
\]
Since $t_i - \phi(\bar{X}_i) \in \widehat{m}^2_{X,x}$, we can see that
\[
\widehat{m}_{X,x} =
t_1 \widehat{\OO}_{X,x} +
\cdots + t_n \widehat{\OO}_{X,x} + m_A \widehat{\OO}_{X,x} + \widehat{m}^2_{X,x}.
\]
Hence, by Nakayama's lemma, we have our desired result.

\medskip
Let us choose an \'{e}tale neighborhood $(U,x')$ of $X$ at $x$ with
the same residue fields such that
$t_1, \ldots, t_n$ are defined over $U$. Here let us define a homomorphism
\[
\psi : A[T_1, \ldots, T_n]/(F_1(T), \ldots, F_r(T)) \to \OO_{U,x'}
\]
to be $\psi(\bar{T}_i) = t_i$ for all $i$.
By the above claim, we can see
\[
\psi^{-1}(m_{U,x'}) = (m_A, \bar{T}_1, \ldots, \bar{T}_n).
\]
Thus, it is sufficient to show that $\psi$ is \'{e}tale. Let
\[
\mu : \widehat{A} \lformal T_1, \ldots, T_n \rformal/(F_1, \ldots, F_r)
\to \widehat{A} \lformal X_1, \ldots, X_n \rformal/(F_1, \ldots, F_r)
\]
be a homomorphism given by the composition of homomorphisms
\[
 \widehat{A} \lformal T_1, \ldots, T_n \rformal/(F_1, \ldots, F_r)
 \overset{\widehat{\psi}}{\longrightarrow} \widehat{\OO}_{U,x'}
 = \widehat{\OO}_{X,x} \overset{\phi^{-1}}{\longrightarrow}
\widehat{A} \lformal X_1, \ldots, X_n \rformal/(F_1, \ldots, F_r).
\]
By the above claim and Lemma~\ref{lem:local:hom:same:residue},
$\mu$ is surjective. Hence, by the following Lemma~\ref{lem:surjective:injective:ring:hom},
it must be an isomorphism. Therefore, so is $\widehat{\psi}$.
This means that $\psi$ is \'{e}tale because
$x'$ and $(m_A, \bar{T}_1, \ldots, \bar{T}_n)$ have the same residue field.
\QED

Finally, we consider the following two lemmas concerning the bijectivity
of ring homomorphisms.

\begin{Lemma}
\label{lem:surjective:injective:ring:hom}
Let $\phi : A \to A$ be an endomorphism of a noetherian ring.
If $\phi$ is surjective, then $\phi$ is injective.
\end{Lemma}

\Proof
We set $I_n = \Ker(\phi^n)$ for $n \geq 1$.
Since $\phi$ is surjective, we can see that
$\phi(I_{n+1}) = I_{n}$ for all $n \geq 1$.
Moreover, there is $N \geq 1$ such that
$I_{N+1} = I_{N}$ because
$A$ is noetherian and $I_n \subseteq I_{n+1}$ for all $n \geq 1$.
Therefore,
\[
\Ker(\phi) = I_1 = \phi^N(I_{N+1}) = \phi^N(I_{N}) = \{ 0 \}.
\]
\QED

\begin{Lemma}
\label{lem:flat:isom}
Let
\[
\xymatrix{
(C, m_C) \ar[rr]^h & & (B, m_B) \\
& (A, m_A) \ar[ul]^g \ar[ur]_f & 
}
\]
be a commutative diagram of local homomorphisms of
noetherian local rings.
We assume that $f$ and $g$ are flat,
$h$ is surjective and that $h$ induces an isomorphism
$C/m_AC \overset{\sim}{\longrightarrow} B/m_AB$.
Then, $h$ is an isomorphism.
\end{Lemma}

\Proof
Let $I$ be the kernel of $h$.
Then, we have an exact sequence
\[
0 \to I \to C \to B \to 0.
\]
Here $B$ and $C$ are flat over $A$.
Thus, so is $I$ over $A$.
Therefore, tensoring $A/m_A$ with
the above exact sequence, it gives rise to an exact sequence
\[
 0 \to I/m_A I \to C/m_AC  \to  B/m_A B \to 0.
\]
Here $C/m_AC \to  B/m_A B$ is an isomorphism.
Thus, $I = m_A I$. We assume that
$I \not= 0$. 
Then, $m_C \in \Supp(I)$ as $C$-modules. 
Thus, $I$ is faithfully flat over $A$ by \cite[Theorem~7.3]{MatComm}.
Hence, by \cite[Theorem~7.2]{MatComm},
$I \not= m_AI$, which is a contradiction.
Therefore, $I = 0$.
\QED

\subsection{Semistable varieties over an algebraically closed field}
\setcounter{Theorem}{0}
Let $k$ be an algebraically closed field and $X$ an algebraic scheme over $k$.
A closed point $x$ of $X$ is called {\em a semistable point of $X$} if
the completion of the local
ring at $x$ is isomorphic to a ring of type
\[
 k \lformal X_1, \ldots, X_n \rformal/(X_1 \cdots X_l).
\]
The number $l$ is called {\em the multiplicity of $X$ at $x$}, and
is denoted by $\mult_{x}(X)$.
Moreover, we say $X$ is a semistable variety over $k$ if
every closed point is a semistable point.
By the following Proposition~\ref{prop:open:semistable:point},
the set of all semistable points of $X$ is a Zariski open set.
Thus, we say a point $x$ of $X$ ($x$ is not necessarily closed) is
{\em a semistable point} if there is a Zariski open set $U$ of $X$
such that $x \in U$ and every closed point of $U$ is a semistable point.
Let $\Omega$ be an algebraically closed field such that
$k$ is a subfield of $\Omega$. Note that if $X$ is a semistable variety over $k$,
then so is $X_{\Omega} = X \times_{\Spec(k)} \Spec(\Omega)$ 
over $\Omega$ (cf. Proposition~\ref{prop:semistable:alg:closed}).

\begin{Proposition}
\label{prop:open:semistable:point}
Let $X$ be an algebraic scheme over an algebraically closed field $k$.
If $x$ is a semistable closed point of $X$, then there is
a Zariski open set $U$ of $X$ such that $x \in U$ and
every closed point of $U$ is a semistable point.
\end{Proposition}

\Proof
By Proposition~\ref{prop:approximation:etale},
there are an \'{e}tale neighborhood $(U,x')$ of $x$ and
an \'{e}tale morphism 
\[
\rho : U \to \Spec(k[T_1, \ldots, T_n]/(T_1 \cdots T_l))
\]
with $\rho(x') = (0, \ldots, 0)$.
Note that $\Spec(k[T_1, \ldots, T_n]/(T_1 \cdots T_l))$
is a semistable variety over $k$. Thus, so is $U$ over $k$.
Therefore, every closed point of $\pi(U)$ is a semistable point.
\QED

\begin{Proposition}
\label{prop:semistable:alg:closed}
Let $X$ be
an algebraic scheme over an algebraically closed field $k$.
Let $\Omega$ be an algebraically closed field such that
$k$ is a subfield of $\Omega$.
Let  $\pi : X_{\Omega} = X \times_{\Spec(k)} \Spec(\Omega) \to X$
be the canonical morphism.
For a point $y \in X_{\Omega}$, if $x = \pi(y)$ is a semistable point,
then so is $y$.
\end{Proposition}

\Proof
Let $U$ be an open set of $X$ containing $x$ such that
every closed point of $U$ is a semistable point.

First we assume that $y$ is a closed point.
Let us choose a closed point $o \in \overline{\{ x \}} \cap U$.
By using Proposition~\ref{prop:approximation:etale} and shrinking $U$ 
around $o$ if necessarily,
there are \'{e}tale morphisms 
\[
f : V \to U
\quad\text{and}\quad
g : V \to W
= \Spec(k[X_1, \ldots, X_n]/(X_1 \cdots X_l))
\]
of algebraic schemes over $k$
and closed points $o' \in V$ and $o'' \in W$ 
such that
$f(o') = o$ and $g(o') = o'' = (0, \ldots, 0)$.
Since $x \in U$, $o \in \overline{\{ x\}} \cap U$ and
$f$ is faithfully flat at $o'$, we can find $x' \in V$ with
$f(x') = x$ and $o' \in \overline{\{ x' \}}$.
Here we set 
\[
\begin{cases}
U_{\Omega} = U \times_{\Spec(k)} \Spec(\Omega), \\
V_{\Omega} = V \times_{\Spec(k)} \Spec(\Omega), \\
W_{\Omega} = \Spec(\Omega[X_1, \ldots, X_n]/(X_1 \cdots X_l))
\end{cases}
\]
and the induced morphisms $V_{\Omega} \to U_{\Omega}$ and
$V_{\Omega} \to W_{\Omega}$ are denoted by $f_{\Omega}$ and
$g_{\Omega}$ respectively.
Then, $y \in U_{\Omega}$.
Let $\tilde{y} : \Spec(\Omega) \to U_{\Omega}$ be the morphism induced
by $y$. Let $\kappa(y)$, $\kappa(x)$ and $\kappa(x')$ be the residue
fields of $y$, $x$ and $x'$ respectively.
Then, there is an embedding $\iota : \kappa(x') \hookrightarrow
\Omega$ over $k$
such that the following diagram is commutative:
\[
 \xymatrix{
 & & \kappa(x') \ar@{_{(}->}[lld]_{\iota}\\
 \Omega & \kappa(y) \ar[l]^{\tilde{y}^*}_{\sim} & 
 \kappa(x) \ar@{_{(}->}[l] \ar@{_{(}->}[u]
}
\]
This gives rise to a morphism
$\beta : \Spec(\Omega) \to V_{\Omega}$ such that
the diagram
\[
 \xymatrix{ & V_{\Omega} \ar[r]^{\pi'} \ar[d]^{f_{\Omega}} & V \ar[d]^{f} \\
   \Spec(\Omega)\ar[ru]^{\beta} \ar[r]_{\tilde{y}} & 
   U_{\Omega} \ar[r]^{\pi} & U 
}
\]
is commutative and the image of $\pi' \circ \beta$ is $x'$.
Let $y'$ be the image of $\beta$. Then, $f_{\Omega}(y') = y$.
Note that $f_{\Omega}$ and $g_{\Omega}$ are \'{e}tale and
the residue fields of $y$, $y'$ and $y'' = g_{\Omega}(y')$ are $\Omega$.
Thus, we can see that
\[
\widehat{\OO}_{X_{\Omega}, y} \simeq
\widehat{\OO}_{V_{\Omega}, y'} \simeq
\widehat{\OO}_{W_{\Omega}, y''}.
\]
We set $y'' = (a_1, \ldots, a_n) \in \AAA^n(\Omega)$ and
$I = \{ i \mid \text{$a_i = 0$ and $i = 1, \ldots, l$} \}$.
Note that $I \not= \emptyset$.
Therefore, if we set $Z_i = X_i - a_i$ and $Z = \prod_{i \in I} Z_i$, then
it is easy to see that
\[
\widehat{\OO}_{W_{\Omega}, y''}
= \Omega\lformal Z_1, \ldots, Z_n \rformal/(Z).
\]
Thus, we get our lemma in the case where $y$ is a closed point.

Next we consider a general case.
We set $U_{\Omega} = \pi^{-1}(U)$. Then, by the previous observation,
every closed point of $U_{\Omega}$ is a semistable point.
On the other hand, $y \in U_{\Omega}$. Thus, $y$ is a semistable point.
\QED

\subsection{Semistable schemes}
\setcounter{Theorem}{0}

Let $S$ be a locally noetherian scheme and $f : X \to S$ a morphism of finite type.
First we assume that $S = \Spec(k)$ for some field $k$.
Let $\bar{k}$ be the algebraic closure of $k$,
$X' = X \times_{\Spec(k)} \Spec(\bar{k})$, and $\pi : X' \to X$ the
canonical morphism.
A point $x$ of $X$ is called {\em a semistable point of $X$} if
every point $x'$ of $X'$ with $\pi(x') = x$ is a semistable point.

For a general $S$, we say  $f : X \to S$ is {\em semistable at $x \in X$}
if $f$ is flat at $x$ and
$x$ is a semistable point of the fiber $f^{-1}(f(x))$ passing through $x$.
Moreover, we say $X$ is a semistable scheme over $S$
if $f$ is semistable at all points of $X$.
By Proposition~\ref{prop:semistable:alg:closed}, for a flat morphism
$f : X \to S$, $X$ is a semistable scheme over $S$
if and only if, for any algebraically closed field $\Omega$, any
morphism $\Spec(\Omega) \to S$ and
any closed point $x' \in X \times_{S} \Spec(\Omega)$,
the completion of the local
ring at $x'$ is isomorphic to a ring of type
\[
 \Omega \lformal X_1, \ldots, X_n \rformal/(X_1 \cdots X_l).
\]
We say a semistable scheme $X$ over $S$ is
proper if $X$ is proper over $S$. Moreover, a proper semistable
scheme $X$ over $S$ is said to be connected if
$f_* (\OO_X) = \OO_S$.
The following proposition is a local description of semistable scheme.

\begin{Proposition}
\label{prop:local:description:semistable:scheme}
Let $f : (A, m_A) \to (B, m_B)$ be a flat local homomorphism of
noetherian local rings such that $f$ is essentially of finite type and
$f$ induces an isomorphism 
$A/m_A \overset{\sim}{\longrightarrow} B/m_B$.
We assume that  there is an isomorphism
\[
(A/m_A) \lformal T_1, \ldots, T_n \rformal/(T_1 \cdots T_l)
 \overset{\sim}{\longrightarrow}
\widehat{B}/m_A \widehat{B}
\]
over $A/m_A$. Then, there are a local homomorphism 
$g : (A, m_A) \to (C, m_C)$ of noetherian local rings and
$F \in C$ with the following properties:
\begin{enumerate}
\renewcommand{\labelenumi}{(\arabic{enumi})}
\item
$g$ is essentially of finite type  and smooth over $A$, and
$g$ induces an isomorphism
$A/m_A \overset{\sim}{\longrightarrow} C/m_{C}$.

\item
There are $x_1, \ldots, x_n \in C$ such that
$dx_1, \ldots, dx_n$ form a free basis of $\Omega_{C/A}$ and
$F - x_1 \cdots x_l \in m_A C$.

\item
There is a local  \'{e}tale homomorphism
$\rho : B \to C/F C$ over $A$.
\end{enumerate}
In particular,
we have an isomorphism
\[
A \lformal X_1, \ldots, X_n \rformal/(X_1 \cdots X_l - G) \overset{\sim}{\longrightarrow}  B
\]
over $A$ for some $G \in m_A\lformal X_1, \ldots, X_n \rformal$.
\end{Proposition}

\Proof
If $l=1$, then $B$ is smooth over $A$.
Thus, if we set $C$ as the localization of $B[X]$ at the maximal ideal
$(X, m_B)$, then we get our desired result, so that
we assume that $l \geq 2$.

We set $k = A/m_A$.
Since $B$ is essentially of finite type over $A$,
there are a polynomial ring $A[ Z_1, \ldots, Z_N]$, a maximal ideal $M$ of $A[ Z_1, \ldots, Z_N]$
and an ideal $I$ of $k[ Z_1, \ldots, Z_N]_M$ such that
$A[Z_1, \ldots, Z_N]/M = k$ and $B = A[Z_1, \ldots, Z_N]_M/I$.
By abuse of notation, we denote $MA[Z_1, \ldots, Z_N]_M$ by $M$.

Let us begin with the following claim:

\begin{Claim}
\label{claim:prop:local:description:semistable:scheme:1}
There are a local $A$-algebra $(D, m_D)$,  a surjective local homomorphism
$\phi : D \to B$ and
an element $F_0$ of $D_0 = D/m_AD$ such that 
$D$ is essentially of finite type and smooth over $A$, 
and that the kernel of $D_0 \to B_0 = B/m_A B$ is $F_0 D_0$.
\end{Claim}

We set 
\[
\begin{cases}
E = A[Z_1, \ldots, Z_N]_M, \\
E_0 = A[Z_1, \ldots, Z_N]_M/ m_A A[Z_1, \ldots, Z_N]_M = k[Z_1, \ldots, Z_N]_{M_0},
\end{cases}
\]
where $M_0$ is the maximal ideal of $k[Z_1, \ldots, Z_N]$ corresponding to
$M$. As before, we also denote $M_0 k[Z_1, \ldots, Z_N]_{M_0}$ by $M_0$.
Let $p : E \to E_0$ be the canonical surjective homomorphism.
Let $I_0$ be the ideals of $E_0$ given by $I_0 = p(I)$. 
Then, $E_0/M_0 = k$ and $B_0 = E_0/I_0$.
Let $u_1, \ldots, u_N$ be a system of parameters of $E_0$.
Since $\widehat{B}_0 \simeq k \lformal T_1, \ldots, T_n \rformal/(T_1 \cdots T_l)$,
$M_0/(M_0^2 + I_0)$ has the dimension $n$ over $k$.
Therefore, after renumbering
$u_1, \ldots, u_N$, we may assume that $u_1, \ldots, u_n$ gives rise to a basis
of $M_0/(M_0^2 + I_0)$.
Thus, for $i > n$, there are $a_{i1}, \ldots, a_{in} \in E_0$
such that
\[
u_i - (a_{i1}u_1 + \cdots + a_{in}u_n) \in M_0^2 + I_0.
\]
Therefore, we can find $g_i \in M_0^2$ and $h_i \in I_0$
with $h_i = u_i - (a_{i1}u_1 + \cdots + a_{in}u_n) + g_i$.
Let us choose that $U_1, \ldots, U_n \in E$ and
$H_{n+1}, \ldots, H_N \in I$ such that $p(U_i) = u_i$ and
$p(H_j) = h_j$ ($i=1, \ldots, n$ and $j = n+1, \ldots, N$).
Here we consider a homomorphism
\[
\rho : E^N \to \Omega_{E/A}
\]
given by
\[
\rho(e_i) = \begin{cases}
d U_i & \text{if $1 \leq i \leq n$} \\
d H_i & \text{if $n+1 \leq i \leq N$},
\end{cases}
\]
where $\{ e_1, \ldots, e_N \}$ is the standard basis of $E^N$.
Since 
\[
d u_1, \ldots, d u_n, d h_{n+1}, \ldots, d h_N
\]
form
a basis of $\Omega_{E_0/k} =  \Omega_{E/A} \otimes_{A} k$,
by using Nakayama's lemma,
we can see that $\rho$ is surjective.
Hence $\rho$ is an isomorphism because $\Omega_{E/A}$
is a free $E$-module of rank $N$.
Therefore, 
\[
d U_1, \ldots, d U_n, d H_{n+1}, \ldots, d H_N
\]
form a free basis of $\Omega_{E/A}$.
Here we set 
\[
D = E/(H_{n+1}, \ldots, H_N)
\quad\text{and}\quad
D_0 = D/m_AD = E_0/(h_{n+1}, \ldots, h_n)
\]
Then, by the previous observation,
$D$ is smooth over $A$ (cf. \cite[Chaper~VII, Theorem~5.7]{AK}), 
and $D_0$ is regular and $\dim D_0 = n$.
On the other hand, since $\dim B_0 = n-1$
and $D_0$ is UFD, there is $F_0 \in D_0$ with
$B_0 = D_0/F_0 D_0$.

\medskip
Let $m_{D_0}$ and $m_{B_0}$ be the maximal ideals of $D_0$ and $B_0$ respectively.
Let $\phi_0 : D_0 \to B_0$ be the natural homomorphism, and let
$\phi_0^h : D_0^h \to B_0^h$ and
$\widehat{\phi}_0 : \widehat{D}_0 \to \widehat{B}_0$ be the induced homomorphisms.

\begin{Claim}
There are $t_1, \ldots, t_n \in D_0^h$ such that
$\phi_0^h(t_1 \cdots t_l) = 0$ and $t_1, \ldots, t_n$ are
a system of parameters of $D_0$.
\end{Claim}

Let us consider the factorization of $F_0$ in $\widehat{D}_0$.
Note that
\[
\widehat{D}_0/F_0 \widehat{D}_0 = \widehat{B}_0 \simeq
 k \lformal T_1, \ldots, T_n \rformal/(T_1 \cdots T_l).
\]
Thus,  we can see that
there is a system of parameters $t'_1, \ldots, t'_n$ of $\widehat{D}_0$
with $F_0 = u t'_1 \cdots t'_l$ for some 
$u \in \widehat{D}_0^{\times}$.
In particular, $F_0 \in m_{D_0}^l \subseteq m_{D_0}^2$.
Note that $B_0^h = D_0^h/F_0 D_0^h$ and
$\widehat{B}_0 = \widehat{D}_0/F_0\widehat{D}_0$
(cf.  [EGA~IV, 18.6.8]).
Since $\widehat{\phi}_0(t'_1) \cdots \widehat{\phi}_0(t'_l) = 0$,
using Artin's approximation theorem \cite{Artin} for $\widehat{B}_0$ and $B_0^h$,
we can find $\bar{t}_1, \ldots, \bar{t}_n \in B_0^h$ such that
$\bar{t}_1 \cdots \bar{t}_l = 0$ and
$\widehat{\phi}_0(t'_i) - \bar{t}_i \in m_{B_0}^2 \widehat{B}_0$ for
$i = 1, \ldots, n$.
Let us choose $t_1, \ldots, t_n$ with
$\phi_0^h(t_i) = \bar{t}_i$ for $i=1, \ldots, n$. Then,
\[
t'_i - t_i \in m^2_{D_0}\widehat{D}_0 + F_0\widehat{D}_0 \subseteq m^2_{D_0}\widehat{D}_0
\]
for all $i$.
Therefore, $t_1, \ldots, t_n$ form a system of parameters of $D_0^h$ and
$\phi^h_0(t_1 \cdots t_l) = \bar{t}_1 \cdots \bar{t}_l = 0$.

\medskip
Let $\pi : B \to B_0$, $\pi' : D \to D_0$ and
$\phi : D \to B$ be the natural
surjective homomorphisms.
These induce the following commutative diagram:
\[
\begin{CD}
D^h @>{\phi^h}>>  B^h \\
@V{{\pi'}^h}VV @VV{\pi^h}V \\
D_0^h @>{\phi_0^h}>>  B_0^h.
\end{CD} 
\]
Let us choose $x_1, \ldots, x_n \in D^h$ with
${\pi'}^h(x_i) = t_i$ for all $i = 1, \ldots, n$.
Since
\[
 \pi^h(\phi^h(x_1 \cdots x_l)) = 0,
\]
$\phi^h(x_1 \cdots x_l) \in m_A B^h$.
Thus, there is
$G \in m_A D^h$ such that
\[
 \phi^h(x_1 \cdots x_1 - G) = 0.
\]

\begin{Claim}
A surjective homomorphism
$D^h/(x_1 \cdots x_l - G)D^h \to B^h$ is an isomorphism.
\end{Claim}

For this purpose, it is sufficient to see that
$\widehat{D}/(x_1 \cdots x_l - G)\widehat{D} \to \widehat{B}$ is an isomorphism.
By using Lemma~\ref{lem:local:hom:same:residue},
there is a surjective homomorphism
\[
\mu : 
\widehat{A} \lformal X_1, \ldots, X_n \rformal \to \widehat{D}
\] 
defined by $\mu(X_i) = x_i$ for $i=1, \ldots, n$.
Let us choose
\[
G' \in \widehat{m}_A  \lformal X_1, \ldots, X_n \rformal
\]
with $\mu(G') = G$.
Then, by  Lemma~\ref{lem:flatness:x:1:x:l}, Lemma~\ref{lem:surjective:injective:ring:hom} and
Lemma~\ref{lem:flat:isom}, a homomorphism
\[
\widehat{A} \lformal X_1, \ldots, X_n \rformal /(X_1 \cdots X_l - G') 
\to \widehat{B}
\]
given by the composition of homomorphisms
\[
\widehat{A} \lformal X_1, \ldots, X_n \rformal /(X_1 \cdots X_l - G') \to
\widehat{D}/(x_1 \cdots x_l - G)\widehat{D} 
\to \widehat{B}
\]
is an isomorphism.
Therefore,
so is $\widehat{D}/(x_1 \cdots x_l - G)\widehat{D} \to \widehat{B}$.

\medskip
By the above claim, $\Ker(\phi)D^h = (x_1 \cdots x_l - G)D^h$.
We choose a local \'{e}tale neighborhood 
$(D, m_D) \to (C, m_C)$ of $\Spec(D)$ such that
$D/m_D \simeq C/m_C$, $T_1, \ldots, T_n, G \in C$ and
$\Ker(\phi)C = (x_1 \cdots x_l - G)C$.
Then, 
\[
B = D \otimes_D (D/\Ker(\phi)) \to  C \otimes_D (D/\Ker(\phi)) = C /(x_1 \cdots x_l - G)C
\]
is \'{e}tale.
Moreover, we can see that $dx_1, \ldots, dx_n$ form
a free basis of $\Omega_{C/A}$ in the same way
as before because $dt_1, \ldots, dt_n$ are
a basis of $\Omega_{C_0/k}$, where $C_0 = C/m_AC$.
\QED

\renewcommand{\theTheorem}{\arabic{section}.\arabic{Theorem}}
\renewcommand{\theClaim}{\arabic{section}.\arabic{Theorem}.\arabic{Claim}}
\renewcommand{\theequation}{\arabic{section}.\arabic{Theorem}.\arabic{Claim}}
\renewcommand{\theTheorem}{\arabic{section}.\arabic{subsection}.\arabic{Theorem}}
\renewcommand{\theClaim}{\arabic{section}.\arabic{subsection}.\arabic{Theorem}.\arabic{Claim}}
\renewcommand{\theequation}{\arabic{section}.\arabic{subsection}.\arabic{Theorem}.\arabic{Claim}}
\section{Preliminaries to log schemes}

\subsection{Pushout in the category of monoids}

In the following, every monoid is a commutative monoid with the unity.
Let $f : Q \to P$ and $g : Q \to R$ be homomorphisms of monoids.
Then, $P \pushout_Q R$ is defined as follows:
We define a relation $\sim$ on $P \times R$ by
\[
 (p, r) \sim (p', r')
 \ \Longleftrightarrow\ 
 \text{$(p, r)\cdot (f(q), g(q')) = (p', r') \cdot (f(q'), g(q))$ 
 for some $q, q' \in Q$}.
\] 
Then, $\sim$ is an equivalence relation on $P \times R$.
Indeed, if $(p, r) \sim (p', r')$ and $(p', r') \sim (p'', r'')$, then
there are $q_1, q_2, q_3, q_4 \in Q$ with
\[
\begin{cases}
(p, r) \cdot (f(q_1), g(q_2)) = (p', r') \cdot (f(q_2), g(q_1)), \\
(p', r') \cdot (f(q_3), g(q_4)) = (p'', r'') \cdot (f(q_4), g(q_3)).
\end{cases}
\]
Thus,
\[
(p,r) \cdot (f(q_1q_3), g(q_2q_4)) = (p'', r'') \cdot (f(q_2q_4), g(q_1q_3)).
\]
Hence we set
\[
 P \pushout_Q R = P \times R/\!\!\sim.
\]
The following are elementary properties of $P \pushout_Q R$.

\begin{Lemma}
\label{lem:pushout:monoid}
\begin{enumerate}
\renewcommand{\labelenumi}{(\arabic{enumi})}
\item
$P \pushout_Q R$ is a monoid in the natural way and we have a commutative
diagram:
\[
\xymatrix{
Q \ar[r]^{g} \ar[d]_f &  R \ar[d]^{\beta} \\
P \ar[r]^<<<<<{\alpha} & P \pushout_Q R.
}
\]

\item
Let
\[
\xymatrix{
Q \ar[r]^{g} \ar[d]_f &  R \ar[d]^{\beta'} \\
P \ar[r]^{\alpha'} & M
}
\]
be a commutative diagram of monoids.
We set $\delta = \alpha' \circ f = \beta' \circ g$.
We assume that $m \cdot \delta(q) = m' \cdot \delta(q)$ implies
$m = m'$ for any $m, m' \in M$ and $q \in Q$.
Then, there is the unique homomorphism
$\gamma : P \pushout_Q R \to M$
such that the following diagram is commutative:
\[
\xymatrix{
Q  \ar[r]^{g} \ar[d]_{f} & R \ar[d]^{\beta} \ar[rdd]^{\beta'} & \\
P \ar[r]^<<<<<{\alpha} \ar[rrd]_{\alpha'} & P \pushout_Q R \ar[rd]^{\gamma} & \\
 & & M
}
\]

\item
If one of $P, Q, R$ is a group, then $P \pushout_Q R$ is the pushout
of $(f : Q \to P, g : Q \to R)$ in the category of monoids.

\item
If $P$ and $R$ are integral,
then so is $P \pushout_Q R$.
Moreover, $P \pushout_Q R$ has the following universality:
In the commutative diagram in \rom{(2)}, if $M$ is integral,
then there is the unique homomorphism
$\gamma : P \pushout_Q R \to M$ with
$\alpha' = \gamma \circ \alpha$ and $\beta' = \gamma \circ \beta$.
Namely, $P \pushout_Q R$ is the pushout of
$(f : Q \to P, g : Q \to R)$ in the category of integral monoids.

\item
Let $h : R \to S$ be another homomorphism of monoids.
We assume that an equation $s \cdot h(r) = s' \cdot h(r)$ for
$s, s' \in S$ and $r \in R$ implies $s = s'$.
Then, there is a natural isomorphism
\[
P \pushout_{Q} S \overset{\sim}{\longrightarrow} (P \pushout_{Q} R) \pushout_{R} S.
\] 
\end{enumerate}
\end{Lemma}

\Proof
(1) We denote the class of $(p,r)$ by $[p,r]$.
It is easy to see that if $(p,r) \sim (p', r')$ and
$(s, t) \sim (s', t')$, then $(p,r)\cdot(s,t) \sim (p',r')\cdot(s', t')$.
Thus, the product of $P \pushout_Q R$ is given by
$[p,r] \cdot [p', r'] = [pp', rr']$. Moreover, $\alpha(p) = [p, 1]$ and
$\beta(r) = [1, r]$.

\medskip
(2) The uniqueness of $\gamma$ is obvious.
We would like to define $\gamma$ to be
\[
\gamma([p,r]) = \alpha'(p) \cdot \beta'(r).
\]
For this purpose, we assume that $[p, r] = [p', r']$, that is,
there are $q, q' \in Q$ with
\[
(p, r) \cdot(f(q), g(q')) = (p', r') \cdot (f(q'), g(q)).
\]
Thus,
\[
(\alpha'(p) \cdot \delta(q)) \cdot (\beta'(r) \cdot \delta(q')) 
=(\alpha'(p') \cdot \delta(q')) \cdot (\beta'(r') \cdot \delta(q)).
\]
Therefore, $\alpha'(p) \cdot \beta'(r) = \alpha'(p') \cdot \beta'(r')$.

\medskip
(3)
We need to show that $m \cdot \delta(q) = m' \cdot \delta(q)$ implies 
$m = m'$.
If $Q$ is a group, then there is $q^{-1} \in Q$. Thus, $m = m'$.
Next we assume that $P$ is a group.
Then, there is $f(q)^{-1} \in P$.
Thus, $m \cdot \alpha'(f(q)) = m' \cdot \alpha'(f(q))$ implies $m = m'$.
If $R$ is a group, then we have our assertion in the same way
as in the case where $P$ is a group.

\medskip
(4)
We assume $[p,r]\cdot[s,t] = [p',r']\cdot[s,t]$.
Then, there are $q, q' \in Q$ with
\[
(ps, rt) \cdot(f(q), g(q')) = (p's, r't) \cdot (f(q'), g(q)).
\]
Here $P$ and $R$ are integral. Thus,
\[
(p, r) \cdot(f(q), g(q')) = (p', r') \cdot (f(q'), g(q)).
\]
Thus, $[p,r] = [p', r']$.
The later part is a consequence of (2)

\medskip
(5) Here we define
\[
\phi : P \pushout_{Q} S \to (P \pushout_{Q} R) \pushout_{R} S
\]
to be $\phi(p,s) = [[p, 1], s]$.
This is well-defined. For, if $(p,s) \sim (p',s')$, then
there are $q_1, q_2 \in Q$ with $(p,s) \cdot (f(q_1), h(g(q_2)) = (p',s') \cdot (f(q_2), h(g(q_1))$.
Thus,
\begin{align*}
([p,1], s) \cdot (\beta(g(q_1)), h(g(q_2))) & = ([p,1], s) \cdot (\alpha(f(q_1)), h(g(q_2))) \\
& = ([p f(q_1), 1], sh(g(q_2))) \\
& = ([p' f(q_2), 1], s'h(g(q_1))) \\
& = ([p', 1], s') \cdot (\alpha(f(q_2)),  h(g(q_1))) \\
& = ([p', 1], s') \cdot (\beta(g(q_2)), h(g(q_1))).
\end{align*}
Moreover, since $([p, r], s) \sim ([p,1], h(r) s)$, $\phi$ is surjective.
Finally, we assume that $\phi(p,s) = \phi(p',s')$.
Then, there are $r_1, r_2 \in R$ such that
$(p, r_1) \sim (p', r_2)$ and $s h(r_2) = s' h(r_1)$.
Thus, there are $q_1, q_2$ with
$p f(q_1) = p' f(q_2)$ and $r_1 g(q_2) = r_2 g(q_1)$.
Hence
\[
s h(g(q_2)) h(r_1) = s h(g(q_1)) h(r_2) = s' h(g(q_1)) h(r_1),
\]
which implies that 
$s h(g(q_2)) = s' h(g(q_1))$ by our assumption.
Therefore, $[p, s] = [p', s']$. 
\QED

\subsection{Semistable structure of monoids}
First, let us recall a semistable structure of a monoid.
For details, see \cite{IwaMw}.
Here, the binary operators of monoids are written in 
the additive way.
Let $f : Q \to P$ be an integral homomorphism of fine and sharp monoids.
Let $\sigma$ be a finite subset of $P$.
For $T \in \NN^{\sigma} = \Map(\sigma, \NN)$, we denote
an element $\sum_{x \in \sigma} T(x) x$ of $P$ by $T \cdot \sigma$.
For $q_0 \in Q$ and $\Delta, B \in \NN^{\sigma}$, we say
$P$ has a semistable structure $(\sigma, q_0, \Delta, B)$ over $Q$
(or $P$ is of semistable type $(\sigma, q_0, \Delta, B)$ over $Q$)
if the following conditions
are satisfied:
\begin{enumerate}
\renewcommand{\labelenumi}{(S\arabic{enumi})}
\item
$q_0 \not= 0$, $\Supp(\Delta) \not= \emptyset$ and 
$\Delta(x)$ is either $0$ or $1$ for all $x \in \sigma$.

\item
$P$ is generated by $\sigma$ and $f(Q)$ and 
the natural homomorphism $\NN^{\sigma} \to P$ given by
$T \mapsto T \cdot \sigma$ is injective.

\item
$\Supp(\Delta) \cap \Supp(B) = \emptyset$ and
$\Delta \cdot \sigma = f(q_0) + B \cdot \sigma$.

\item
If we have a relation
\[
T \cdot \sigma = f(q) + T' \cdot \sigma
\quad(T, T' \in \NN^{\sigma})
\]
with $q \not= 0$,  then $T(x) > 0$ for all $x \in \Supp(\Delta)$.
\end{enumerate}
The following are basic properties of a semistable structure
(cf. \cite[Proposition~2.2]{IwaMw}):
\begin{enumerate}
\renewcommand{\labelenumi}{(\Alph{enumi})}
\item
Let $T, T' \in \NN^{\sigma}$ and $q \in Q$.
If  $T \cdot \sigma = f(q) + T' \cdot \sigma$,
$q \not= 0$ and $\Supp(T) \cap \Supp(T') = \emptyset$,
then $q = n q_0$, $T = n \Delta$ and
$T' = nB$ for some $n \in \NN$.

\item
Let $T, T' \in \NN^{\sigma}$ and $q, q' \in Q$.
If $T \cdot \sigma + f(q) = T' \cdot \sigma + f(q')$ and
there are $x, x' \in \Supp(\Delta)$ with
$T(x) = T'(x') = 0$, then
$T = T'$ and $q = q'$.

\item
Let $\NN \to Q \times \NN^{\sigma \setminus \Supp(\Delta)}$ and
$\NN \to \NN^{\Supp(\Delta)}$ be homomorphisms
given by $1 \mapsto (f(q_0), \rest{B}{\sigma \setminus \Supp(\Delta)})$ and
$1 \mapsto \rest{\Delta}{\Supp(\Delta)}$ respectively.
Then, 
\[
P \simeq (Q \times \NN^{\sigma \setminus \Supp(\Delta)} ) \pushout_{\NN}
\NN^{\Supp(\Delta)}.
\]
\end{enumerate}

Let us begin with the following lemma.

\begin{Lemma}
\label{lem:no:splitable:equiv:semistbale}
Let $f : Q \to P$ be an integral homomorphism of fine and sharp monoids.
We assume that $P$ has a semistable structure
$(\sigma, q_0, \Delta, B)$
over $Q$ for some $\sigma \subseteq P$,
$q_0 \in Q$ and $\Delta, B \in \NN^{\sigma}$. Then, the following are equivalent.
\begin{enumerate}
\renewcommand{\labelenumi}{(\arabic{enumi})}
\item
$f : Q \to P$ does not split, that is, there is no sub-monoid $N$ of $P$ with
$P = f(Q) \times N$.

\item $\#(\Supp(\Delta)) \geq 2$.

\item
For all $x \in \sigma$, $x$ is irreducible and $x \not \in f(Q)$.
\end{enumerate}
\end{Lemma}

\Proof
(1) $\Longrightarrow$ (2):
We assume $\#(\Supp(\Delta)) = 1$. We set $\Supp(\Delta) = \{ p \}$
and $\sigma' = \sigma \setminus \{ p \}$.
Then, since $p = f(q_0) + B \cdot \sigma$,
$P$ is generated by $\sigma'$ and $f(Q)$.
Let $N$ be a sub-monoid of $P$ generated by $\sigma'$.
Let us see that $P = f(Q) \times N$. For this purpose, it is sufficient to see that
if $I \cdot \sigma + f(q) = I' \cdot \sigma + f(q')$ with $I(p) = I'(p) = 0$, then
$I = I'$ and $q = q'$. This is nothing more than the above property (B).

\medskip
(2) $\Longrightarrow$ (3):
Let $x$ be an element of $\sigma$.
First of all, let us see that if $x = I \cdot \sigma + f(q)$ for some $I \in \NN^{\sigma}$ and $q \in Q$, then
$I = I_x$ and $q = 0$, where $I_x \in \NN^{\sigma}$ is given by
\[
I_x(y) = \begin{cases} 1 & \text{if $y = x$,} \\ 0 & \text{otherwise.} \end{cases}
\]
Indeed, if $q = 0$, then $I = I_x$. Thus, we assume that $q \not= 0$. 
If $x \not\in \Supp(I)$, then, by the property (A),
$I_x = n \Delta$ for some $n \in \NN$. This is impossible because
$\#(\Supp(\Delta))  \geq 2$. Thus, $x \in \Supp(I)$.
Therefore, we have $(I - I_x) \cdot \sigma + f(q) = 0$.
Here $P$ is sharp. Hence $I = I_x$ and $q = 0$.
This is a contradiction.

We assume that $x \in f(Q)$. Then, there is $q \in Q$ with $x = f(q)$.
Thus, by the above claim, we have a contradiction.

Next we assume that $x = y + z$ for some $y, z \in P$.
We set $y = I \cdot \sigma + f(q)$ and $z = I' \cdot \sigma + f(q')$.
Then, $x = (I + I') \cdot \sigma + f(q + q')$. Thus, $I_x = I + I'$ and
$q + q' = 0$. Therefore, we can see either $y = 0$ or $z = 0$.

\medskip
(3) $\Longrightarrow$ (1):
We assume that  there is a sub-monoid $N$ of $P$ with
$P = f(Q) \times N$.
First, let us see that $\sigma \subseteq N$. For $x \in \sigma$,
we can set $x = y + z$ for some $y \in f(Q)$ and $z \in N$.
Since $x$ is irreducible, either $y = 0$ or $z = 0$.
If $z = 0$, then $x \in f(Q)$. Thus, $y = 0$, that is, $x \in N$.
Let us consider the relation
\[
\Delta \cdot \sigma = f(q_0) + B \cdot \sigma.
\]
Then, since $\sigma \subseteq N$,
we have $\Delta \cdot \sigma = B \cdot \sigma$ and $q_0 = 0$.
This is a contradiction.
\QED

Let us consider the uniqueness of a semistable structure.

\begin{Proposition}
\label{prop:unique:semistable:non:split}
Let $f : Q \to P$ be an integral homomorphism of fine and sharp monoids.
Then, we have the following:
\begin{enumerate}
\renewcommand{\labelenumi}{(\arabic{enumi})}
\item
We assume that $P$ has semistable structures
\[
(\sigma, q_0, \Delta, B)
\quad\text{and}\quad
(\sigma', q'_0, \Delta', B')
\]
over $Q$ for some $\sigma, \sigma' \subseteq P$,
$q_0, q'_0 \in Q$, $\Delta, B \in \NN^{\sigma}$ and $\Delta', B' \in \NN^{\sigma'}$.
If $\sigma = \sigma'$, then
$q_0 = q'_0$, $\Delta = \Delta'$ and $B = B'$.

\item
We assume that $f : Q \to P$ does not split.
Then, there is the unique semistable structure of $P$ over $Q$ if it exists.
\end{enumerate}
\end{Proposition}

\Proof
(1) Note that $\Supp(\Delta') \cap \Supp(B') = \emptyset$.
Thus, by the property (B),
\[
\Delta' \cdot \sigma = f(q'_0) + B' \cdot \sigma
\]
implies that
$\Delta' = n \Delta$, $B' = nB$ and $q'_0 = nq_0$ for some $n \in \NN$.
Considering $\Delta' = n \Delta$, we can see $n=1$.

(2)
By Lemma~\ref{lem:no:splitable:equiv:semistbale}, 
for all $x \in \sigma$, 
$x$ is irreducible and $x \not\in f(Q)$.
Moreover, since $P$ is generated by $\sigma$ and $f(Q)$, 
we can see that $\sigma$ is the set of all irreducible elements of $P$ not lying in $f(Q)$.
Therefore, (1) implies (2).
\QED

The above proposition gives rise to the following definition.

\begin{Definition}
\label{def:marking:semistable:structure}
Let $f : Q \to P$ be an integral homomorphism of fine and sharp monoids such that
$f : Q \to P$ does not split.
We assume that $P$ has the semistable structure $(\sigma, q_0, \Delta, B)$
over $Q$. Then, $q_0$ is called {\em the marking of $f : Q \to P$}, and
is denoted by $\marking(P/Q)$.
\end{Definition}

Here we consider the invariance of the marking of the semistable structure by
a generalization.

\begin{Proposition}
\label{prop:quotient:semistable}
Let $f : Q \to P$ be a homomorphism of fine and sharp monoids.
Let $N$ be a sub-monoid of $P$ and $\bar{P} = P/N$ the quotient monoid of $P$ by $N$,
i.e., $\bar{P}$ is the quotient set given by the following equivalence relation $\sim$:
\[
\text{$p \sim p'$ $\Longleftrightarrow$ $p + n = p'+ n'$ \rom{(}$\exists n, n' \in N$\rom{)}.}
\]
Let $\pi : P \to \bar{P}$ be the natural homomorphism and
$\bar{f} = \pi \circ f : Q \to \bar{P}$ the composition of maps.
We assume that 
\rom{(i)} $\bar{P}$ is fine and sharp,
\rom{(ii)} $f$ and $\bar{f}$ are integral and that
\rom{(iii)} $P$ has a semistable structure $(\sigma, q_0, \Delta, B)$ over $Q$ for some
$\sigma \subseteq P$, $q_0 \in Q$ and
$\Delta, B \in \NN^{\sigma}$.
If we set $\sigma' = \{ x \in \sigma \mid \pi(x) \not= 0 \}$ and
$\bar{\sigma} = \pi(\sigma')$, then the natural map $\alpha : \sigma' \to \bar{\sigma}$
induced by $\pi$
is bijective. Moreover, if we set $\bar{\Delta} = \Delta \circ \alpha^{-1}$ and
$\bar{B} = B \circ \alpha^{-1}$ as elements of $\NN^{\bar{\sigma}}$, then
$(\bar{\sigma}, q_0, \bar{\Delta}, \bar{B})$ gives rise to
a semistable structure of $\bar{P}$ over $Q$.
\end{Proposition}

\Proof
First of all, note that $N \subseteq \pi^{-1}(0)$.
We set $\sigma'' = \{ x \in \sigma \mid \pi(x) = 0 \}$.

\begin{Claim}
\label{claim:prop:quotient:semistable:1}
$\pi^{-1}(0)$ is generated by $\sigma''$.
\end{Claim}

Let $N'$ be the sub-monoid generated by $\sigma''$. 
Clearly, $N' \subseteq \pi^{-1}(0)$.
Conversely, 
let $a$ be an element of $\pi^{-1}(0)$. We can write $a$ as a form
$T \cdot \sigma + f(q)$,
where $T \in \NN^{\sigma}$ and $q \in Q$.
Then, $\sum_{x \in \sigma} T(x) \pi(x) + \bar{f}(q) = 0$.
Thus, since $\bar{P}$ is sharp, we have
$T(x) \pi(x) = 0$ for all $x \in \sigma$ and $\bar{f}(q) = 0$.
Therefore, $q = 0$ and
$x \in \sigma''$ if $T(x) > 0$. This means that $c \in N'$.

\begin{Claim}
\label{claim:prop:quotient:semistable:2}
If $\pi(a) = \pi(b)$ for $a = \sum_{x \in \sigma'} a_x x$ and 
$b = \sum_{x \in \sigma'} b_x x$ of $P$, 
then $a = b$.
\end{Claim}

We choose $A, B \in \NN^{\sigma}$ such that
$A \cdot \sigma = a$, $B \cdot \sigma = b$ and
$\Supp(A), \Supp(B) \subseteq \sigma'$.
Then, by using the above claim,
there are $T$ and $T'$ of $\NN^{\sigma}$
with $\Supp(T), \Supp(T') \subseteq \sigma''$ and
$A \cdot \sigma + T \cdot \sigma = B \cdot \sigma + T' \cdot \sigma$.
Thus, $A + T = B + T'$.
For $x \in \sigma'$,
\[
A(x) = A(x) + T(x) = B(x) + T'(x) = B(x).
\]
Therefore, $A = B$.

\medskip
By the above claim, we can see $\alpha : \sigma' \to \bar{\sigma}$ is bijective.
From now on, let us see that $\bar{P}$ is of semistable type 
$(\bar{\sigma}, q_0, \bar{\Delta}, \bar{B})$ over $Q$ according to the conditions
(S1) -- (S4) of semistable structure.

(S1): We need to see $\Supp(\bar{\Delta}) \not= \emptyset$.
We assume the contrary.
Then, $\pi(x) = 0$ for all $x \in \Supp(\Delta)$. Thus,
\[
0 = \bar{f}(q_0) + \sum_{x \in \Sigma} B(x) \pi(x).
\]
Here $\bar{P}$ is sharp and $\bar{f}$ is integral. 
Therefore, we have $q_0 = 0$.
This is a contradiction.

(S2): Let us consider a natural commutative diagram:
\[
\begin{CD}
\NN^{\sigma'} @<{\sim}<{- \circ \alpha}< \NN^{\bar{\sigma}} \\
@VVV @VVV \\
P @>{\pi}>> \bar{P}
\end{CD}
\]
Then, Claim~\ref{claim:prop:quotient:semistable:2},
implies that $\NN^{\bar{\sigma}} \to \bar{P}$ is injective.

(S3): This is obvious in our case.

(S4): Finally, let us consider a relation
\[
\bar{T} \cdot \bar{\sigma} = \bar{f}(q) + \bar{T}' \cdot \bar{\sigma}
\quad(\bar{T}, \bar{T}' \in \NN^{\bar{\sigma}})
\]
with $q \not= 0$. We choose $T, T' \in \NN^{\sigma}$
such that $\Supp(T), \Supp(T') \subseteq \sigma'$,
$\bar{T} = T \circ \alpha^{-1}$ and $\bar{T}' = T' \circ \alpha^{-1}$.
Then, $\pi(T \cdot \sigma) = \pi(f(q) + T \cdot \sigma)$.
Thus, there are $S, S' \in \NN^{\sigma}$ such that
$\Supp(S), \Supp(S') \subseteq \sigma''$ and
\[
T \cdot \sigma + S \cdot \sigma = f(q) + T \cdot \sigma + S' \cdot \sigma.
\]
Therefore, $(T+S)(x) > 0$ for all $x \in \Supp(\Delta)$. In particular,
\[
T(x) = (T+S)(x) > 0
\]
for all $x \in \Supp(\Delta) \cap \sigma'$.
This means that $\bar{T}(\bar{x}) > 0$ for all $\bar{x} \in \Supp(\bar{\Delta})$.
\QED

Finally, we consider the following proposition.

\begin{Proposition}
\label{prop:semistable:monoid:saturated}
Let $f : Q \to P$ be an integral homomorphism of fine and sharp monoids.
If $P$ has a semistable structure
for some $\sigma \subseteq P$,
$q_0 \in Q$, $\Delta, B \in \NN^{\sigma}$,
then $\Coker(Q^{gr} \to P^{gr})$ is torsion free.
\end{Proposition}

\Proof
First of all, note that every element of
$P$ can be written as a form
$I \cdot \sigma + f(q)$ with $I \in \NN^{\sigma}$,
$q \in Q$ and $\Supp(\Delta) \not\subseteq \Supp(I)$.
Let us choose $x \in \Coker(Q^{gr} \to P^{gr})$
with $nx = 0$ for some positive integer $n$.
Then, there are $I_1, I_2 \in \NN^{\sigma}$,
$q_1, q_2 \in Q$ such that
$\Supp(\Delta) \not\subseteq \Supp(I_1)$, 
$\Supp(\Delta) \not\subseteq \Supp(I_2)$ and
$x$ is equal to 
\[
 (I_1 \cdot \sigma + f(q_1)) - (I_2 \cdot \sigma + f(q_2))
\]
as an element of $P^{gr}$.
Since $n x = 0$ in $\Coker(Q^{gr} \to P^{gr})$, 
there are $q_3, q_4 \in Q$ with
\[
 n(I_1 \cdot \sigma + f(q_1)) - n(I_2 \cdot \sigma + f(q_2)) =
 f(q_3) - f(q_4).
\]
Thus,
\[
 n I_1 \cdot \sigma + f(nq_1 + q_4) = n I_2 \cdot \sigma + f(nq_2 + q_3).
\]
Therefore, $n I_1 = n I_2$ because 
$\Supp(\Delta) \not\subseteq \Supp(n I_1)$ and
$\Supp(\Delta) \not\subseteq \Supp(n I_2)$. Hence $I_1 = I_2$.
Thus,
\[
 x = (I_1 \cdot \sigma + f(q_1)) - (I_2 \cdot \sigma + f(q_2))
 = f(q_1) - f(q_2),
\]
which means that $x = 0$ in $\Coker(Q^{gr} \to P^{gr})$.
\QED

\subsection{Extension of a ring for a good chart}
\setcounter{Theorem}{0}

Here we consider a ring extension to get a good chart.

\begin{Proposition}
\label{prop:split:log:structure}
Let $(A,m)$ be a noetherian local ring, $S = \Spec(A)$ and $s$ the
closed point of $S$.
Let $M_S$ be a fine log structure on $S$.
Then, there is
a local homomorphism $f : (A, m)  \to (B, n)$
of noetherian local rings  with the following properties:
\begin{enumerate}
\renewcommand{\labelenumi}{(\arabic{enumi})}
\item
$B/n$ is algebraic over $A/m$, and $f$ is flat and quasi-finite.

\item
Let $f^a : S' = \Spec(B) \to S = \Spec(A)$ be the induced homomorphism,
$s'$ the closed point of $S' = \Spec(B)$, and $M_{S'} = (f^a)^*(M_S)$.
There are  a fine and sharp monoid $Q$ and
a homomorphism $\pi_Q : Q \to M_{S',s'}$ such that
$Q \to M_{S',\bar{s}'} \to \overline{M}_{S',\bar{s}'}$ is bijective.
\end{enumerate}
\end{Proposition}

\Proof
Let us begin with the following lemma:

\begin{Lemma}
\label{lem:ext:class:zero}
Let $G$ be a finitely generated abelian group and
$R$ a ring.
Let us fix an element $\delta$ of $\Ext^1(G, R^{\times})$.
Then, there are $u_1, \ldots, u_l \in R^{\times}$ and
integers $a_1, \ldots, a_l \geq 2$ with the following property:
\begin{enumerate}
\renewcommand{\labelenumi}{(\arabic{enumi})}
\item
$a_1 \cdots a_l$ is equal to the order of the torsion part of $G$.

\item
For any homomorphism $f : R \to S$ of rings, if
there are $v_1, \ldots, v_l$ with
$v_i^{a_i} = f(u_i)$ for all $i$, then
the image of $\delta$ via the canonical homomorphism
\[
\Ext^1(G, R^{\times}) \to \Ext^1(G, S^{\times})
\]
is zero.
\end{enumerate}
\end{Lemma}

\Proof
By the fundamental theorem of abelian groups,
we have the following exact sequence:
\[
\begin{CD}
0 @>>> \ZZ^l @>{\phi}>> \ZZ^{l'} @>>> G @>>> 0,
\end{CD}
\]
where $\phi$ is given by $\phi(x_1, \ldots, x_l) = (a_1 x_1, \ldots, a_l x_l, 0, \ldots, 0)$
for some integers $a_1, \ldots, a_l \geq 2$. Note that
$a_1 \cdots a_l$ is equal to the order of the torsion part of $G$.
The above exact sequence gives rise to an exact sequence
\[
\begin{CD}
 \Hom(\ZZ^{l'}, R^{\times}) @>{\phi_R^*}>> \Hom(\ZZ^l, R^{\times}) @>{\alpha_R}>> \Ext^1(G, R^{\times}) @>>> \Ext^1(\ZZ^{l'}, R^{\times}).
\end{CD}
\]
Note that $\Ext^1(\ZZ^{l'}, R^{\times}) = \{ 0 \}$. Thus, there is $h \in \Hom(\ZZ^l, R^{\times})$ with
$\alpha_R(h) = \delta$.
We set $u_i = h(e_i)$ for $i=1, \ldots, l$.

Let $f : R \to S$ be any homomorphism of rings with
$v_i^{a_i} = f(u_i)$ ($i=1, \ldots, l$) for some $v_1, \ldots, v_l \in S$. Let us consider
the following commutative diagram:
\[
\begin{CD}
 \Hom(\ZZ^{l'}, R^{\times}) @>{\phi_R^*}>> \Hom(\ZZ^l, R^{\times}) @>{\alpha_R}>> \Ext^1(G, R^{\times}) @>>> 0 \\
 @V{g_1}VV @V{g_2}VV @V{g_3}VV @. \\
 \Hom(\ZZ^{l'}, S^{\times}) @>{\phi_{S}^*}>> \Hom(\ZZ^l, S^{\times}) @>{\alpha_S}>> \Ext^1(G, S^{\times}) @>>> 0 \\
\end{CD}
\]
Note that $g_2(h)(e_i) = f(u_i)$ for $i=1, \ldots, l$. Thus, if we set
$h' \in  \Hom(\ZZ^{l'}, S^{\times})$ by
\[
h' (e_i) = \begin{cases}
v_i & \text{if $i = 1, \ldots. l$} \\
0   & \text{if $i > l$}
\end{cases}
\]
then $\phi_S^*(h') = g_2(h)$. Therefore,
\[
g_3(\delta) = g_3(\alpha_R(h)) = \alpha_S(g_2(h)) = \alpha_S(\phi_S^*(h')) = 0.
\]
\QED

\bigskip
Let us start the proof of Proposition~\ref{prop:split:log:structure}.
Let $\delta \in \Ext^1(\overline{M}_{S,\bar{s}}^{gr}, \OO_{S,\bar{s}}^{\times})$ 
be the extension class of
\[
0 \to \OO_{S,\bar{s}}^{\times} \to M_{S,\bar{s}}^{gr} \to 
\overline{M}_{S,\bar{s}}^{gr} \to 0.
\]
Then, by Lemma~\ref{lem:ext:class:zero},
there are $u_1, \ldots, u_l \in \OO_{S,\bar{s}}^{\times}$ and
integers $a_1, \ldots, a_l$ with the properties as in Lemma~\ref{lem:ext:class:zero}.
Let us choose an \'{e}tale neighborhood $(U,u)$ of $s$ such that
$u_1, \ldots, u_l \in \OO_{U,u}^{\times}$.
Let $B$ be the localization of
\[
\OO_{U,u}[X_1, \ldots, X_l]/(X_1^{a_1} - u_1, \ldots, X_l^{a_l} - u_l).
\]
at a closed point over $u$. Then, $B$ is flat and quasi-finite over $A$.
Let $v_i$ be the class of $X_i$ in $B$. 
Note that $v_i^{a_i} = u_i$ in $B$ for all $i$.
Let $s'$ be the closed point of $S' = \Spec(B)$,
$\pi : S' \to S$ the canonical morphism, and $M_{S'} = \pi^*(M_S)$.
Then, we have an exact sequence
\[
0 \to \OO_{S',\bar{s}'}^{\times} \to M_{S', \bar{s}'}^{gr} \to 
\overline{M}_{S', \bar{s}'}^{gr} \to 0.
\]
Since $M^{gr}_{S',\bar{s}'}$ is the push-out 
$\OO_{S',\bar{s}'}^{\times} \times_{\OO_{S,\bar{s}}^{\times}} M^{gr}_{S,\bar{s}}$,
we can see that
$\overline{M}_{S', \bar{s}'}^{gr} = \overline{M}_{S, \bar{s}}^{gr}$ and
the extension class $\delta'$ of the above exact sequence is the image of the canonical
homomorphism $\Ext^1(\overline{M}_{S, \bar{s}}^{gr}, \OO_{S,\bar{s}}^{\times}) \to
\Ext^1(\overline{M}_{S', \bar{s}'}^{gr}, \OO_{S',\bar{s}'}^{\times})$.
Thus, by Lemma~\ref{lem:ext:class:zero}, $\delta' = 0$.
Therefore, we have a splitting
$s : \overline{M}_{S', \bar{s}'}^{gr} \to M_{S', \bar{s}'}^{gr}$ of
$M_{S', \bar{s}'}^{gr} \to \overline{M}_{S', \bar{s}'}^{gr}$.
Here we set $Q = \overline{M}_{S',\bar{s}'}$.
Let us see that $s(q) \in M_{S',\bar{s}'}$ for all $q \in Q$.
Indeed, if we denote $M_{S',\bar{s}'}^{gr} \to \overline{M}_{S',\bar{s}'}^{gr}$ by $\pi$,
then $\pi(s(q)) = q$. Thus, there are $u \in \OO_{S', \bar{s}'}^{\times}$ and
$m \in M_{S',\bar{s}'}$ with $s(q) = m \cdot u$, which implies
$s(q) \in M_{S',\bar{s}'}$.
Moreover, $Q \to  M_{S', \bar{s}'} \to  \overline{M}_{S', \bar{s}'}$ is the identity map.
Furthermore, changing $S'$ by an \'{e}tale neighborhood of $S'$, we can see that
$Q \to M_{S',\bar{s}}$ is defined on $S'$.
\QED

\subsection{Local structure theorem}
\setcounter{Theorem}{0}

Let $\alpha : M_X \to \OO_X$ be a log structure on a scheme $X$.
For $x \in X$, an element $p \in \overline{M}_{X,\bar{x}}$ is said to
be regular if there is $m \in M_{X,\bar{x}}$ such that
$p$ is $m$ modulo $\OO_{X,\bar{x}}^{\times}$ and
$\alpha(m)$ is a regular element of $\OO_{X,\bar{x}}$.
Note that the regularity of $p$ does not depend on the choice of $m$.

\begin{Theorem}
\label{thm:local:structure:theorem}
Let $(f, h) : (X, M_X) \to (S, M_S)$ be a smooth and  integral morphism of fine log schemes.
Let $x$ be an element of $X$ and $s = f(x)$.
We assume that $f : X \to S$ is semistable at $x$.
Then, we have the following:
\begin{enumerate}
\renewcommand{\labelenumi}{(\arabic{enumi})}
\item
If $f$ is smooth at $x$, then
there is a submonoid $N$
of $\overline{M}_{X,\bar{x}}$ such that
$\overline{M}_{X,\bar{x}} = 
\bar{h}_{\bar{x}}(\overline{M}_{S,\bar{s}}) \times N$ and
$N$ is isomorphic to $\NN^a$ 
for some non-negative integer $a$.
Further, every element of $N$ is regular.

\item
If $f$ is not smooth at $x$ and 
$\bar{h}_{\bar{x}} : \overline{M}_{S,\bar{s}} \to
\overline{M}_{X,\bar{x}}$ splits, there is a submonoid $N$
of $\overline{M}_{X,\bar{x}}$ such that
$\overline{M}_{X,\bar{x}} = 
\bar{h}_{\bar{x}}(\overline{M}_{S,\bar{s}}) \times N$ and
$N$ is isomorphic to the monoid arising from monomials of 
\[
\ZZ[U_1,U_2, \ldots, U_a]/(U_1^2 - U_2^2)
\]
for some $a \geq 2$.
In this case, the characteristic of the residue
field of $\OO_{X,\bar{x}}$ is not equal to $2$.
Further, every element of $N$ is regular.

\item
If $f$ is not smooth at $x$ and
$\bar{h}_{\bar{x}} : \overline{M}_{S,\bar{s}} \to
\overline{M}_{X,\bar{x}}$ does not split, then
$\overline{M}_{X,\bar{x}}$ has the unique semistable structure 
$(\sigma, q_0, \Delta, B)$ over $\overline{M}_{S,\bar{s}}$
for some $\sigma \subseteq \overline{M}_{X,\bar{x}}$ with $\#(\sigma) \geq 2$,
$q_0 \in \overline{M}_{S,\bar{s}}$ and $\Delta, B \in \NN^{\sigma}$.
More precisely, $\sigma$ is the set of all irreducible elements
of $\overline{M}_{X,\bar{x}}$ not lying in
$\bar{h}_{\bar{x}}(\overline{M}_{S,\bar{s}})$.
Further, every element of $\sigma \setminus \Supp(\Delta)$ is regular.
\end{enumerate}
\end{Theorem}

\Proof
Let us begin with the following lemma.

\begin{Lemma}
\label{lem:comp;pullback:monoid}
Let $f : X \to Y$ be a morphism of schemes and
$M_Y$ a fine log structure on $Y$.
If we set $M_X = f^*(M_Y)$, then, for any $x \in X$
and $y \in Y$ with $y = f(x)$,
the induced homomorphism
$\overline{M}_{Y,\bar{y}} \to \overline{M}_{X,\bar{x}}$
is bijective.
\end{Lemma}

\Proof
Let $P$ be a chart of $M_{Y, \bar{y}}$. Then,
$\overline{M}_{Y, \bar{y}}$ and $\overline{M}_{X, \bar{x}}$ is given
by
\[
P/\pi^{-1}(\OO_{Y,\bar{y}}^{\times})
\quad\text{and}\quad
P/{\pi'}^{-1}(\OO_{X,\bar{x}}^{\times})
\]
respectively, where $\pi : P \to \OO_{Y, \bar{y}}$ is the canonical morphism
and $\pi' : P \to \OO_{Y, \bar{y}} \to \OO_{X,\bar{x}}$
(cf. \cite{KatoLog} and \cite{IwaMw}).
Thus, it is sufficient to see that
$\pi^{-1}(\OO_{Y,\bar{y}}^{\times}) = {\pi'}^{-1}(\OO_{X,\bar{x}}^{\times})$.
Indeed, letting $m_{\bar{x}}$ and $m_{\bar{y}}$ be the maximal
ideals of $\OO_{X,\bar{x}}$ and $\OO_{Y,\bar{y}}$, and
$\alpha : \OO_{Y,\bar{y}} \to \OO_{X,\bar{x}}$ be the canonical homomorphism,
\[
p \in \pi^{-1}(\OO_{Y,\bar{y}}^{\times}) \Longleftrightarrow
\pi(p) \in \OO_{Y,\bar{y}}^{\times} \Longleftrightarrow \alpha(\pi(p)) \in\OO_{X,\bar{x}}^{\times}
\Longleftrightarrow p \in {\pi'}^{-1}(\OO_{X,\bar{x}}^{\times})
\]
because $\alpha(m_{\bar{y}}) \subseteq m_{\bar{x}}$ and
$\alpha(\OO_{Y,\bar{y}}^{\times}) \subseteq \OO_{X,\bar{x}}^{\times}$.
\QED

Let us go back to the proof of Theorem~\ref{thm:local:structure:theorem}.
Let us consider the geometric fiber $X_{\bar{s}}$ over $s$.
Then, by using Lemma~\ref{lem:comp;pullback:monoid},
we may assume that $S = \Spec(k)$ for 
some algebraically closed field $k$.
Thus, the theorem follows from \cite[Theorem~3.1]{IwaMw}
except the following facts:
\begin{enumerate}
\renewcommand{\labelenumi}{(\roman{enumi})}
\item
$N$ is isomorphic to the monoid $T$ arising from monomials of 
\[
 \ZZ[U_1,U_2, \ldots, U_a]/(U_1^2 - U_2^2)
\]
in the case (2).

\item
The regularity of elements of either $N$ or $\sigma \setminus \Supp(\Delta)$.
\end{enumerate}

\medskip
(i)
Let $T_k$ be the the monoid arising from monomials of 
\[
k[U_1,U_2, \ldots, U_a]/(U_1^2 - U_2^2).
\]
In order to see (i),
we need to show the natural homomorphism
$T \to T_k$ is bijective. 
Let $\bar{U}_1^{e_1}\bar{U}_2^{e_2} \cdots \bar{U}_{a}^{e_a}$ and 
$\bar{U}_1^{e'_1}\bar{U}_2^{e'_2} \cdots \bar{U}_{a}^{e'_a}$ be
elements of $T$. Clearly, we may assume that $e_1, e_1' \in \{ 0, 1\}$.
We suppose that $\bar{U}_1^{e_1}\bar{U}_2^{e_2} \cdots \bar{U}_{a}^{e_a} = \bar{U}_1^{e'_1}\bar{U}_2^{e'_2} \cdots \bar{U}_{a}^{e'_a}$
in $k[U_1,U_2, \ldots, U_a]/(U_1^2 - U_2^2)$.
Then, there is $\phi \in k[U_1, \ldots. U_a]$ with
\[
U_1^{e_1}U_2^{e_2} \cdots U_{a}^{e_a} - U_1^{e'_1}U_2^{e'_2} \cdots U_{a}^{e'_a}= 
(U_1^2 - U_2^2) \phi.
\]
Comparing the degrees with respect to $U_1$ of both sides,
we can see that $\phi = 0$. Therefore, $(e_1, \ldots, e_a) = (e'_1, \ldots, e'_a)$.

\medskip
(ii) Let $(\OO_{S,s}, m_{S,s})  \to (A,m)$ be a flat local homomorphism
of local rings.
We set $S' = \Spec(A)$, $X' = X \times_{S} S'$ and
the induced morphisms as follows:
\[
 \begin{CD}
  X' @>{\pi'}>> X \\
  @V{f'}VV @VV{f}V \\
  S' @>{\pi}>> S.
 \end{CD}
\]
Let us choose $x' \in X'$ with $f'(x') = m$ and $\pi'(x') = x$.
Then, since $\OO_{X,x} \to \OO_{X',x'}$ is faithfully flat,
using Lemma~\ref{lem:comp;pullback:monoid},
if regularity holds at $x'$, then so does at $x$.

Let $k$ be the algebraic closure of the residue field at $x$.
Note that by virtue of [EGA~III, Chapter~0, 10.3.1],
there are a noetherian local ring $(A,m)$ and a local homomorphism
$(\OO_{S, s}, m_{S,s}) \to (A, m)$ such that 
$m_{S,s}A = m$, $A/m$ is isomorphic to $k$ over $\OO_{S, s}/m_{S,s}$ 
and that $A$ is flat over $\OO_{S,s}$. 
Therefore, we may assume that $\OO_{S,s}/m_{S,s}$ is
algebraically closed and $x$ is a closed point.
Moreover, by using Proposition~\ref{prop:split:log:structure},
we may further assume that
there are a fine and sharp monoid $Q$ and a homomorphism
$\pi_Q : Q \to M_{S,s}$ such that
$Q \to M_{S, \bar{s}} \to \overline{M}_{S,\bar{s}}$ is bijective. 
Hence, there is a fine and sharp monoids $P$ together with
homomorphisms $f : Q \to P$ and
$\pi_P : P \to M_{X,\bar{s}}$ such that
the following properties are satisfied:
\begin{enumerate}
\renewcommand{\labelenumi}{(\alph{enumi})}
\item
The diagram
\[
 \begin{CD}
  Q @>f>> P \\
  @V{\pi_Q}VV @VV{\pi_P}V \\
  M_{S, \bar{s}} @>>> M_{X, \bar{x}},
 \end{CD}
\]
is commutative.

\item
The induced homomorphism
$P \to M_{X,\bar{x}} \to \overline{M}_{X, \bar{x}}$ is bijective.

\item
The natural homomorphism
\[
\OO_{S,\bar{s}} \otimes_{\OO_{S,\bar{s}}[Q]} \OO_{S,\bar{s}}[P] 
\to \OO_{X,\bar{s}}
\] 
is smooth.
\end{enumerate}
Since 
$\OO_{S,\bar{s}} \otimes_{\OO_{S,\bar{s}}[Q]} \OO_{S,\bar{s}}[P] 
\to \OO_{X,\bar{s}}$
is smooth,
it is sufficient to see the regularity
of each element in 
$\OO_{S,\bar{s}} \otimes_{\OO_{S,\bar{s}}[Q]}\OO_{S,\bar{s}}[P]$.

If there is a submonoid $N$ of $P$ with $P = f(Q) \times N$,
then
\[
  \OO_{S,\bar{s}} \otimes_{\OO_{S,\bar{s}}[Q]}\OO_{S,\bar{s}}[P]
= \OO_{S,\bar{s}}[N].
\]
Thus, the assertions follow from Lemma~\ref{x2:y2:regular} below.

Next, we assume that
$f : Q \to P$ does not splits.
Let us set $\sigma = \{ p_1, \ldots, p_r\}$ such that
$\Supp(\Delta) = \{ p_1, \ldots, p_l \}$.
Moreover, we set $x_i = \alpha(\pi_P(p_i))$ and
$t = \beta(\pi_Q(q_0))$, where 
$\alpha : M_X \to \OO_X$ 
and $\beta : M_{S,\bar{s}} \to
\OO_{S,\bar{s}}$
are the canonical homomorphisms.
Then,
\[
 \OO_{S,\bar{s}} \otimes_{\OO_{S,\bar{s}}[Q]}\OO_{S,\bar{s}}[P] =
 \OO_{S,\bar{s}}[X_1, \ldots, X_r]/
 (X_1 \cdots X_l - t X_{l+1}^{b_{l+1}} \cdots X_r^{b_r}),
\]
where $b_i = B(p_i)$ and $x_i$ is the class of $X_i$.
Thus, the assertions follow from Lemma~\ref{x2:y2:regular} below.
\QED

\begin{Lemma}
\label{x2:y2:regular}
Let $A$ be a ring.
Then, we have the following:
\begin{enumerate}
\renewcommand{\labelenumi}{(\arabic{enumi})}
\item
Let $A[ X ]$ be the polynomial ring of
one variable over $A$. For a regular element $a \in A$, $X$ is regular
in $A[X]/(X^2 -a)$,
that is,
the multiplication of $X$ in $A[X]/(X^2 -a)$
is injective.

\item
We assume $A$ is a local ring with the maximal ideal $m$.
Let $A[X_1, \ldots, X_l]$ be the polynomial ring of
$l$-variables over $A$. For $a \in m$, let us consider a ring $R$
given by
$R = A[X_1, \ldots, X_l]/(X_1 \cdots X_l - a)$.
If $\alpha$ is a regular element of $A$, then so is $\alpha$ in $R$.
\end{enumerate}
\end{Lemma}

\Proof
(1) We assume that $X f(X) = (X^2 - a)g(X)$ for some
$f(X), g(X) \in A[ X ]$.
We set $g(X) = Xh(X) + c$ for some $h(X) \in A[ X ]$
and $c \in A$.
Then, 
\[
ca = X(h(X)(X^2-a) + cX - f(X)).
\]
Thus, $ca = 0$. Since $a$ is regular, $c$ must be zero.
Therefore,
\[
Xf(X) = X(X^2-a)h(X),
\]
which implies $f(X) = (X^2 - a)h(X)$
because $X$ is regular in $A[X]$.

\medskip
(2)
Let $\widehat{R}$ be the completion with respect to $(m, X_1, \ldots, X_n)$.
Since $R \to \widehat{R}$ is faithfully flat,
it is sufficient to see the homomorphism $\widetilde{\alpha}
: \widehat{R} \to \widehat{R}$ given by the multiplication of $\alpha$
is injective. Note that $\widehat{R}$ is the direct products of many copies of
$\widehat{A}$ by Lemma~\ref{lem:flatness:x:1:x:l}.
Thus, $\widetilde{\alpha}$ is injective.
\QED

\subsection{The support of log structures}
\setcounter{Theorem}{0}
Here we consider the support of log structures.
The main result of this subsection is the following proposition:

\begin{Proposition}
\label{prop:closed:support:log:morphism}
Let $X$ be a scheme and let $M$ and $N$ be fine log structures on $X$.
Let $h : N \to M$ be a homomorphism of log structures, i.e.,
a homomorphism of sheaves of monoids with the following diagram commutative:
\[
\xymatrix{
  N \ar[rr]^{h} \ar[dr] & & M \ar[dl] \\
   & \OO_X & 
 }
\]
Then, the set 
$\{ x \in X \mid \text{$h_{\bar{x}} : N_{\bar{x}} \to M_{\bar{x}}$
is surjective} \}$
is open.
\end{Proposition}

\Proof
It is sufficient to show that if
$h_{\bar{x}} : N_{\bar{x}} \to M_{\bar{x}}$
is surjective, then there is an \'{e}tale neighborhood $U$ of $x$ such that,
for all $y \in U$, $h_{\bar{y}} : N_{\bar{y}} \to M_{\bar{y}}$
is surjective.
By virtue of \cite[(2.8)]{KatoLog}, for a suitable \'{e}tale neighborhood $U$ of $x$,
there are finitely generated monoids $P$ and $Q$ together with
homomorphisms $\pi : P \to \rest{M}{U}$,
$\mu : Q \to \rest{N}{U}$ and $f : Q \to P$ such that
$\pi$ and $\mu$ give rise to local charts of $M$ and $N$ respectively and
the diagram
\[
\begin{CD}
Q @>{f}>> P \\
@V{\mu}VV @VV{\pi}V \\
\rest{N}{U} @>{h_{U}}>> \rest{M}{U}
\end{CD}
\]
is commutative.
Let $\{ p_1, \ldots, p_n\}$ and $\{ q_1, \ldots, q_r \}$ be
generators of $P$ and $Q$ respectively.
Renumbering $p_1, \ldots, p_n$ and $q_1, \ldots, q_r$, we may assume that
\[
\begin{cases}
\pi(p_1), \ldots, \pi(p_{n'}) \in \OO_{X,\bar{x}}^{\times}, \ 
\pi(p_{n'+1}), \ldots, \pi(p_{n}) \not\in \OO_{X,\bar{x}}^{\times}, \\
\mu(q_1), \ldots, \mu(q_{r'}) \in \OO_{X,\bar{x}}^{\times}, \ 
\mu(q_{r'+1}), \ldots, \mu(q_{r}) \not\in \OO_{X,\bar{x}}^{\times}.
\end{cases}
\]
Let $P_0$ and $Q_0$ be submonoids of $P$ and $Q$ generated by
$p_1, \ldots, p_{n'}$ and $q_1, \ldots, q_{r'}$ respectively.
Here let us see the following:

\begin{Claim}
\label{claim:prop:closed:support:log:morphism:1}
$\pi_{\bar{x}}^{-1}(\OO_{X,\bar{x}}^{\times}) = P_0$ and
$\mu_{\bar{x}}^{-1}(\OO_{X,\bar{x}}^{\times}) = Q_0$.
\end{Claim}

Clearly, $P_0 \subseteq \pi_{\bar{x}}^{-1}(\OO_{X,\bar{x}}^{\times})$.
Conversely, we assume that $w \in \pi_{\bar{x}}^{-1}(\OO_{X,\bar{x}}^{\times})$.
We set
$w = p_1^{a_1} \cdots p_n^{a_n}$.
Then, 
\[
\pi(w) = \pi(p_1)^{a_1} \cdots \pi(p_n)^{a_n} \in \OO_{X,\bar{x}}^{\times}.
\]
Therefore, if $a_i > 0$, then $\pi(p_i) \in \OO_{X,\bar{x}}^{\times}$.
Hence, $a_i = 0$ for all $i > n'$, which means that $w \in P_0$.
In the same way, we can see that
$\mu_{\bar{x}}^{-1}(\OO_{X,\bar{x}}^{\times}) = Q_0$.

\begin{Claim}
\label{claim:prop:closed:support:log:morphism:2}
For $y \in U$, $h_{\bar{y}} : N_{\bar{y}} \to M_{\bar{y}}$ is
surjective if and only if $\bar{h}_{\bar{y}} : 
\overline{N}_{\bar{y}} \to \overline{M}_{\bar{y}}$ is
surjective.
\end{Claim}

Clearly, if $h_{\bar{y}} : N_{\bar{y}} \to M_{\bar{y}}$ is
surjective, then so is  $\bar{h}_{\bar{y}} : 
\overline{N}_{\bar{y}} \to \overline{M}_{\bar{y}}$.
Conversely, we assume that $\bar{h}_{\bar{y}} : 
\overline{N}_{\bar{y}} \to \overline{M}_{\bar{y}}$ is
surjective.
Let $m$ be an element of $M_{\bar{y}}$.
Then, there is $n \in N_{\bar{y}}$ such that
$m \equiv h_{\bar{y}}(n) \mod \OO_{X, \bar{x}}^{\times}$, i.e.,
$m = u h_{\bar{y}}(n)$ for some $u \in \OO_{X,\bar{y}}^{\times}$.
Thus, $m = u h_{\bar{y}}(n) = h_{\bar{y}}(un)$.

\medskip
Shrinking $U$ if necessarily, we may assume that
\[
\pi(p_1), \ldots, \pi(p_{n'}) \in \OO_{X,\bar{y}}^{\times}
\quad\text{and}\quad
\mu(q_1), \ldots, \mu(q_{r'}) \in \OO_{X,\bar{y}}^{\times}
\]
for all $y \in U$.
Here let us see that
$h_{\bar{y}} : N_{\bar{y}} \to M_{\bar{y}}$ is surjective for all $y \in U$,
which is equivalent to show that
$\bar{h}_{\bar{y}} : \overline{N}_{\bar{y}} \to \overline{M}_{\bar{y}}$ is surjective by
Claim~\ref{claim:prop:closed:support:log:morphism:2}.
Note that the commutative diagram
\[
\begin{CD}
Q @>{f}>> P \\
@V{\overline{\mu}_{\bar{y}}}VV @VV{\overline{\pi}_{\bar{y}}}V \\
\overline{N}_{\bar{y}} @>{\bar{h}_{\bar{y}}}>> \overline{M}_{\bar{y}}
\end{CD}
\]
gives rise to the commutative diagram
\[
\begin{CD}
Q/\mu_{\bar{y}}^{-1}(\OO_{X,\bar{y}}^{\times})  @>>> 
P/\pi_{\bar{y}}^{-1}(\OO_{X,\bar{y}}^{\times})  \\
@VVV @VVV \\
\overline{N}_{\bar{y}} @>{\bar{h}_{\bar{y}}}>> \overline{M}_{\bar{y}}
\end{CD}
\]
such that the vertical homomorphisms are bijective (cf. \cite{KatoLog} and \cite{IwaMw}).
Therefore, it is sufficient to see that
\[
Q/\mu_{\bar{y}}^{-1}(\OO_{X,\bar{y}}^{\times})  \to
P/\pi_{\bar{y}}^{-1}(\OO_{X,\bar{y}}^{\times})
\]
is surjective. Note that
$Q_0 \subseteq \mu_{\bar{y}}^{-1}(\OO_{X,\bar{y}}^{\times})$ and
$P_0 \subseteq \pi_{\bar{y}}^{-1}(\OO_{X,\bar{y}}^{\times})$.
Thus, we get the following commutative diagram:
\[
\begin{CD}
Q/Q_0 @>>> P/P_0 \\
@VVV @VVV \\
Q/\mu_{\bar{y}}^{-1}(\OO_{X,\bar{y}}^{\times})  @>>> 
P/\pi_{\bar{y}}^{-1}(\OO_{X,\bar{y}}^{\times}).
\end{CD}
\]
Here, by Claim~\ref{claim:prop:closed:support:log:morphism:1}, 
$Q/Q_0 \to P/P_0$ is surjective
because $\overline{N}_{\bar{x}} \to \overline{M}_{\bar{x}}$ is surjective.
Hence, so is $Q/\mu_{\bar{y}}^{-1}(\OO_{X,\bar{y}}^{\times})  \to
P/\pi_{\bar{y}}^{-1}(\OO_{X,\bar{y}}^{\times})$.
\QED

\begin{Corollary}
\label{cor:support:single:log}
Let $X$ be a scheme and $M$ a fine log structure on $X$.
Then, the set 
$\Supp(M) = \{ x \in X \mid \text{$M_{\bar{x}}$ is not trivial\,} \}$
is closed.
\end{Corollary}

\Proof
There is a natural homomorphism
$\OO_X^{\times} \to M$. Thus, this is a consequence of the above proposition.
\QED

\begin{Corollary}
\label{cor:support:log;morphism}
Let $X$ and $Y$ be schemes and let $M$ and $N$ be fine log structures on $X$
and $Y$ respectively.
Let $(f, h) : (X, M) \to (Y, N)$ be a log morphism.
Then, the set 
\[
\Supp(M/N) = \{ x \in X \mid \text{
$N_{\overline{f(x)}} \times \OO_{X,\bar{x}}^{\times} \to M_{\overline{x}}$
is not surjective} \}
\]
is closed.
\end{Corollary}

\Proof
Note that
the surjectivity of $N_{\overline{f(x)}} \times \OO_{X,\bar{x}}^{\times} \to M_{\overline{x}}$
is equivalent to the sujectivity of
$f^{*}(N)_{\bar{x}} \to M_{\bar{x}}$.
Thus, it follows from Proposition~\ref{prop:closed:support:log:morphism}.
\QED

\renewcommand{\theTheorem}{\arabic{section}.\arabic{Theorem}}
\renewcommand{\theClaim}{\arabic{section}.\arabic{Theorem}.\arabic{Claim}}
\renewcommand{\theequation}{\arabic{section}.\arabic{Theorem}.\arabic{Claim}}

\section{Rigidity theorem}

First, we would like to define the admissibility of morphisms.
Let $k$ be an algebraically closed field, and let $\phi : X \to Y$ be
a morphism of algebraic schemes over $k$.
Let $Z$ be a subscheme of $X$.
We say $\phi$ is admissible with respect to $Z$
if, for any irreducible component $X'$ of $X$,
$\phi(X') \not\subset Z$.

Let $M_Y$ and $M_k$ be fine log structures of $Y$ and $\Spec(k)$ respectively, and let
$(Y, M_Y) \to (\Spec(k), M_k)$ be a log morphism.
Then, as in Corollary~\ref{cor:support:log;morphism},
$\Supp(M_{Y}/M_k)$ is given by
\[
 \{ y \in Y \mid
\text{$M_k \times \OO_{Y, \bar{y}}^{\times} \to M_{Y,\bar{y}}$ is
not surjective} \}. 
\]
We say $\phi : X \to Y$ is admissible with respect to $M_Y/M_k$ if
$\phi : X \to Y$ is admissible with respect to $\Supp(M_{Y}/M_k)$.

In general, let $X \to S$ and $Y \to S$ be schemes of finite type
over a locally noetherian scheme $S$, and
let $M_Y$ and $M_S$ be fine log structures of $Y$ and $S$.
Let $\phi : X \to Y$ be a morphism over $S$.
For a point $s \in S$, we say
$\phi : X \to Y$ is admissible over $s$ with
respect to $M_Y/M_S$, if
$\phi_{\bar{s}} : X_{\bar{s}} \to Y_{\bar{s}}$ is admissible with
respect to $(\rest{M_Y}{Y_{\bar{s}}})/(\rest{M_S}{\bar{s}})$,
where $X_{\bar{s}}$ and $Y_{\bar{s}}$ are the geometric fibers
of $X$ and $Y$ over $s$.
If $\phi : X \to Y$ is admissible over any points of $S$ with respect to
$M_Y/M_S$, then $\phi$ is said to be admissible with respect to $M_Y/M_S$.

The following theorem is one of the main theorems of this paper.

\begin{Theorem}
\label{thm:rigid:log:morphism}
Let $X$, $Y$ and $S$ be locally noetherian schemes, and
let $M_X$, $M_Y$ and $M_S$ be fine log structures of $X$, $Y$ and $S$
respectively. Let
$(X, M_X) \to (S, M_S)$ and $(Y, M_Y) \to (S, M_S)$
be integral and log smooth morphisms, and 
let $\phi : X \to Y$ be a morphism over $S$.
Let us fix a point $s \in S$.
We assume that the geometric fibers $X_{\bar{s}}$ and $Y_{\bar{s}}$
over $s$ are semistable varieties and that
$\phi : X \to Y$ is admissible over $s$ with respect to $M_Y/M_S$.
If 
\[
(\phi, h) : (X, M_X) \to (Y, M_Y)
\quad\text{and}\quad
(\phi,  h') :  (X, M_X) \to (Y, M_Y)
\]
are extensions of
$\phi : X \to Y$ as log morphisms over $(S, M_S)$,
then, for all closed points $x$ lying over $s$,
$h_{\bar{x}} = h'_{\bar{x}}$ as homomorphisms
$M_{Y, \overline{\phi(x)}} \to M_{X,\bar{x}}$ of the germs of \'{e}tale topology.
\end{Theorem}

\Proof
Since this is a local problem,
we may assume that $S = \Spec(A)$ for a noetherian local ring $(A, m)$. 
Let $\rho : (A,m) \to (B, n)$ be a local homomorphism of local rings
such that $B/n$ is algebraic over $A/m$.
We denote the closed point of $S$ by $s$ and the closed point of $S' = \Spec(B)$ by $s'$.
We set
$X' = X \times_{S} S'$, $Y' = Y \times_{S} S'$, $M_X'  = \pi_X^*(M_X)$,
$M_{Y'} = \pi_Y^*(M_Y)$, and $M_{S'} = \pi_S^*(M_S)$,
where $\pi_X : X' \to X$, $\pi_Y : Y' \to Y$ and $\pi_S : S' \to S$
are the canonical morphisms. 
\[
\begin{CD}
X' @>{\pi_X}>> X \\
@V{\phi_{S'}}VV @VV{\phi}V \\
Y' @>{\pi_Y}>> Y
\end{CD}
\]
Then, we have log morphisms
\[
(\phi_{S'}, h_{S'}), (\phi_{S'}, h'_{S'}) :
(X', M_{X'}) \to (Y', M_{Y'})
\]
over $(S', M_{S'})$.

\begin{Claim}
\label{claim:thm:rigid:log:morphism:1}
If $\rho$ is flat and $h_{S', \bar{x}'} = h'_{S', \bar{x}'}$ for
all closed points $x'$ lying over $s'$, then $h_{\bar{x}} = h'_{\bar{x}}$
for all closed points $x$ lying over $s$.
\end{Claim}

Let us choose a closed point $x \in X$ over $s$.
Then, there is a $x' \in X'$ such that $\pi_X(x') = x$ and $x'$ is lying over $s'$.
If we set $y = \phi(x)$ and $y' = \phi_{S'}(x')$, then
$\pi_{Y}(y') = y$. Here we consider the natural commutative diagram:
\[
\xymatrix{
\overline{M}_{Y,\bar{y}} \ar[r] \ar@<-1ex>[d]_{\overline{h}_{\bar{x}}} 
\ar@<1ex>[d]^{\overline{h}'_{\bar{x}}} & 
\overline{M}_{Y',\bar{y}'} \ar[d]^{\bar{h}_{S',\bar{x}'} = \bar{h}'_{S',\bar{x}'}}\\
\overline{M}_{X,\bar{x}} \ar[r] & \overline{M}_{X',\bar{x}'}
}
\]
By Lemma~\ref{lem:comp;pullback:monoid}, 
$\overline{M}_{Y,\bar{y}} \to
\overline{M}_{Y',\bar{y}'}$ and
$\overline{M}_{X,\bar{x}} \to
\overline{M}_{X',\bar{x}'}$ are bijective.
Thus, we can see that $\bar{h}_{\bar{x}} = \bar{h}'_{\bar{x}}$.
Let us pick up $w \in M_{Y,\bar{y}}$.
Then, since $\bar{h}_{\bar{x}} = \bar{h}'_{\bar{x}}$, there is $u \in \OO_{X,\bar{x}}^{\times}$
with $h_{\bar{x}}(w) = h_{\bar{x}}(w) \cdot u$.
Here $h_{S', \bar{x}'} = h'_{S', \bar{x}'}$. Thus, $u$ must be $1$ in $\OO_{X', \bar{x}'}$.
Note that $\OO_{X',\bar{x}'}$ is flat over $\OO_{X,\bar{x}}$.
Therefore, $u$ is the identity in $\OO_{X,\bar{x}}$.

\bigskip
Let $I$ be an ideal of $A$ with $I^2 = \{ 0 \}$, and $B = A/I$. 
Next we consider a case where $\rho$ is given by 
the natural homomorphism $A \to B$.

\begin{Claim}
\label{claim:thm:rigid:log:morphism:2}
We assume that $k = A/m$ is algebraically closed and
there are a fine and sharp monoid $Q$ and
a homomorphism $\pi_Q : Q \to M_{S,s}$ such that
$Q \to M_{S,\bar{s}} \to \overline{M}_{S,\bar{s}}$ is bijective.
If $h_{S', \bar{x}'}= h'_{S', \bar{x}'}$ for all closed points $x'$ lying over $s'$,
then $h_{\bar{x}} = h'_{\bar{x}}$
for all closed points $x$ lying over $s$.
\end{Claim}

Let $x$ be a closed point of $X$ lying over $s$, and $y = \phi(x)$.
First of all, by \cite{IwaMw},
there are finite and sharp monoids $P$ and $P'$ and
homomorphisms
$P \to M_{X, \bar{x}}$, $Q \to P$,
$P' \to M_{Y, \bar{y}}$, $Q \to P'$
with the following properties:
\begin{enumerate}
\renewcommand{\labelenumi}{(\arabic{enumi})}
\item
The induced homomorphisms $P \to M_{X,\bar{x}} \to \overline{M}_{X, \bar{x}}$ and
$P' \to M_{Y,\bar{y}} \to \overline{M}_{Y, \bar{y}}$ are bijective.

\item
The following diagrams are commutative:
\[
 \begin{CD}
  Q @>f>> P \\
  @VVV @VVV \\
  M_{S, \bar{s}} @>>> M_{X, \bar{x}},
 \end{CD}\qquad
 \begin{CD}
  Q @>f'>> P' \\
  @VVV @VVV \\
  M_{S, \bar{s}} @>>> M_{Y, \bar{y}}.
 \end{CD}
\]

\item
There are \'{e}tale neighborhoods $(U, x')$ and $(V, y')$ of $x$ and $y$
such that $P \to M_{X, \bar{x}}$ and $P' \to M_{Y, \bar{y}}$ are defined over $U$ and $V$
respectively, and that the natural morphisms
\[
U \to \Spec(A \otimes_{A[Q]} A[P] ) 
\quad\text{and}\quad
V \to \Spec(A \otimes_{A[Q']} A[P'] )
\]
are smooth at $x'$ and $y'$
respectively.
\end{enumerate}
Clearly we may assume that $P$, $P'$ and $Q$ are submonoids
of $M_{X,\bar{x}}$, $M_{Y,\bar{y}}$ and $M_{S,\bar{s}}$ respectively.
Moreover, by Lemma~\ref{lem:etale:not:subset} below,
the admissibility of $\phi_s$ guarantees that
for any irreducible components $T$ of $U_s$,
$\phi_s(T) \not\subseteq  \Supp(M_{V_s}/M_k)$,
where $U_s$ and $V_s$ are fibers of $U \to S$ and $V \to S$ over $s$.

Let $\sigma$ (resp. $\sigma'$) be the set of all irreducible
elements of $P$ not lying in $f(Q)$ (resp. the set of all irreducible
elements of $P'$ not lying in $f'(Q)$).
For $j \in \sigma$ and $i \in \sigma'$,
we denote $\alpha(j)$ by $x_j$ and $\alpha'(i)$ by $y_{i}$,
where $\alpha : M_{X,\bar{x}} \to \OO_{X,\bar{x}}$ and
$\alpha' : M_{Y,\bar{y}} \to \OO_{Y,\bar{y}}$ are the canonical
homomorphisms.
Moreover, $\rest{x_j}{U_s}$ and $\rest{y_i}{V_s}$ are denoted by
$x_{js}$ and $y_{is}$ respectively.
Let us consider $h$ and $h'$ on the geometric fibers $X_{\bar{s}}$ and $Y_{\bar{s}}$ over $s$.
Using Lemma~\ref{lem:comp;pullback:monoid} and \cite[Theorem~4.1]{IwaMw}, 
$\bar{h}_{x} = \bar{h}'_{x}$ as $P' \to P$. 
Thus, we can set as follows:
\addtocounter{Claim}{1}
\begin{equation}
h(i) = u_{i} \cdot (I_{i} \cdot \sigma + f(q_i))
\quad\text{and}\quad
h'(i) = u'_{i}\cdot (I_{i} \cdot \sigma + f(q_i)),
\end{equation}
where $q_i \in Q$, $I_{i} \in \NN^{\sigma}$ and $u_{i}, u'_{i}
\in \OO^{\times}_{X, \bar{x}}$.
Then, we have
\addtocounter{Claim}{1}
\begin{equation}
\phi^*(y_{i}) = f^*(\beta(q_i)) \cdot x^{I_{i}} \cdot u_{i} = 
f^*(\beta(q_i)) \cdot x^{I_{i}} \cdot u'_{i},
\end{equation}
where $\beta : M_{S,\bar{s}} \to \OO_{S,\bar{s}}$ is the canonical
homomorphism.
Here we consider the following cases:
\begin{enumerate}
\renewcommand{\labelenumi}{(\Alph{enumi})}
\item
$f : Q \to P$ splits and $f' : Q \to P'$ splits.

\item
$f : Q \to P$ does not split and $f' : Q \to P'$ splits.

\item
$f : Q \to P$ splits and $f' : Q \to P'$ does not split.

\item
$f : Q \to P$ does not split and $f' : Q \to P'$ does not split.
\end{enumerate}

\medskip
(Case A): In this case, 
there are submonoids $N$ and $N'$ of $P$ and $P'$ respectively
such that $P = f(Q) \times N$ and $P' = f'(Q) \times N'$. Note that
$\sigma$ and $\sigma'$ are nothing more than
the set of all irreducible elements of $N$ and $N'$ respectively.
Then, by the local structure theorem,
\[
 \Supp(M_{Y_s}/M_k) = \bigcup_{i \in \sigma'} 
\{ y_{is} = 0\}.
\]
around $\bar{y}$ on $Y_{s}$. Thus, using the admissibility of $\phi_s$,
$\phi_s^*(y_{is}) \not= 0$.
Hence, we can see that $q_i = 0$ for all $i \in \sigma'$.
Therefore,
\[
 x^{I_{i}} \cdot u_{i} = x^{I_{i}} \cdot u'_{i}
\]
for all $i \in \sigma'$.
Here, $x_j$'s are regular elements. Therefore, $u_{i} = u'_{i}$ for all
$i \in \sigma'$.

\medskip
(Case B): In this case,
there is a submonoid $N'$ of $P'$
such that $P' = f'(Q) \times N'$.
Moreover, $P$  is of semistable type
\[
(\sigma, q_0, \Delta, B)
\]
over $Q$ for some $q_0 \in Q$ and $\Delta, B \in \NN^{\sigma}$.
By the local structure theorem,
\[
 \Supp(M_{Y_s}/M_k) = \bigcup_{i \in \sigma'} \{ y_{is} = 0\} 
\]
around $\bar{y}$ on $Y_s$. Thus, by the admissibility of $\phi_s$
$\phi_s^*(y_{is}) \not= 0$, which implies that $q_i = 0$ for all $i \in \sigma'$.
Hence,
\[
 x^{I_{i}} \cdot u_{i} = x^{I_{i}} \cdot u'_{i}
\]
for all $i \in \sigma'$.
By virtue of the admissibility of $\phi_s$ again, we can see that
$\Supp(I_i) \cap \Supp(\Delta) = \emptyset$.
Thus, $x^{I_{i}}$'s are regular elements.
Therefore, $u_i = u'_i$ for all $i \in \sigma'$.

\medskip
(Case C): In this case, there is a submonoid $N$ of $P$
with $P = f(Q) \times N$.
$P'$  is of semistable type
\[
(\sigma', q'_0, \Delta', B')
\]
over $Q$ for some $q'_0 \in Q$ and $\Delta', B' \in \NN^{\sigma'}$.
Note that
\[
 \Supp(M_{Y_s}/M_k) = \Sing(Y_s) \cup 
 \bigcup_{i \in \sigma' \setminus \Supp(\Delta')} \{ y_{is} = 0 \}.
\]
around $\bar{y}$ on $Y_s$.

Let us see that if $\phi_s^*(y_{is}) \not= 0$ for some $i \in \sigma'$,
then $u_i = u'_i$.
Indeed, we have $q_i = 0$. Thus,
\[
 x^{I_{i}} \cdot u_{i} = x^{I_{i}} \cdot u'_{i}.
\]
for all $i \in \sigma'$.
Here, $x^{I_{i}}$'s are regular elements.
Hence, $u_i = u'_i$.

Therefore, we may assume that there is $i \in \sigma'$ with
$\phi_s^*(y_{is}) = 0$. Then, by the admissibility of $\phi_s$,
there is $i_0$ such that $\phi_s^*(y_{is}) \not=  0$
for all $i \in \sigma' \setminus \{ i_0 \}$. 
Note that $i_0 \in \Supp(\Delta')$.
Then, $u_i = u'_i$ for all $i \in \sigma' \setminus \{ i_0 \}$.
Let us consider a relation
\[
\Delta' \cdot \sigma' = f'(q'_0) + B' \cdot \sigma'. 
\]
Then, we have
\[
 \begin{cases}
 \sum_{i \in \Supp(\Delta')} h(i) = f(q'_0) + \sum_{i \in \Supp(B')} B'(i)h(i) \\
 \sum_{i \in \Supp(\Delta')} h'(i) = f(q'_0) + \sum_{i \in \Supp(B')} B'(i)h'(i)
 \end{cases}
\]
Here $h(i) = h'(i)$ for all $i \not= i_0$. Thus,
we can see that $h(i_0) = h'(i_0)$.

\medskip
(Case D):
In the final case, $P$ and $P'$ are of semistable type
\[
(\sigma, q_0, \Delta, B)
\quad\text{and}\quad
(\sigma',q'_0, \Delta', B')
\]
over $Q$ for some $q_0, q'_0 \in Q$,
$\Delta, B \in \NN^{\sigma}$ and $\Delta', B' \in \NN^{\sigma'}$.
For $j \in \Supp(\Delta)$ and $i \in \Supp(\Delta')$,
let $U_{js}$ and $V_{is}$ be the irreducible components
of $U_s$ and $V_s$ given by $x_{js} = 0$ and $y_{is} = 0$
respectively. By the admissibility of $\phi_s$, 
for each $j \in \Supp(\Delta)$, there is a unique
$i \in \Supp(\Delta')$ with $\phi_s(U_{js}) \subseteq V_{is}$.
This $i$ is denoted by $\mu(j)$.
Note that
\[
 \Supp(M_{Y_s}/M_k) = \Sing(Y_s) \cup 
\bigcup_{i \in \sigma' \setminus \Supp(\Delta')}
\{ y_{is} = 0 \} 
\]
around $\bar{y}$ on $Y_s$. Therefore,
if $i \not= \mu(j)$, then $\rest{\phi^*(y_i)}{U_j} \not = 0$.
Thus, we can see $q_i = 0$ for $i \not= \mu(j)$.

First of all, let us see that
$q_i = 0$ and $u_i = u'_i$ for all 
$i \in \sigma' \setminus \Supp(\Delta')$.
By the previous observation, $q_i = 0$.
Thus,
\[
 x^{I_{i}} \cdot u_{i} = x^{I_{i}} \cdot u'_{i}.
\]
By using the admissibility of $\phi_s$,
we can see that $\Supp(I_i) \cap \Supp(\Delta) = \emptyset$.
Thus, $x^{I_i}$ is a regular element.
Therefore, $u_i = u'_i$.

Here, we consider the case where
$\mu(\Supp(\Delta)) = \{ i_0 \}$.
Then, for $i \not=  i_0$,
$q_i = 0$. Thus,
\[
 x^{I_{i}} \cdot u_{i} = x^{I_{i}} \cdot u'_{i}.
\]
By virtue of the admissibility of $\phi_s$,
$\Supp(I_i) \cap \Supp(\Delta) = \emptyset$.
Hence, $x^{I_i}$ is a regular element.
Therefore, $u_i = u'_i$.
Let us consider a relation
\[
\Delta' \cdot \sigma' = f'(q'_0) + B' \cdot \sigma'. 
\]
Then, we have
\[
 \begin{cases}
 \sum_{i \in \Supp(\Delta')} h(i) = f(q'_0) + \sum_{i \in \Supp(B')} B'(i)h(i) \\
 \sum_{i \in \Supp(\Delta')} h'(i) = f(q'_0) + \sum_{i \in \Supp(B')} B'(i)h'(i)
 \end{cases}
\]
Since  $h(i) = h'(i)$ for all $i \not= i_0$,
we can see that $h(i_0) = h'(i_0)$.

Finally, let us consider the case
where
$\# \mu(\Supp(\Delta)) \geq 2$.
In this case, $q_i = 0$ for all $i \in \sigma'$.
Thus,
\[
 x^{I_{i}} \cdot u_{i} = x^{I_{i}} \cdot u'_{i}.
\]
By our assumption, $u_{i} \equiv u'_{i} \mod I\OO_{X,\bar{x}}$.
Note that $x_j$ $(j \not\in \Supp(\Delta))$ is regular. Thus, 
if we set $I'_{i} = \rest{I_{i}}{\Supp(\Delta)} 
\in \NN^{\Supp(\Delta)}$, then
\[
 x^{I'_{i}} \cdot u_{i} = x^{I'_{i}} \cdot u'_{i}
\]
for all $i \in \Supp(\Delta')$.
Moreover, by the admissibility of $\phi_s$,
$\Supp(I'_{i}) \cap \Supp(I'_{i'}) = \emptyset$ for all $i \not= i'$.
Further, let us consider a relation
\[
\Delta' \cdot \sigma' = f'(q'_0) + B' \cdot \sigma'. 
\]
Since $h(i) = h'(i)$ for all $\sigma' \setminus \Supp(\Delta')$, 
we have
\[
 \sum_{i \in \Supp(\Delta')} h(i) =
 \sum_{i \in \Supp(\Delta')} h'(i). 
\]
Thus, 
$\prod_{i \in \Supp(\Delta')} u_i = \prod_{i \in \Supp(\Delta')} u'_i$.
Here we set $v_{i} = u_{i}/u'_{i}$ for $i \in \Supp(\Delta)$.
Then, gathering the above observations, we have seen that
\[
 \begin{cases}
  \text{$x^{I'_{i}} = x^{I'_{i}} \cdot v_{i}$ for all 
        $i \in \Supp(\Delta')$,} \\
  \text{$v_{i} \equiv 1 \mod I \OO_{X,\bar{x}}$ for all 
        $i \in \Supp(\Delta')$,} \\
  \text{$\prod_{i \in \Supp(\Delta')} v_i = 1$,} \\
  \text{$\Supp(I'_{i}) \cap \Supp(I'_{i'}) = \emptyset$ for all $i \not= i'$.}
 \end{cases}
\]
Since $A \otimes_{A[Q]} A[P]  \to \OO_{U, u}$ is smooth,
$A \otimes_{A[Q]} A[P \times \NN^e]  \to \OO_{U,x'}$ is \'{e}tale
for some $e\geq 0$.
Let $o$ be the origin of $\Spec(A \otimes_{A[Q]} A[P \times \NN^e])$.
Then, the residue field of $A \otimes_{A[Q]} A[P \times \NN^e]$ at $o$ is $k$.
Moreover, the residue fields of $\OO_{U,x'}$ and $\OO_{X,\bar{x}}$ are $k$
because $k$ is algebraically closed.
Therefore, the completion of $A \otimes_{A[Q]} A[P \times \NN^e] $ at $o$ is isomorphic
to the completion of $\OO_{X,\bar{x}}$.
Thus, by Lemma~\ref{lem:unit:one} below, $v_{i} = 1$, that is,
$u_{i} = u'_{i}$ for all $i$.

\bigskip
Let $k = A/m$ and $\bar{k}$ the algebraic closure of $k$.
By virtue of [EGA~III, Chapter~0, 10.3.1],
there are a noetherian local ring $(B,n)$ and a local homomorphism
$A \to B$ such that $mB = n$, $B/n$ is isomorphic to $\bar{k}$ over $k = A/m$ and
that $B$ is flat over $A$. Thus, 
by Claim~\ref{claim:thm:rigid:log:morphism:1}, we may assume that
the residue field $k = A/m$ is algebraically closed.
Moreover, by Proposition~\ref{prop:split:log:structure}, we may further assume that
there are a find and sharp monoid $Q$ and
a homomorphism $\pi_Q : Q \to M_{S,s}$ such that
$Q \to M_{S,\bar{s}} \to \overline{M}_{S,\bar{s}}$ is bijective.

\medskip
Let $A_i = A/m^{i+1}$, $\rho_i : A_i \to A_{i-1}$ the canonical homomorphism and
$I_i = \Ker(\rho_i)$. Then, $A_0 = k$  and $I_i^2 = \{ 0 \}$ for $i \geq 1$.
We set $X_i = X \times_S \Spec(A_i)$, $M_{X_i} = \rest{M_X}{X_i}$,
$Y_i = Y \times_{S} \Spec(A_i)$, $M_{Y_i} = \rest{M_Y}{Y_i}$.
Moreover, the induced morphisms
$M_{Y_i} \to M_{X_i}$ and $M_{Y_i} \to M_{X_i}$ via $h$ and $h'$ are
denoted by $h_i$ and $h'_i$ respectively.
Note that $h_0 = h'_0$ at closed points of $X_s$ by \cite{IwaMw}.
Moreover, by Claim~\ref{claim:thm:rigid:log:morphism:2},
$h_n = h_n$ at closed points lying over $s$ implies that $h_{n+1} = h'_{n+1}$
at closed points lying over $s$.
Therefore, we have $h_n = h'_n$ at closed points of $X_s$ for all $n \geq 0$.
Let $x$ be a closed point of $X$ over $s$ and $y = \phi(x)$.
Since $\bar{h}_{\bar{x}} = \bar{h}'_{\bar{x}}$ as a homomorphism
$\overline{M}_{Y,\bar{y}} \to \overline{M}_{X, \bar{x}}$,
for $w \in M_{Y,\bar{y}}$, there is $u \in \OO_{X,\bar{x}}^{\times}$
with $h_{\bar{x}}(w) = h'_{\bar{x}}(w) \cdot u$.
Since $h_n = h'_n$, we can see that $u - 1 \in m^{n+1} \OO_{X,\bar{x}}$.
Note that $\OO_{X,\bar{x}}$ is noetherian, which implies that
$\bigcap_{n=0} m^{n+1} \OO_{X,\bar{x}} = \{ 0 \}$.
Therefore, $u = 1$.
\QED

As corollary of Theorem~\ref{thm:rigid:log:morphism},
we have the following:

\begin{Corollary}[Rigidity theorem]
\label{cor:rigid:log:morphism}
Let $f : X \to S$ and $g : Y \to S$ be semistable schemes over a locally
noetherian scheme $S$, and let $\phi : X \to Y$ be a morphism over $S$.
Let $M_X$, $M_Y$ and $M_S$ be fine log structures on $X$, $Y$ and $S$
respectively.
We assume that $(X, M_X)$ and $(Y, M_Y)$ are log smooth and integral over
$(S, M_S)$ and $\phi$ is admissible with respect to $M_Y/M_S$.
If we have log morphisms
\[
(\phi, h) : (X, M_X) \to (Y, M_Y)
\quad\text{and}\quad
(\phi, h') : (X, M_X) \to (Y, M_Y)
\]
over $(S, M_S)$ as extensions of $\phi : X \to Y$,
then $h = h'$.
\end{Corollary}

The following two lemmas was needed for the proof of
Theorem~\ref{thm:rigid:log:morphism}.

\begin{Lemma}
\label{lem:etale:not:subset}
Let
\[
\begin{CD}
X' @>{\pi'}>> Y' \\
@V{\mu}VV @VV{\nu}V \\
X @>{\pi}>> Y
\end{CD}
\]
be  a commutative diagram of reduced algebraic schemes 
over an algebraically closed field such that $X$ and $X'$ is equi-dimensional
and $\mu$ is flat.
Let $Z$ be a closed subset of $Y$.
If $\pi(T) \not\subseteq Z$ for any irreducible components $T$ of $X$,
then $\pi'(T') \not\subseteq \nu^{-1}(Z)$ for 
any irreducible components $T'$ of $X'$.
\end{Lemma}

\Proof
We assume that $\pi'(T') \subseteq \nu^{-1}(Z)$ for 
an irreducible component $T'$ of $X'$. Then,
\[
\pi(\mu(T')) = \nu(\pi'(T')) \subseteq \nu(\nu^{-1}(Z)) \subseteq Z.
\]
Let $T$ be the Zariski closure of $\mu(T')$.
If $\dim T < \dim X$, then
\[
\dim \mu^{-1}(x) \geq \dim T' - \dim T >  \dim X' - \dim X
\]
for $x \in \mu(T')$, which is a contradiction because
$\mu$ is flat.
Thus, we have $\dim T = \dim X$, which means that
$T$ is an irreducible component of $X$.
On the other hand, we know $\pi(T) \subseteq Z$.
This is contradict to our assumption. Therefore, we get our lemma.
\QED

\begin{Lemma}
\label{lem:unit:one}
Let $(A, m)$ be a noetherian complete
local ring and $A\lformal X_1, \ldots, X_n \rformal$
the ring of formal power series of $n$-variables over $A$.
For a fixed $a \in m$, let
\[
R = A\lformal X_1, \ldots, X_n \rformal /(X_1 \cdots X_n - a)
\]
and $J$ an ideal of $R$ with $J^2 = 0$.
Let $u_1, \ldots u_l$ be elements of $R$ and
$I_1, \ldots, I_l$ elements of $\NN^n$ with 
$\Supp(I_i) \cap \Supp(I_j) = \emptyset$ for $i \not= j$.
We assume that \rom{(1)} $u_1 \cdots u_l = 1$,
\rom{(2)} $X^{I_i} u_i = X^{I_i}$ in $R$ for all $i$, and
that \rom{(3)} $u_i \equiv 1 \mod J$.
Then, we have $u_1 = \cdots = u_l = 1$.
\end{Lemma}

\Proof
We set $\Sigma = \{ I \in \NN^n \mid \Delta \not\leq I \}$
and
\[
 A\lformal X_1, \ldots, X_n \rformal_{\Sigma} =
\left\{ \sum_{I \in \Sigma }a_I X^I \mid a_I \in A \right\},
\]
where $\Delta = (1, \ldots, 1)$.
Then, by Lemma~\ref{lem:flatness:x:1:x:l},
the natural map $A\lformal X_1, \ldots, X_n \rformal_{\Sigma}
\to R$ is bijective. Here we claim the following:

\begin{Claim}
Let $T$ be an element of $\NN^n$.
We set $\Sigma_T = \{ I \in \Sigma \mid I + T \geq \Delta\}$.
Then, for $f \in A\lformal X_1, \ldots, X_n \rformal_{\Sigma}$,
if $X^T f = 0$ in $R$, then
$f$ can be written by a form
\[
 f = \sum_{I \in \Sigma_T} b_I X^I.
\]
\end{Claim}

If either $T = (0, \ldots, 0)$ or $T \geq \Delta$, then our assertion
is trivial. Thus, we may assume that $T \not= (0, \ldots, 0)$ and
$T \not\geq \Delta$.
For $ I \in \NN^n$, we can find a
non-negative integer $a$ and $J \in \Sigma$
with $I = a\Delta + J$.
We denote $a$ and $J$ by $a(I)$ and $J(I)$ respectively.
Here let us see that
$J(I+T) \not\in \{ S + T \mid S \in \Sigma \setminus \Sigma_T\}$
for $I \in \Sigma_T$. Indeed, since $I \in \Sigma_T$,
we can find $i$ with $I(i) = 0$ and $T(i) > 0$.
Thus,
\[
 J(I + T)(i) = T(i) - a(I+T) < T(i).
\]
Hence $J(I+T) \not\in \{ S + T \mid S \in \Sigma \setminus \Sigma_T\}$.

Here we set $f = \sum_{I \in \Sigma} a_I X^I$.
Then,
\begin{align*}
 X^T f & = \sum_{I \in \Sigma_T} a_I X^{I + T} + 
         \sum_{I \in \Sigma \setminus \Sigma_T} a_I X^{I + T} \\
       & = \sum_{I \in \Sigma_T} a_I a^{a(I+T)}X^{J(I+T)} + 
         \sum_{I \in \Sigma \setminus \Sigma_T} a_I X^{I+T}.
\end{align*}
Thus, $a_I = 0$ for $I \in \Sigma \setminus \Sigma_T$.

\bigskip
Since $u_i \equiv 1 \mod J$, there is $a_i \in J$ with $u_i = 1 + a_i$.
Then, $X^{I_i} a_i = 0$. Moreover, since $J^2 = 0$,
\[
 u_1 \cdots u_l = 1 + a_1 + \cdots + a_l = 1.
\]
Hence $a_1 + \cdots + a_l = 0$.
Since $X^{I_i} a_i = 0$, by the above claim,
$a_i = \sum_{I \in \Sigma_{I_i}} c_{i, I} X^I$,
where $\Sigma_{I_i} = \{ I \in \Sigma \mid I + I_i \geq \Delta\}$.
Therefore,
\[
 \sum_{i=1}^l \sum_{I \in \Sigma_{I_i}} c_{i,I} X^{I} = 0.
\]
Note that if $I \in \Sigma_{I_i}$ and $I' \in \Sigma_{I_j}$ for
$i \not= j$, then $I \not= I'$ because 
$\Supp(I_i) \cap \Supp(I_j) = \emptyset$.
Thus, we can see that $c_{i, I} = 0$, 
which shows us $a_i = 0$ for all $i$.
\QED

\section{Descent of log morphisms}

Let $A$ be a ring.
Note that $A$ gives rise to the commutative monoid 
$(A, \times)$ with respect to the multiplication.
{\em A pre-log monoid over $A$} is a monoid $M$ together with
a homomorphism $\alpha : M \to (A, \times)$.
Moreover, a pre-log monoid $\alpha : M \to (A, \times)$ over $A$
is called {\em a log monoid} if
$\alpha$ induces
an isomorphism $\alpha^{-1}(A^{\times}) 
\overset{\sim}{\longrightarrow} A^{\times}$.
We often identify $\alpha^{-1}(A^{\times})$ with $A^{\times}$.
First, let us see the following lemma:

\begin{Lemma}
\label{lem:associated:log:monoid}
Let $A$ be a ring and $\alpha : M \to (A, \times)$ a pre-log monoid over $A$.
Then, we have the following:
\begin{enumerate}
\renewcommand{\labelenumi}{(\arabic{enumi})}
\item
A homomorphism $\alpha' : M \pushout_{\alpha^{-1}(A^{\times})} A^{\times} \to A$
given by $\alpha'(m, a) = \alpha(m) a$ is well-defined.

\item
$(\alpha')^{-1}(A^{\times}) \overset{\sim}{\longrightarrow} A^{\times}$, that is,
$\alpha' : M \pushout_{\alpha^{-1}(A^{\times})} A^{\times} \to A$ is a log monoid over $A$.
\end{enumerate}
The log monoid $\alpha' : M \pushout_{\alpha^{-1}(A^{\times})} A^{\times} \to A$
as above is called the associated log monoid of $\alpha : M \to A$.
\end{Lemma}

\Proof
(1) This is obvious because the diagram
\[
\begin{CD}
\alpha^{-1}(A^{\times}) @>{\alpha}>> A^{\times} \\
@VVV @VVV \\
M @>{\alpha}>> A
\end{CD}
\]
is commutative.

(2) If the class of $(m, a)$ belongs to $(\alpha')^{-1}(A^{\times})$,
then it is easy to see that $(m, a) \sim (1, \alpha(m)a)$.
Moreover, if $(1, a) \sim (1, a')$, then $a = a'$.
Thus, we get our assertion.
\QED

Let $\alpha : M \to A$ be a pre-log monoid over $A$ and
$f : A \to B$ a homomorphism of rings.
The associated log monoid 
$M \pushout_{(f \circ \alpha)^{-1}(B^{\times})} B^{\times}$
of $f \circ \alpha : M \to (B, \times)$ is denoted by $M \boxtimes_A B$
and the canonical homomorphism $M \boxtimes_A B \to B$ is denoted by
$\alpha_B$, that is, $\alpha_B(m, b) = f(\alpha(m)) b$.
Note that the associated log monoid of $\alpha : M \to A$
is nothing more than $M \boxtimes_A A$.

We say $f$ is {\em quasi-local}  if $f^{-1}(B^{\times}) = A^{\times}$.
Note that a local homomorphism of local rings is quasi-local.
More generally, if, for any maximal ideal $m$ of $A$,
there is a maximal ideal $n$ of $B$ with $f^{-1}(n) = m$,
then $f$ is quasi-local.
If $M$ is a log monoid over $A$ and $f : A \to B$ is quasi-local, then
$M \boxtimes_A B = M \pushout_{A^{\times}} B^{\times}$, namely,
$(m, b) \sim (m', b')$ if and only if
there is $a \in A^{\times}$ with $(m, b)(a, f(a^{-1})) = (m', b')$.

Let $\beta : N \to A$ be another pre-log monoid over $A$.
A homomorphism $\phi : M \to N$ is called a homomorphism over $A$
if $\alpha = \beta \circ \phi$, i.e.,
the following diagram is commutative:
\[
 \xymatrix{
  M \ar[rr]^{\phi} \ar[dr]_{\alpha}& & N \ar[dl]^{\beta} \\
   & A & 
 }
\]
Note that if $\alpha : M \to A$ and $\beta : N \to A$ are log monids over $A$, and
$\phi$ is a homomorphism over $A$,
then $\phi(am) = a\phi(m)$ for all $m \in M$ and $a \in A^{\times}$.
We denote the set of all homomorphisms $M \to N$ over $A$ by
$\Hom_A(M,N)$.

More generally, let $f : A \to B$ be a homomorphism
of rings and let $\alpha : M \to A$ and
$\beta : N \to B$ be pre-log monoids over $A$ and $B$ respectively.
Then, a homomorphism $\phi : M \to N$ is called a homomorphism over $f$
if $f \circ \alpha = \beta \circ \phi$, i.e.,
the diagram
\[
 \xymatrix{
  M \ar[r]^{\phi} \ar[d]_{\alpha}&  N \ar[d]^{\beta} \\
   A \ar[r]^f & B
 }
\]
is commutative. The set of all homomorphisms $\phi : M \to N$ over $f$
is denoted by $\Hom_f(M, N)$.
 
\begin{Lemma}
\label{lem:adjoint:log:monoid}
Let $f : A \to B$ be a homomorphism
of rings, and
let $\alpha : M \to A$ and $\beta : N \to B$ be log monoids over $A$ and $B$ respectively.
Then, there is a natural isomorphism
\[
 \Hom_f(M, N) \overset{\sim}{\longrightarrow} 
 \Hom_B( M \boxtimes_{A} B, N).
\]
\end{Lemma}

\Proof
We define 
\[
\gamma : \Hom_f(M, N) \to \Hom_B( M \boxtimes_{A} B, N)
\]
to be $\gamma(\phi)(m, b) = \phi(m)b$.
This is well defined because
the following diagram 
\[
\begin{CD}
(f \circ \alpha)^{-1}(B^{\times}) @>{f \circ \alpha}>> B^{\times} \\
@VVV @VV{\beta^{-1}}V \\
M @>{\phi}>>  N
\end{CD}
\]
is commutative.
We need to see that $\gamma(\phi)$ is a homomorphism over $B$.
Indeed,
\[
\beta \circ \gamma(\phi) (m, b) = \beta(\phi(m)b) =\beta(\phi(m)) b
= f(\alpha(m)) b = \alpha_B(m, b).
\]

Next, we define
\[
\delta : \Hom_B( M \boxtimes_{A} B, N) \to \Hom_f(M, N)
\]
to be $\delta(\phi')(m) = \phi'(m, 1)$. Since
\[
\beta \circ \delta(\phi')(m) = \beta(\phi'(m,1)) = \alpha_B(m,1) = f \circ\alpha(m),
\]
$\delta(\phi')$ is a homomorphism over $f$.

With this notation,
\[
(\delta \circ \gamma)(\phi)(m) = \delta(\gamma(\phi))(m) = \gamma(\phi)(m, 1)
= \phi(m)
\]
and
\[
(\gamma \circ \delta)(\phi')(m, b) = \gamma(\delta(\phi'))(m) =
\delta(\phi')(m)b = \phi'(m, 1)b = \phi'(m, b).
\]
Therefore, we get our claim.
\QED

Let $\alpha : M \to A$ and
$\beta : N \to A$ be pre-log monoids over $A$.
Let $\phi : M \to N$ be a homomorphism over $A$ and
$h : B \to C$ be a homomorphism of $A$-algebras.
Then, there is a homomorphism
\[
 \phi \boxtimes h : 
 M \boxtimes_{A} B \to N \boxtimes_{A} C
\]
given by $(\phi \boxtimes h)(m, b) = (\phi(m), h(b))$.
This is well-defined. Indeed, let $f : A \to B$ and $g : A \to C$ be
the canonical homomorphisms with $h \circ f = g$.
Here we consider $\mu : M \to N \boxtimes_A C$ and
$\nu : B^{\times} \to N \boxtimes_A C$ given by
$\mu(m) = (\phi(m), 1)$ and $\nu(b) = (1, h(b))$.
Then, for $m \in (f \circ \alpha)^{-1}(B^{\times})$,
\begin{align*}
(1, h((f\circ \alpha)(m))) \cdot (\phi(m), 1) & =
(\phi(m), 1) \cdot (1, h(f(\alpha(m)))) \\
& = (\phi(m), 1) \cdot (1, g(\alpha(m))) \\
& = (\phi(m), 1) \cdot (1, (g \circ \beta)(\phi(m))),
\end{align*}
which means that $(1, h((f\circ \alpha)(m))) \sim (\phi(m), 1)$
in $N \boxtimes_A C$. Thus, the diagram
\[
\begin{CD}
(f \circ \alpha)^{-1}(B^{\times}) @>>> B^{\times} \\
@VVV @VV{\nu}V \\
M @>{\mu}>> N \boxtimes_A C
\end{CD}
\]
is commutative.
Therefore, we have the homomorphism
$\phi \boxtimes h : 
 M \boxtimes_{A} B \to N \boxtimes_{A} C$.

\begin{Lemma}[Descent lemma for log monoids]
\label{lem:descent:lemma:log:strcuture}
Let $f : A \to B$ be a faithfully flat homomorphism of rings.
Here we consider the descent exact sequence
\[
 0\to A \overset{f}{\longrightarrow} B 
\labelrightrightarrow{p}{q} B \otimes_A B,
\]
where $p(b) = b \otimes 1$ and $q(b) = 1 \otimes b$.
Then, we have the following:
\begin{enumerate}
\renewcommand{\labelenumi}{(\arabic{enumi})}
\item
Let  $\alpha : M \to A$ be a log monoid over $A$ such that
$M$ is integral.
Then, the sequence
\[
 0 \to M = M \boxtimes_{A} A
   \overset{\operatorname{id}\boxtimes f}{\longrightarrow}
   M \boxtimes_{A} B \labelrightrightarrow{\operatorname{id} \boxtimes p}%
{\operatorname{id} \boxtimes q} 
   M \boxtimes_{A} (B \otimes_{A} B)
\]
is exact.

\item
Let $\alpha : M \to A$ and $\beta : N \to A$ be log monoids over $A$ such that
$M$ and $N$ are integral.
Then, the sequence
\begin{multline*}
\qquad \quad 0 \to \Hom_A(M, N) \overset{F}{\longrightarrow}
\Hom_B(M \boxtimes_{A} B, N \boxtimes_{A} B)
\labelrightrightarrow{P}{Q} \\
\Hom_{B \otimes_A B}(M \boxtimes_{A} (B \otimes_A B), 
N \boxtimes_{A} (B \otimes_A B))
\end{multline*}
induced by
the descent exact sequence is exact, where
$F$, $P$ and $Q$ are given as follows:
First of all, $F(\phi) = \phi \boxtimes \operatorname{id}$.
Note that  $f$, $p$ and  $q$ are quasi-local.
Thus, using the homomorphism $B^{\times} \to (B \otimes_A B)^{\times}$
via $p$ and \rom{(5)} of Lemma~\rom{\ref{lem:pushout:monoid}},
\[
\begin{cases}
  M \boxtimes_{A} (B \otimes_A B)
  \overset{\sim}{\longrightarrow}
 (M \boxtimes_{A} B) \boxtimes_{B}
 (B \otimes_A B), \\
  N \boxtimes_{A} (B \otimes_A B)
  \overset{\sim}{\longrightarrow}
 (N \boxtimes_{A} B) \boxtimes_{B}
 (B \otimes_A B).
\end{cases}
\]
Then, $P(\phi') = \phi' \boxtimes_p 
\operatorname{id}_{(B \otimes_A B)}$.
The subscript $p$ means that we use the homomorphism
$B^{\times} \to (B \otimes_A B)^{\times}$ in terms of $p$.
In the same way, using the homomorphism
$B^{\times} \to (B \otimes_A B)^{\times}$ via $q$,
$Q(\phi') = \phi' \boxtimes_q 
\operatorname{id}_{(B \otimes_A B)}$.
\end{enumerate}
\end{Lemma}

\Proof
(1) First,
let us see that
$M \overset{\operatorname{id}\boxtimes f}{\longrightarrow} M \boxtimes_{A} B$
is injective. Indeed, we assume that $(m, 1) \sim (m', 1)$ in $M \boxtimes_{A} B$.
Then, there is $a \in A^{\times}$ with
$(m, 1)(a, f(a^{-1})) = (m', 1)$ because $f$ is quasi-local.
Thus, $f(a) = 1$. Since $f$ is faithfully flat, 
$f$ is injective. Therefore, $a = 1$.
Hence $m = m'$.

It is easy to see that
$(\operatorname{id} \boxtimes p) \circ (\operatorname{id}\boxtimes f) = 
(\operatorname{id} \boxtimes q) \circ (\operatorname{id}\boxtimes f)$.

Finally, we assume $(\operatorname{id} \boxtimes p)(m,b)
= (\operatorname{id} \boxtimes q)(m,b)$ for some $m \in M$ and 
$b \in B^{\times}$. 
We set $r = p \circ f = q \circ f$.
Then, there is $a \in A^{\times}$
with $(m, b \otimes 1)\cdot(a, r(a^{-1})) = (m, 1 \otimes b)$
because $r$ is quasi-local.
Thus, $m a  = m$ and $(b \otimes 1)r(a^{-1}) = 1 \otimes b$.
Here since $M$ is integral, $a = 1$.
Therefore, $b \otimes 1 = 1 \otimes b$. Hence
$b = f(a')$ for some $a' \in A$. Here $f$ is quasi-local.
Thus, $a' \in A^{\times}$.
Hence, $(m, b) = (\operatorname{id} \boxtimes f)(m, a)$ with $a \in A^{\times}$.

\medskip
(2) We need to show that, for 
$\phi' \in \Hom_B(M \boxtimes_{A} B, N \boxtimes_{A} B)$,
if $P(\phi') = Q(\phi')$, then there is a $\phi \in  \Hom_A(M, N)$
with $F(\phi) = \phi'$.
Since the following diagrams
\[
\xymatrix{
 M \boxtimes_{A} B
   \ar[r]^<<<<<{\operatorname{id} \boxtimes p} \ar[d]^{\phi'} & 
 M \boxtimes_{A} (B \otimes_A B) \ar[d]^{P(\phi')} \\
 N \boxtimes_{A} B
   \ar[r]^<<<<<{\operatorname{id} \boxtimes p} & 
 N \boxtimes_{A} (B \otimes_A B)
}\quad
\xymatrix{
 M \boxtimes_{A} B
   \ar[r]^<<<<<{\operatorname{id} \boxtimes q} \ar[d]^{\phi'} & 
 M \boxtimes_{A} (B \otimes_A B) \ar[d]^{Q(\phi')} \\
 N \boxtimes_{A} B
   \ar[r]^<<<<<{\operatorname{id} \boxtimes q} & 
 N \boxtimes_{A} (B \otimes_A B)
}
\]
are commutative, we have
\begin{align*}
(\operatorname{id} \boxtimes p)\circ \phi' \circ 
(\operatorname{id} \boxtimes f )(m)
& = 
P(\phi')\circ (\operatorname{id} \boxtimes p) \circ
(\operatorname{id} \boxtimes f)(m) \\
& = Q(\phi')\circ(\operatorname{id} \boxtimes q)\circ
(\operatorname{id} \boxtimes f)(m)\\
& = (\operatorname{id} \boxtimes q)\circ\phi'\circ
(\operatorname{id} \boxtimes f )(m).
\end{align*}
Therefore, there is a unique $n \in N$ such that
$\phi'\circ (\operatorname{id} \boxtimes f )(m) 
= (\operatorname{id} \boxtimes f )(n)$.
We denote this $n$ by $\phi(m)$.
It is easy to see that $\phi$ is a homomorphism and
the following diagram is commutative:
\[
 \xymatrix{
 M \ar[r]^<<<<<{\operatorname{id} \boxtimes f} \ar[d]_{\phi} & 
 M \boxtimes_{A} B \ar[d]_{\phi'}\\
 N \ar[r]^<<<<<{\operatorname{id} \boxtimes f} & 
 N \boxtimes_{A} B.
 }
\]
Finally, let us see that $\phi$ is a homomorphism over $A$.
Indeed, since $\phi'$ is a homomorphism over $B$ and
the following diagrams
\[
 \xymatrix{
 M \ar[r]^<<<<<{\operatorname{id} \boxtimes f} \ar[d]_{\alpha}& 
 M \boxtimes_{A} B \ar[d]^{\alpha_B} \\
 A \ar[r]^f & B
 }\quad
\xymatrix{
 N \ar[r]^<<<<<{\operatorname{id} \boxtimes f} \ar[d]_{\beta}& 
 N \boxtimes_{A} B \ar[d]^{\beta_B} \\
 A \ar[r]^f & B
 }
\]
are commutative,
\begin{align*}
 f \circ \beta \circ \phi & = 
 \beta_B \circ (\operatorname{id} \boxtimes f)\circ \phi=
 \beta_B \circ \phi' \circ (\operatorname{id} \boxtimes f)\\
 &  = \alpha_B \circ (\operatorname{id} \boxtimes f)
 = f \circ \alpha.
\end{align*}
Thus, $\alpha = \beta \circ \phi$ because $f$
is injective.
\QED

As a consequence of the above descent lemma, we have the following proposition:

\begin{Proposition}
\label{prop:descent:etale:monoid}
Let $\pi : X' \to X$ be a faithfully flat and quasi-compact morphism
of schemes and $X'' = X' \times_X X'$.
Let $p : X'' \to X'$ and $q : X'' \to X'$ be the projections to the first factor and
the second factor respectively.
Let $M$ and $N$ be fine
log structures on $X$, and let $r = \pi \circ p = \pi \circ q$.
Then,
\[
0 \to \Hom_{\OO_X}(M, N) \overset{\pi^*}{\longrightarrow}
\Hom_{\OO_{X'}}(\pi^*(M), \pi^*(N))
\labelrightrightarrow{p^*}{q^*}
\Hom_{\OO_{X''}}(r^*(M),  r^*(N))
\]
is exact.
\end{Proposition}

\Proof
Let us begin with the following claim.

\begin{Claim}
\label{prop:descent:etale:monoid:claim:1}
The map $\pi^* : \Hom_{\OO_X}(M, N) \to
\Hom_{\OO_{X'}}(\pi^*(M), \pi^*(N))$ is injective.
\end{Claim}

Let $\phi, \psi \in \Hom_{\OO_X}(M, N)$ with
$\pi^*(\phi) = \pi^*(\psi)$.
Since $M$ and $N$ are fine log structures in the \'{e}tale topology,
it is sufficient to see that $\phi_{\bar{x}} = \psi_{\bar{x}}$
for all $x \in X$.
Choose $x' \in X'$ with $\pi(x') = x$.
Then,  we can see
\[
 \begin{cases}
 \pi^*(M)_{\bar{x}'} = 
 M_{\bar{x}} \boxtimes_{\OO_{X,\bar{x}}} \OO_{X',\bar{x}'}, \\
 \pi^*(N)_{\bar{x}'} = 
 N_{\bar{x}} \boxtimes_{\OO_{X,\bar{x}}} \OO_{X',\bar{x}'}.
 \end{cases}
\]
Therefore, the claim follows from
Lemma~\ref{lem:descent:lemma:log:strcuture}.

\medskip
Next we consider the descent problem, namely,
if $\phi' \in \Hom_{\OO_{X'}}(\pi^*(M), \pi^*(N))$ with
$p^*(\phi') = q^*(\phi')$, then there is
$\phi \in \Hom_{\OO_{X}}(M, N)$ with $\pi^*(\phi) = \phi'$.
For this purpose, let us see the following claim:

\begin{Claim}
\label{prop:descent:etale:monoid:claim:2}
Let $\{ \rho_i : U_i \to X \}_{i \in I}$ be an \'{e}tale covering of $X$.
If the descent problem holds on $U_i$ for each $i \in I$,
then so does on $X$.
\end{Claim}

Let $\phi' \in \Hom_{\OO_{X'}}(\pi^*(M), \pi^*(N))$
with $p^*(\phi') = q^*(\phi')$.
Let $s_i : U_i \times_X U_i \to U_i$ and
$t_i : U_i \times_X U_i \to U_i$ be the projection to the first factor
and the second factor respectively.
We set the induced morphisms as follows:
\[
 \xymatrix{
 X' \ar[d]^{\pi} & X'_{U_i} \ar[l]_{\rho'_i} \ar[d]^{\pi_i}
 & X'_{U_i \times_X U_i} \ar[d]^{\pi'_i}
 \ar@<1ex>[l]^<<<<{t'_i} \ar@<-1ex>[l]_<<<<{s'_i} \\
 X  & U_i \ar[l]_{\rho_i}
 & U_i \times_X U_i \ar@<1ex>[l]^<<<<{t_i} \ar@<-1ex>[l]_<<<<{s_i}
}
\]
Then, by our assumption, for each $i \in I$,
there is $\phi_{U_i} \in  \Hom_{\OO_{U_i}}(M_{U_i}, N_{U_i})$
with $\pi_i^*(\phi_{U_i}) = {\rho'_i}^*(\phi')$.
Here
\[
{\pi'_i}^*s_i^*(\phi_{U_i}) =
{s'_i}^* \pi_i^*(\phi_{U_i}) = {s'_i}^* {\rho'_i}^*(\phi')
= {t'_i}^* {\rho'_i}^*(\phi')
= {t'_i}^* \pi_i^*(\phi_{U_i}) = {\pi'_i}^* t_i^*(\phi_{U_i}).  
\]
Thus, by the previous claim, $s_i^*(\phi_{U_i}) = t_i^*(\phi_{U_i})$
for each $i \in I$.
Since $M$ and $N$ are sheaves on the \'{e}tale topology,
there is $\phi_i \in \Hom_{\OO_{V_i}}(M_{V_i}, N_{V_i})$
with $\rho_i^*(\phi_i) = \phi_{U_i}$,
where $V_i = \rho_i(U_i)$.
Moreover, by the previous claim,
$\rest{\phi_{i}}{V_i \cap V_j} = \rest{\phi_{j}}{V_i \cap V_j}$.
Thus, there is $\phi \in \Hom_{\OO_X}(M, N)$ such that
$\rest{\phi}{V_i} = \phi_i$ for all $i \in I$.

\medskip
By the above claims, we may assume that
$X = \Spec(A)$ for some ring $A$.
Thus, using the quasi-compactness of $\pi$ and
the standard techniques as in the case of modules,
we may assume that $X' = \Spec(B)$ for some ring $B$
faithfully flat over $A$.
Therefore, by Lemma~\ref{lem:descent:lemma:log:strcuture},
we have our assertion as in Claim~\ref{prop:descent:etale:monoid:claim:1}.
\QED

\begin{Theorem}
\label{thm:descent:log:smooth}
Let $f : X \to S$ and $g : Y \to S$ be semistable schemes over a locally
noetherian scheme $S$, and let $\phi : X \to Y$ be a morphism over $S$.
Let $M_X$, $M_Y$ and $M_S$ be fine log structures on $X$, $Y$ and $S$
respectively.
We assume that $(X, M_X)$ and $(Y, M_Y)$ are log smooth and integral over
$(S, M_S)$ and $\phi$ is admissible with respect to
$M_Y/M_S$.
Then, we have the following:
\begin{enumerate}
\renewcommand{\labelenumi}{(\arabic{enumi})}
\item
Let $\pi : X' \to X$ be a surjective and smooth morphism.
If there is a log morphism
$(\phi \circ \pi, h') : (X', \pi^*(M_X)) \to (Y, M_Y)$
over $(S, M_S)$, then it descends to a log morphism
$(\phi, h) : (X, M_X) \to (Y, M_Y)$ over $(S, M_S)$.

\item
Let $\pi : S' \to S$ be a faithfully flat and quasi-compact morphism.
Let $X' = X \times_S S'$, $Y' = Y \times_S  S'$ and 
$\phi' = \phi \times_S \operatorname{id}_{S'}$. We set the induced morphisms as
follows:
\[
\begin{CD}
X @<{\pi_X}<< X' \\
@V{f}VV @VV{f'}V \\
S @<{\pi}<< S'
\end{CD}
\qquad\qquad
\begin{CD}
Y @<{\pi_Y}<< Y' \\
@V{g}VV @VV{g'}V \\
S @<{\pi}<< S'
\end{CD}
\]
If there is a log morphism
$(\phi', h') : (X', \pi_X^*(M_X)) \to (Y', \pi_Y^*(M_Y))$
over $(S', \pi^*(M_S))$, then it descends to
a log morphism
$(\phi, h) : (X, M_X) \to (Y, M_Y)$ over $(S, M_S)$.
\end{enumerate}
\end{Theorem}

\Proof
In the case of (1), it is easy to see that
$f \circ \pi : X' \to S$ is a semistable scheme and
$\phi \circ \pi : X' \to Y$ is admissible with respect to $\Supp(M_Y/M_S)$
by Lemma~\ref{lem:etale:not:subset}.
Thus, (1) and (2) are consequences of Theorem~\ref{thm:rigid:log:morphism}, 
Proposition~\ref{prop:descent:etale:monoid} and
Lemma~\ref{lem:adjoint:log:monoid}.
\QED

\begin{Corollary}
\label{cor:smooth:local:const}
Let $f : X \to S$ and $g : Y \to S$ be semistable schemes over a locally
noetherian scheme $S$, and let $\phi : X \to Y$ be a morphism over $S$.
Let $M_X$, $M_Y$ and $M_S$ be fine log structures on $X$, $Y$ and $S$
respectively.
We assume that $(X, M_X)$ and $(Y, M_Y)$ are log smooth and integral over
$(S, M_S)$ and $\phi$ is admissible with respect to
$M_Y/M_S$.
Let $\{ \pi_i : X_i \to X \}_{i \in I}$ and $\{ \mu_j : Y_j \to Y \}_{j \in J}$
be families of smooth morphisms such that
$X = \bigcup_{i \in I} \pi_i(X_i)$ and $Y = \bigcup_{j \in J} \mu_j(Y_j)$.
We assume that, for each $i \in I$, there are $j \in J$ and
a log morphism $(\phi_i, h_i) : (X_i, \pi_i^*(M_X))  \to (Y_j, \mu_j^*(M_Y))$ with
the following diagram commutative:
\[
\begin{CD}
X_i @>{\phi_i}>> Y_j \\
@V{\pi_i}VV @VV{\mu_j}V \\
X @>>{\phi}> Y.
\end{CD}
\]
Then, there is a log morphism
\[
(\phi, h) : (X, M_X) \to (Y, M_Y)
\]
as an extension of $\phi : X \to Y$ such that
the diagram
\[
\begin{CD}
(X_i,  \pi_i^*(M_X))@>{(\phi_i, h_i)}>> (Y_j, \mu_j^*(M_Y))\\
@V{(\pi_i, \operatorname{nat})}VV @VV{(\mu_j,  \operatorname{nat})}V \\
(X, M_X) @>>{(\phi, h)}> (Y, M_Y)
\end{CD}
\]
is commutative, where $\operatorname{nat}$ is the natural homomorphism.
\end{Corollary}

\Proof
We set $U_i = \pi_i(X_i)$.
Then, by (1) of Theorem~\ref{thm:descent:log:smooth},
there is a log morphism
\[
(\rest{\phi}{U_i}, h'_i) : (U_i, \rest{M_X}{U_i}) \to (Y, M_Y)
\]
such that the diagram
\[
\begin{CD}
(X_i,  \pi_i^*(M_X))@>{(\phi_i, h_i)}>> (Y_j, \mu_j^*(M_Y))\\
@V{(\pi_i, \operatorname{nat})}VV @VV{(\mu_j,  \operatorname{nat})}V \\
(U_i, \rest{M_X}{U_i}) @>>{(\rest{\phi}{U_i}, h'_i)}> (Y, M_Y)
\end{CD}
\]
is commutative. By using the rigidity theorem,
we can construct 
\[
(\phi, h) : (X, M_X) \to (Y, M_Y)
\]
such that $\rest{(\phi, h)}{U_i} = (\rest{\phi}{U_i}, h'_i)$.
\QED

\section{Dualizing sheaves of log semistable schemes}

Let $(f, h) : (X, M_X) \to (Y, M_Y)$ be a morphism of fine
log schemes. We say $(f, h) : (X, M_X) \to (Y, M_Y)$ is free (or
$(X, M_X)$ is free over $(S, M_S)$) if, for any $x \in X$,
$\Coker(\overline{M}^{gr}_{Y,\bar{f(x)}} \to \overline{M}^{gr}_{X,\bar{x}})$
is a free abelian group.

\begin{Proposition}
\label{prop:hom;dualizing:to:det:log:form}
Let $S$ be a locally noetherian scheme and
$f : X \to S$ a semistable scheme over $S$.
Let $M_X$ and $M_S$ be fine log structures of 
$X$ and $S$ respectively.
We assume that $(X, M_X)$ and $(Y, M_Y)$ are free, integral and 
log smooth over $(S, M_S)$. Let $\omega_{X/S}$ be the dualizing sheaf
of $X \to S$.
Then, there is the canonical injective homomorphism
\[
\phi : \omega_{X/S} \to \det(\Omega^1_{X/S}(\log(M_X/M_S)))
\]
with the following properties:
\begin{enumerate}
\renewcommand{\labelenumi}{(\arabic{enumi})}
\item
$\phi$ is the identity on the outside of $\Supp(M_X/M_S)$.

\item
Let us take the effective Cartier divisor $\BB_{M_X/M_S}$ with
\[
\phi(\omega_{X/S}) \otimes \OO_X(\BB_{M_X/M_S}) =  \det(\Omega^1_{X/S}(\log(M_X/M_S))),
\]
that is, $x \in \Supp(\BB_{M_X/M_S})$ if and only if $\phi$ is not surjective at $x$.
Then, $\BB_{M_X/M_S}$ is flat over $S$.
\rom{(}$\BB_{M_X/M_S}$ is called the boundary divisor of $(X, M_X) \to (S, M_S)$.\rom{)}

\item
Let $x$ be a closed point of $X$ and $s = f(x)$.
We set $Q = \overline{M}_{S,\bar{s}}$ and $P = \overline{M}_{X,\bar{x}}$.
Let $h : Q \to P$ be the induced homomorphism.
Let $\alpha : M_X \to \OO_X$ be the canonical homomorphism.
Here we define $t_{\bar{x}} \in \OO_{X,\bar{x}}$ as follows:
\begin{enumerate}
\renewcommand{\labelenumii}{(\arabic{enumi}.\alph{enumii})}
\item
If $h : Q \to P$ splits, i.e.,
there is a submonoid $N$ of $P$ with $P = f(Q) \times N$, then
\[
t_{\bar{x}} =  \prod_{p \in \Irr(N)} \alpha(\tilde{p}),
\]
where  $\Irr(N)$ is the set of all irreducible elements of $N$ and
$\tilde{p}$'s are elements of $M_{X,\bar{x}}$ with 
$\tilde{p} \equiv p \mod \OO_{X,\bar{x}}^{\times}$.

\item
If $h : Q \to P$ does not split, i.e.,
$h : Q \to P$ has a semistable structure $(\sigma, q_0, \Delta, B)$
for some $\sigma \subseteq P$,  $q_0 \in Q$ and $\Delta, B \in \NN^{\sigma}$,
then
\[
t_{\bar{x}} =  \prod_{p \in \sigma \setminus 
\Supp(\Delta)} \alpha(\tilde{p}),
\]
where
$\tilde{p}$'s are elements of $M_{X,\bar{x}}$ with 
$\tilde{p} \equiv p \mod \OO_{X,\bar{x}}^{\times}$.
\end{enumerate}
Then, $\OO_{X,\bar{x}}(-\BB_{M_X/M_S}) = t_{\bar{x}} \OO_{X,\bar{x}}$.
\end{enumerate}
\end{Proposition}

\Proof
Let us begin with the following lemma:

\begin{Lemma}
\label{lem:canonical:hom:dualizing}
Let $\phi : (A, m_A) \to (B, m_B)$ be a homomorphism of noetherian local rings
such that $\phi$ is essentially of finite type and smooth.
Let $g_1, \ldots, g_r$ be a regular sequence of $B$.
We set $I = (g_1, \ldots, g_r)$ and $C = B/I$.
Note that the dualizing sheaf $\omega_{C/A}$ of $C$ over $A$
is given by
\[
\bigwedge^n \Omega_{B/A} \otimes_B \bigwedge^r (I/I^2)^{\vee} = 
\left(\bigwedge^n \Omega_{B/A} \otimes_B C\right) \otimes_C \bigwedge^r (I/I^2)^{\vee},
\]
where $(I/I^2)^{\vee}$ is the dual as a $C$-module.
Then, there is the canonical homomorphism
\[
c : \bigwedge^{n-r} \Omega_{C/A} \to \omega_{C/A}
\]
with the following properties:
Let $f_1, \ldots, f_n \in B$ such that
$df_1, \ldots, df_n$ form a free basis of $\Omega_{B/A}$.
For a subset $S = \{ j_{r+1}, \ldots, j_{n}\}$ of $\{ 1, \ldots, n \}$,
\[
c (d\bar{f}_{j_{r+1}} \wedge \cdots \wedge d\bar{f}_{j_{n}}) = \Delta_{S} \cdot
(\bar{g}_1 \wedge \cdots \wedge \bar{g}_r)^{\vee} \otimes_C (df_1 \wedge \cdots \wedge
df_n),
\]
where $\bar{f}_{j_{r+1}}, \ldots, \bar{f}_{j_n}$ are the classes of
$f_{j_{r+1}}, \ldots, f_{j_n}$ in $C$,
$\bar{g}_1, \ldots, \bar{g}_r$ are the classes of $g_1, \ldots, g_r$
in $I/I^2$ and $\Delta_S$ is given by the following equation:
\[
dg_1 \wedge \cdots \wedge dg_r \wedge df_{j_{r+1}} \wedge \cdots \wedge
df_{j_n} = \Delta_S \cdot (df_1 \wedge \cdots \wedge df_n).
\]
We say $\omega_0 = (\bar{g}_1 \wedge \cdots \wedge \bar{g}_r)^{\vee} 
\otimes_C (df_1 \wedge \cdots \wedge
df_n)$ is the basis of $\omega_{C/A}$ with respect to
$g_1, \ldots, g_r, f_1, \ldots, f_r$.
\end{Lemma}

\Proof
This lemma is essentially proved in \cite[Lemma~4.12 in Chapter~6]{LiuAlg}.
Let $N$ be a submodule of $\Omega_{B/A}$ generated by $dg_1, \ldots, dg_r$.
Note that $N \otimes_{B} C$ does not depend on the choice of regular
sequences which generate $I$. Moreover, $(\Omega_{B/A}/N) \otimes_B C =
\Omega_{C/A}$.
Here, we have a homomorphism
\[
\rho : \bigwedge^r N \otimes_B \bigwedge^{n-r} \Omega_{B/A} \to \bigwedge^n \Omega_{B/A}
\]
given by $\rho((x_1 \wedge \cdots \wedge x_r) \otimes (x_{r+1} \wedge \cdots \wedge
x_n)) = x_1 \wedge \cdots \wedge x_n$.
Since $\bigwedge^{r+1} N = 0$, the above homomorphism
induces 
\[
\bar{\rho} : \bigwedge^r N \otimes_B \bigwedge^{n-r} (\Omega_{B/A}/N) \to \bigwedge^n \Omega_{B/A}.
\]
Therefore, by tensoring $\otimes_B C$ and by composing
$I/I^2 \to N \otimes_B C$, we get
\[
\bar{\rho}_C : \bigwedge^r (I/I^2)
\otimes_C \bigwedge^{n-r} \Omega_{C/A} \to
\bigwedge^r (N \otimes_B C )
\otimes_C \bigwedge^{n-r} \Omega_{C/A}
\to \bigwedge^n \Omega_{B/A} \otimes_B C.
\]
Thus, $\bar{\rho}_C$ gives rise to the canonical homomorphism
\[
c : \bigwedge^{n-r} \Omega_{C/A} \to 
\left(\bigwedge^n \Omega_{B/A} \otimes_B C\right) 
\otimes_C \bigwedge^r (I/I^2)^{\vee} = \omega_{C/A}.
\]
Let $f_1, \ldots, f_n \in B$ such that
$df_1, \ldots, df_n$ form a free basis of $\Omega_{B/A}$, and
let $S = \{ j_{r+1}, \ldots, j_n \}$ be a subset of
$\{ 1, \ldots, n \}$.
Then, note that
\begin{multline*}
\bar{\rho}_C((\bar{g}_1 \wedge \cdots \wedge \bar{g}_r) \otimes
(d\bar{f}_{j_{r+1}} \wedge \cdots \wedge d\bar{f}_{j_{n}})) \\
= \rho((dg_1 \wedge \cdots \wedge dg_r) \otimes
(df_{j_{r+1}} \wedge \cdots \wedge df_{j_{n}})) \otimes 1\\
= (dg_1 \wedge \cdots \wedge dg_r \wedge
df_{j_{r+1}} \wedge \cdots \wedge df_{j_{n}}) \otimes 1 \\
= \Delta_S \cdot (df_1 \wedge \cdots \wedge df_n ) \otimes 1.
\end{multline*}
Thus, we get the lemma.
\QED

\bigskip
Let us go back to the proof of Proposition~\ref{prop:hom;dualizing:to:det:log:form}.
We set $X_0 = X \setminus \Supp(M_X/M_S)$.
Then, $f$ is smooth on $X_0$ and the log structure $M_X$ on $X_0$ is trivial over $M_S$.
Thus, $\rest{\omega_{X/S}} {X_0}=  \det(\Omega^1_{X_0/S})$ and
$\rest{\det(\Omega^1_{X/S}(\log(M_X/M_S)))}{X_0} = \det(\Omega^1_{X_0/S})$.
Therefore, we can take the canonical homomorphism
\[
\phi_0 : \rest{\omega_{X/S}} {X_0} \to \rest{\det(\Omega^1_{X/S}(\log(M_X/M_S)))}{X_0}
\]
as the identity map, so that our problem is how we can extend the homomorphism
$\phi_0$ on $X$. Note that an extension of $\phi$ is uniquely determined if it exists.

\begin{Claim}
\label{claim:prop:hom;dualizing:to:det:log:form:1}
Let $g : S' \to S$ be a faithfully flat morphism of schemes. We set  $X' = X \times_S S'$
and the induced morphisms as follows:
\[
\begin{CD}
X' @>{g'}>> X \\
@V{f'}VV @VV{f}V \\
S' @>{g}>> S.
\end{CD}
\]
Moreover, we set $M_{S'} = g^*(M_S)$ and $M_{X'} = {g'}^*(M_X)$. If there is an extension
\[
\phi' : \omega_{X'/S'} \to \det(\Omega^1_{X'/S'}(\log(M_{X'}/M_{S'})))
\]
of
\[
 \phi'_0 : \rest{\omega_{X'/S'}}{X'_0} \to \rest{\det(\Omega^1_{X'/S'}(\log(M_{X'}/M_{S'})))}{X'_0},
\]
then $\phi'$ descends to an extension
\[
\phi : \omega_{X/S} \to \det(\Omega^1_{X/S}(\log(M_{X}/M_{S})))
\]
over $X$, where $X'_0 = \Supp(M_{X'}/M_{S'}) = {g'}^{-1}(X_0)$.
Further, if $\BB_{M_{X'}/M_{S'}}$ is flat over $S'$, then
so is $\BB_{M_{X}/M_S}$ over $S$.
\end{Claim}

Let us consider $S'' = S' \times_S S'$ and $X'' = X' \times_X X'$. 
Let $p : S'' \to S'$ and $p' : X'' \to X'$ (resp. $q : S'' \to S'$ and $q' : X'' \to X'$) be the
projections to the first (resp. second) factor:
\[
S'' \labelrightrightarrow{p}{q} S' \overset{g}{\longrightarrow} S,
\qquad\qquad
X'' \labelrightrightarrow{p'}{q'} X' \overset{g'}{\longrightarrow} X.
\]
We set $\rho = g \circ p = g \circ q$ and
$\rho' = g' \circ p' = g' \circ q'$. 
Note that
\[
\begin{cases}
\omega_{X'/S'} = {g'}^*(\omega_{X/S}) \\
\det(\Omega^1_{X'/S'}(\log(M_{X'}/M_{S'}))) =
{g'}^*(\det(\Omega^1_{X/S}(\log(M_{X}/M_{S})))) \\
\omega_{X''/S''} = {\rho'}^*(\omega_{X/S}) \\
\det(\Omega^1_{X''/S''}(\log(M_{X''}/M_{S''}))) =
{\rho'}^*(\det(\Omega^1_{X/S}(\log(M_{X}/M_{S})))).
\end{cases}
\]
Thus, by virtue of the uniqueness of extension,
we can see that ${p'}^*(\phi') = {q'}^*(\phi')$.
Hence, by using descent theory, $\phi'$ descends to an extension
\[
\phi : \omega_{X/S} \to \det(\Omega^1_{X/S}(\log(M_{X}/M_{S})))
\]
over $X$.
Further, by the definition of boundary divisors, 
$\BB_{M_{X'}/M_{S'}} = {g'}^{-1}(\BB_{M_{X}/M_S})$.
Thus, the natural morphism $\BB_{M_{X'}/M_{S'}} \to \BB_{M_{X}/M_S}$
is faithfully flat. Therefore, if $\BB_{M_{X'}/M_{S'}}$ is flat over $S'$,
then so is $\BB_{M_{X}/M_S}$ over $S$.

\medskip
Moreover, let us consider one more similar claim:

\begin{Claim}
\label{claim:prop:hom;dualizing:to:det:log:form:2}
Let $u : U \to X$ be an \'{e}tale morphism.
We set $M_{U} = {u}^*(M_X)$. If there is an extension
\[
\phi' : \omega_{U/S} \to \det(\Omega^1_{U/S}(\log(M_{U}/M_{S})))
\]
of
\[
 \phi'_0 : \rest{\omega_{U/S}}{U_0} \to \rest{\det(\Omega^1_{U/S}(\log(M_{U}/M_{S})))}{U_0},
\]
then $\phi'$ descends to an extension
\[
\phi : \rest{\omega_{X/S}}{u(U)} \to \rest{\det(\Omega^1_{X/S}(\log(M_{X}/M_{S})))}{u(U)},
\]
where $U_0 = {u}^{-1}(X_0)$.
Further, if $\BB_{M_{U}/M_{S}}$ is flat over $S$, then
so is $\rest{\BB_{M_{X}/M_S}}{u(U)}$ over $S$.
\end{Claim}

Let $U' = U \times_X U$ and $p : U' \to U$ (resp. $q : U' \to U$) be the projection
to the first (resp. second) factor.  Note that
\[
\begin{cases}
\omega_{U/S} = {u}^*(\omega_{X/S}) \\
\det(\Omega^1_{U/S}(\log(M_{U}/M_{S}))) =
{u}^*(\det(\Omega^1_{X/S}(\log(M_{X}/M_{S})))) \\
\omega_{U'/S} = {u'}^*(\omega_{X/S}) \\
\det(\Omega^1_{U'/S}(\log(M_{U'}/M_{S}))) =
{u'}^*(\det(\Omega^1_{X/S}(\log(M_{X}/M_{S})))),
\end{cases}
\]
where $u' = u \circ p = u \circ q$ and $M_{U'} = {u'}^*(M_X)$.
Thus, using the uniqueness of extension and descent theory,
we can see our claim in the same way as in the previous claim.

\medskip
Finally, we claim the following:

\begin{Claim}
\label{claim:prop:hom;dualizing:to:det:log:form:3}
$\phi_0$ extends to $\phi : \omega_{X/S} \to \det(\Omega^1_{X/S}(\log(M_X/M_S)))$.
Moreover, if we take $t_{\bar{x}}$ as in \rom{(3)} of
Proposition~\rom{\ref{prop:hom;dualizing:to:det:log:form}} for $x \in X$,
then 
\[
\phi(\omega_{X/S})_{\bar{x}} = 
t_{\bar{x}} \det(\Omega^1_{X/S}(\log(M_X/M_S)))_{\bar{x}}
\]
and $\OO_{X,\bar{x}}/t_{\bar{x}} \OO_{X,\bar{x}}$ is flat over $\OO_{S,\bar{s}}$.
\end{Claim}

Clearly we may assume that $S = \Spec(A)$ for some noetherian local ring
$(A, m)$. 
Let $x$ be a closed point of $X$ lying over $m$.
Then, by using Proposition~\ref{prop:split:log:structure} and
[EGA~III, Chapter~0, 10.3.1] together with
Claim~\ref{claim:prop:hom;dualizing:to:det:log:form:1} and
Lemma~\ref{lem:flat:free:comp} below, we may further assume the following:
\begin{enumerate}
\renewcommand{\labelenumi}{(\roman{enumi})}
\item
There is a fine and sharp monoid $Q$ and
a homomorphism $\pi_Q : Q \to M_{S, s}$
such that $Q \to M_{S,s} \to \overline{M}_{S,\bar{s}}$ is
bijective.

\item
$k = A/m$ is algebraically closed and the natural homomorphism
$A/m$ to the residue field of $X$ at $x$ is an isomorphism.
\end{enumerate}
Then, we have a good chart of
$(X, M_X) \to (S, M_S)$ at $x$, namely, 
there are a fine and sharp monoid $P$ and
homomorphisms $\pi_P : P \to \overline{M}_{X,\bar{s}}$ and
$h : Q \to P$ such that
$P \to M_{X,\bar{x}} \to \overline{M}_{X,\bar{x}}$ is bijective,
the diagram
\[
\begin{CD}
Q @>{h}>> P \\
@V{\pi_Q}VV @VV{\pi_P}V \\
M_{S,\bar{s}} @>>> M_{X,\bar{x}}
\end{CD}
\]
is commutative and that
the natural homomorphism
\[
\OO_{S,\bar{s}} \otimes_{\OO_{S,\bar{s}}[Q]} \OO_{S,\bar{s}}[P]
\to \OO_{X,\bar{x}}
\]
is smooth.
Note that it is sufficient to construct the extension of $\phi_0$ over $\OO_{X,\bar{x}}$
by Claim~\ref{claim:prop:hom;dualizing:to:det:log:form:2}.

Here we need to take care of the following two cases because
$(X, M_X)$ is free over $(S, M_S)$:
\begin{enumerate}
\item[(A)] $f$ is smooth at $x$.

\item[(B)] $f$ is not smooth at $x$ and $h : Q \to P$ does not split.
\end{enumerate}

For the case (A), there is a submonoid $N$ of $P$ such that
$P = h(Q) \times N$ and 
$N$ is isomorphic to $\NN^a$ for some non-negative integer $a$.
Then,
\[
\OO_{S,\bar{s}} \otimes_{\OO_{S,\bar{s}}[Q]} \OO_{S,\bar{s}}[P] \simeq
\OO_{S,\bar{s}}[\NN^a] = \OO_{S,\bar{s}}[X_1, \ldots, X_a].
\]
Thus, adding more indeterminates $X_{a+1}, \ldots, X_n$,
$\OO_{X,\bar{x}}$ is \'{e}tale over
\[
\OO_{S,\bar{s}}[X_1, \ldots, X_n].
\]
Therefore, it is sufficient to see our assertion on 
$X' = \Spec(\OO_{S,\bar{s}}[X_1, \ldots, X_n])$ around
the origin $o = (m_{S,\bar{s}}, X_1, \ldots, X_n)$.
Then, the bases of 
\[
\omega_{X'/S, \bar{o}}
\quad\text{and}\quad
\det(\Omega^1_{X'/S}(\log(M_{X'}/M_S))_{\bar{o}}
\]
are
\[
d X_1 \wedge \cdots \wedge d X_n
\quad\text{and}\quad
\frac{d X_1}{X_1} \wedge \cdots \wedge \frac{d X_a}{X_a}
\wedge d X_{a+1} \wedge \cdots \wedge d X_{n}
\]
respectively.
Hence, we get the first and second assertions.
Moreover, in this case, $t_{\bar{x}} = X_1 \cdots X_a$.
Thus, the last assertion follows from Remark~\ref{rem:flatness:x:1:x:l}.

\medskip
Finally, we consider the case (B).
Here we set $\sigma = \{ p_1, \ldots, p_r \}$ with
$\Supp(\Delta) = \{ p_1, \ldots, p_l \}$.
Then, $\OO_{X,\bar{x}}$ is \'{e}tale over
\[
O_{S,\bar{s}}[ X_1, \ldots, X_l, X_{l+1}, \ldots, X_r, X_{r+1}. \ldots, X_n]/
(X_1 \cdots X_l - a X_{l+1}^{b_{l+1}} \cdots X_{r}^{b_r}),
\]
where the class of $X_i$ ($i=1, \ldots, l$) is $\alpha(\pi_{P}(p_i))$,
$b_j = B(p_j)$ ($j=l+1, \ldots, r$) and
$a = \alpha(\pi_{Q}(q_0))$.
We denote the class of $X_i$ by $x_i$. Moreover,
the polynomial $X_1 \cdots X_l - a X_{l+1}^{b_{l+1}} \cdots X_{r}^{b_r}$ is
denoted by $F$.
As in the case (A), it is sufficient to see our assertion on 
$X' = \Spec(\OO_{S,\bar{s}}[X_1, \ldots, X_n]/(F))$ around
the origin $o = (m_{S,\bar{s}}, X_1, \ldots, X_n)$.
Here we set
\[
\omega_i = \begin{cases}
{\displaystyle d\log(p_1) \wedge \cdots \wedge \widehat{d\log(p_i)} \wedge
\cdots \wedge d\log(p_r) \wedge dx_{r+1} \wedge \cdots \wedge dx_n}  &
\text{if $1 \leq i \leq r$} \\
{\displaystyle d\log(p_1) \wedge \cdots \wedge d\log(p_r) 
\wedge dx_{r+1} \wedge \cdots \wedge \widehat{dx_i} \wedge \cdots \wedge dx_n} &
\text{if $r < i \leq n$}
\end{cases}
\]
as an element of $\det(\Omega^1_{X'/S}(\log(M_{X'}/M_S))$.
Then, it is easy to see that
\[
\omega_i = \begin{cases}
(-1)^{i-1} \omega_1 & \text{if $1 \leq i \leq l$} \\
(-1)^{i} b_i \omega_1 & \text{if $l < i \leq r$} \\
0 & \text{if $r < i \leq n$}
\end{cases}
\]
because
\[
d\log(p_1) + \cdots + d\log(p_l) = b_{l+1} d\log(p_{l+1}) + \cdots + b_r d\log(p_r).
\]
Let
\[
\lambda = (\overline{F})^{\vee} \otimes \rest{d X_1 \wedge \cdots \wedge dX_n}{X_0}
\in ((F)/(F^2))^{\vee} \otimes \rest{\det(\Omega^1_{\AAA_S/S})}{X'}
\]
be the basis of $\omega_{X'/S}$ with respect to $F, X_1, \ldots, X_n$ as in 
Lemma~\ref{lem:canonical:hom:dualizing}, and let
\[
\phi : \omega_{X'/S} \to  \det(\Omega^1_{X'/S}(\log(M_{X'}/M_S))
\]
be a homomorphism given by
$\phi(\lambda) = x_{l+1} \cdots x_{r} \cdot \omega_1$.
Moreover, let 
\[
c : \bigwedge^{n-1}  \Omega^1_{X'/S}
\to \omega_{X'/S}
\]
be the canonical homomorphism described in Lemma~\ref{lem:canonical:hom:dualizing} and 
let
\[
c' : \bigwedge^{n-1}  \Omega^1_{X'/S} \to \det(\Omega^1_{X'/S}(\log(M_{X'}/M_S))
\]
be the natural homomorphism. 
More precisely, $c$ is given  by
\[
c(dx_1 \wedge \cdots \wedge \widehat{dx_i} \wedge \cdots \wedge dx_n)
= \Delta_i \cdot \lambda,
\]
where
$\Delta_i$ is determined by 
\[
dF \wedge d X_1 \wedge \cdots \wedge \widehat{dX_i} \wedge \cdots \wedge dX_n
= \Delta_i \cdot dX_1 \wedge \cdots \wedge dX_n.
\]
Thus, $\Delta_i = (-1)^{i-1} (\partial F/\partial X_i)$.
In order to see our assertion, it is sufficient to show that
the following diagram
\[
\xymatrix{
& \det(\Omega^1_{X'/S}(\log(M_{X'}/M_S)) \\
\bigwedge^{n-1}  \Omega^1_{X'/S} \ar[ur]^{c'} \ar[dr]^{c} &  \\
& \omega_{X'/S} \ar[uu]_{\phi}
}
\]
is commutative. Indeed,
$\bigwedge^{n-1}  \Omega^1_{X'/S}$ is generated by
$\{ dx_1 \wedge \cdots \wedge \widehat{d x_i} \wedge \cdots d x_n \}_{i=1}^n$.
Thus, we need to see that
\[
\phi(c(dx_1 \wedge \cdots \wedge \widehat{d x_i} \wedge \cdots d x_n))
= x_1 \cdots \widehat{x_i} \cdots x_r \cdot \omega_i
\]
for all $i$.
In the case where $1 \leq i \leq l$,
\begin{align*}
\beta(c(dx_1 \wedge \cdots \wedge \widehat{d x_i} \wedge \cdots d x_n)) &
= \beta((-1)^{i-1} (\partial F/\partial X_i) \lambda) \\
& = (-1)^{i-1} x_1 \cdots \widehat{x_i} \cdots x_l \cdot x_{l+1} \cdots x_r \cdot \omega_1 \\
& = x_1 \cdots \widehat{x_i} \cdots x_r \cdot \omega_i.
\end{align*}
In the case where $l < j \leq r$,
\begin{align*}
\beta(c(dx_1 \wedge \cdots \wedge \widehat{d x_i} \wedge \cdots d x_n)) &
= \beta((-1)^{i-1} (\partial F/\partial X_i) \lambda) \\
& = (-1)^{i} b_i a x_{l+1}^{b_{l+1}} \cdots x_i^{b_i-1} \cdots x_r^{b_r} \cdot
x_{l+1} \cdots x_r \cdot \omega_1 \\
& = a x_{l+1}^{b_{l+1}} \cdots x_i^{b_i} \cdots x_r^{b_r} \cdot
x_{l+1} \cdots \widehat{x_i} \cdots x_r \cdot \omega_i \\
& = x_1 \cdots \widehat{x_i} \cdots x_r \cdot \omega_i.
\end{align*}
Finally,  in the case where $j > r$,
\[
\beta(c(dx_1 \wedge \cdots \wedge \widehat{d x_i} \wedge \cdots d x_n))
= \beta((-1)^{i-1} (\partial F/\partial X_i) \lambda) = 0 = x_1 \cdots \widehat{x_i} \cdots x_r \cdot \omega_i .
\]

Moreover, in this case, $t_{\bar{x}} = x_{l+1} \cdots x_r$.
Therefore, $\OO_{X,\bar{x}}/t_{\bar{x}}\OO_{X,\bar{x}}$ is \'{e}tale over a ring
\[
R = \OO_{S,\bar{s}}[X_1, \ldots, X_n]/(X_1 \cdots X_l - a X_{l+1}^{b_{l+1}} \cdots X_r^{b_r}, X_{l+1} \cdots X_r).
\]
Note that if we set $D = \OO_{S,\bar{s}}[X_{l+1}, \ldots, X_n]/(X_{l+1} \cdots X_r)$, then
\[
R \simeq D[X_1, \ldots, X_l]/(X_1 \cdots X_l - a x_{l+1}^{b_{l+1}} \cdots x_r^{b_r}),
\]
where $x_{l+1}, \ldots, x_r$ are the classes of $X_{l+1}, \ldots, X_{r}$ in
$D$. By using Remark~\ref{rem:flatness:x:1:x:l},
$D$ is flat over $\OO_{S,\bar{s}}$ and $R$ is flat over $D$. Thus, $R$ is flat over $\OO_{S,\bar{s}}$.
\QED

\begin{Lemma}
\label{lem:flat:free:comp}
Let $A$ be a  ring,
$M$ an $A$-module, and $N_1$ and $N_2$ $A$-submodules of $M$.
Let $f : A \to B$ be a ring homomorphism
such that $B$ is faithfully flat over $A$.
If $N_1 \otimes_A B = N_2 \otimes_A B$ as $B$-submodules of $M \otimes_A B$,
then $N_1 = N_2$.
\end{Lemma}

\Proof
For each $i = 1, 2$,  let us consider a sequence
\addtocounter{Claim}{1}
\begin{equation}
\label{lem:flat:free:comp:eqn:1}
 0 \to N_i \to N_1 + N_2 \to 0,
\end{equation}
which gives rise to
\addtocounter{Claim}{1}
\begin{equation}
\label{lem:flat:free:comp:eqn:2}
0 \to N_i \otimes_A B \to (N_1 + N_2) \otimes_A B \to 0.
\end{equation}
Here, since $N_1 \otimes_A B = N_2 \otimes_A B$, we have
$(N_1 + N_2) \otimes_A = N_i \otimes_A B$.
Thus, \eqref{lem:flat:free:comp:eqn:2} is exact,
which implies that  so is \eqref{lem:flat:free:comp:eqn:1}.
Thus, $N_1 = N_2 = N_1 + N_2$.
\QED

\section{Extension of log morphisms}

Let $S$ be a locally noetherian scheme and
$f : X \to S$ a semistable scheme over $S$.
Let $M_X$ and $M_S$ be fine log structures of 
$X$ and $S$ respectively.
We assume that $(X, M_X)$ are free, integral and 
log smooth over $(S, M_S)$.
Note that, for a point $x \in \Sing(f)$ and $s = f(x)$,
$\overline{M}_{S,\bar{s}} \to \overline{M}_{X,\bar{x}}$ does not split.
Thus,
it has the unique semistable structure $(\sigma, q_0, \Delta, B)$,
where $\sigma$ is the set of irreducible elements of
$\overline{M}_{X,\bar{x}}$ not lying in the image of
$\overline{M}_{S,\bar{s}}$, $q_0 \in \overline{M}_{S,\bar{s}} $ and
$\Delta, B \in \NN^{\sigma}$.
Recall that
this $q_0$ is called the marking of
$\overline{M}_{S,\bar{s}} \to \overline{M}_{X,\bar{x}}$
(cf. Definition~\ref{def:marking:semistable:structure}).
We dente this by
$\marking_{M_X/M_S}(x)$.
Thus, we have a map
\[
\marking_{M_X/M_S} : \Sing(f) \to \coprod_{s \in S} \overline{M}_{S,\bar{s}}.
\]
If we fix $s \in S$, then we get
\[
\rest{\marking_{M_X/M_S}}{X_{\bar{s}}} \to  \overline{M}_{S,\bar{s}},
\]
where $X_{\bar{s}}$ is the geometric fiber over $s$.
Note that $\rest{\marking_{M_X/M_S}}{X_{\bar{s}}}$ is locally constant
(cf Proposition~\ref{prop:quotient:semistable}).
We say $\marking_{M_X/M_S}$ is {\em regular} if
$\marking_{M_X/M_S}(x)$ is regular for every $x \in \Sing(X_{\bar{s}})$.

\begin{Theorem}
\label{thm:extension:log:mor}
Let $(A, tA)$ be a discrete valuation ring and
$f : X \to \Spec(A)$ a generically smooth semistable scheme over $S = \Spec(A)$.
Let $M_S$ be a fine log structure on $S$ and let
$M_X$ and $M'_X$ be fine log structures on $X$.
Let
$(f, h) : (X, M_X) \to (S, M_S)$ and $(f, h') : (X, M'_X) \to (S, M_S)$
be free, smooth and integral morphisms.
If $\BB_{M_X/M_S} = \BB_{M'_X/M_S}$,
$\marking_{M_X/M_S} = \marking_{M'_X/M_S}$ and
$\marking_{M_X/M_S}$ is regular,
then there is
an isomorphism
\[
(\operatorname{id}, h) :
(X, M'_{X}) \overset{\sim}{\longrightarrow} (X, M_{X})
\]
over $(S, M_{S})$.
\end{Theorem}

\Proof
Before starting the proof of Theorem~\ref{thm:extension:log:mor},
we need several preparations.
Let $(A, m)$ be a noetherian regular local ring and
$R$ the ring of formal power series of $n$-variables
over $A$, i.e.,
$R = A \lformal X_1, \ldots, X_n \rformal$.
An element of $f = \sum_{I} a_I X^I$ of $R$ is said to be primitive if
there is no prime element $p$ of $A$ such that
$p \mid a_I$ for all $I$.
Then, we have the following:

\begin{Lemma}
\label{lem:prime:X:a:UFD}
\begin{enumerate}
\renewcommand{\labelenumi}{(\arabic{enumi})}
\item
For a non-zero $f \in R$, there are $a \in A$ and $f' \in R$
such that $f = af'$ and $f'$ is primitive.

\item
If $f$ and $g$ are primitive elements of $R$,
then so is $fg$.

\item
Let $a$ and $b$ be non-zero elements of $A$, and
let $f$ and $g$ be primitive elements of $R$.
If $af = bg$, then there is $u \in A^{\times}$ with $b = ua$.

\item
If $a$ is a non-zero element of $m$, then
$X_1 \cdots X_n - a$ is a prime element of $R$.
\end{enumerate}
\end{Lemma}

\Proof
(1) For $g = \sum_{I \in \NN^n} g_I X^I \in R$, the ideal of $A$ generated by
$\{ g_I \}_{I \in \NN^n}$ is denoted by $I(g)$.
If $g \not= 0$ and there is a prime element $p \in A$ with
$g = p g'$ for some $g' \in R$, then
$I(g) = pI(g')$. In particular, $I(g) \subsetneq I(g')$.
Thus, we have our assertion because $A$ is noetherian.

(2) We assume that there is a prime element $p$ such that
$p$ divides all coefficient of $fg$.
Then, $fg = 0$ in $(A/pA)\lformal X_1, \ldots, X_n \rformal$.
Thus, either $f = 0$ or $g = 0$ in
$(A/pA)\lformal X_1, \ldots, X_n \rformal$ because
$(A/pA)\lformal X_1, \ldots, X_n \rformal$ is an integral domain.
This is a contradiction.

(3) It is sufficient to see that if $a$ is divisible by a prime element $p$,
then so is $b$ by $p$.
Indeed, since $p$ divides $a$, 
$af = 0$ in $(A/pA)\lformal X_1, \ldots, X_n \rformal$.
Thus, $bg = 0$ in $(A/pA)\lformal X_1, \ldots, X_n \rformal$.
We set $g = \sum_I g_I X^I$. Then, since $g$ is primitive,
there is $I$ such that $g_I \not= 0$ in $A/pA$.
Then, $bg_I = 0$ in $A/pA$, which implies that $b = 0$ in $A/pA$.

(4) We assume that there is a decomposition
\[
X_1 \cdots X_n - a = fg
\]
with $f, g \not\in R^{\times}$.
We consider the above decomposition in $(A/m)\lformal X_1, \ldots, X_n \rformal$.
Then, renumbering $X_1, \ldots, X_n$ and
replacing $f$ by $(\text{unit}) \cdot f$, we may set
\[
\bar{f} = X_1 \cdots X_{i}
\quad\text{and}\quad
\bar{g} = X_{i+1} \cdots X_n,
\]
where $\bar{f}$ and $\bar{g}$ are the classes of $f$ and $g$
in $(A/m)\lformal X_1, \ldots, X_n \rformal$.
If $i = n$, then $g$ is a unit. Thus, $i < n$.
Here we can find $h, k \in m\lformal X_1, \ldots, X_n \rformal$
with $f = X_1 \cdots X_{i} + h$ and
$g = X_{i+1} \cdots X_n + k$.
Therefore,
\[
X_1 \cdots X_n - a  = (X_1 \cdots X_{i} + h)(X_{i+1} \cdots X_n + k).
\]
Thus,
\[
- a = (X_1 \cdots X_{i} + h(X_1, \ldots, X_{n-1}, 0))k(X_1, \ldots, X_{n-1}, 0).
\]
Here, by (1), 
there is $b \in A$ and $k' \in A \lformal X_1, \ldots, X_{n-1}\rformal$
such that 
\[
k(X_1, \ldots, X_{n-1}, 0) = b k'
\]
and
$k'$ is primitive.
Thus, by (2) and (3), we have $u \in A^{\times}$ with $b = -au$.
Hence,
\[ 
1 = (X_1 \cdots X_{i} + h(X_1, \ldots, X_{n-1}, 0))k'u.
\]
Therefore, $1 = h(0, \ldots, 0)k'(0, \ldots, 0)u$, which means that
$h(0, \ldots, 0) \in A^{\times}$. This is a contradiction
because $h \in m\lformal X_1, \ldots, X_n \rformal$.
\QED

Next, we consider primary decompositions in the special rings.

\begin{Proposition}
\label{prop:primary:decomp}
Let $(A, tA)$ be a discrete valuation ring and
\[
R = A \lformal X_1, \ldots, X_l, Y_1, \ldots, Y_n  \rformal/(X^{\Delta} - u t^aY^B),
\]
where $u$ is a unit, $a$ is a positive integer,
$\Delta = (1, \ldots, 1) \in \NN^l$ and
$B \in \NN^n$.
Let $I$ be an element of $\NN^n$ such that
$I(j) \leq B(j)$ for all $j=1, \ldots, n$.
We denote the class of $X_i$ and the class of $Y_j$ in $R$ by $x_i$
and $y_j$ respectively.
Then, we have the following:
\begin{enumerate}
\renewcommand{\labelenumi}{(\arabic{enumi})}
\item
$R$ is an integral domain.

\item
The ideals $(x_i,t)$ \rom{(}$1\leq i \leq l$\rom{)} and
$(x_i, y_j)$ \rom{(}$1\leq i \leq l$,  $j \in \Supp(B)$\rom{)} of $R$
are distinct prime ideals such that
$\dim R/(x_i, t) = \dim R/(x_i, y_j) = l + n -1$
for all $i=1, \ldots, l$ and all $j \in \Supp(B)$.

\item
$(t^a y^{I}) = (x_1, t^a y^{I}) \cap \cdots \cap (x_l, t^a y^{I})$.

\item
 For each $i=1, \ldots, l$, 
 \[
 (x_i,  t^a y^{I}) = (x_i, t^a) \cap 
\bigcap_{j \in \Supp(I)}(x_i, y_j^{I(j)})
\]
is a primary decomposition
with no embedded primes.
Moreover, $\sqrt{(x_i, t^a)} = (x_i, t)$ and $\sqrt{(x_i, y_j^{I(j)})} = (x_i, y_j)$.

\item
Let $j_1, \ldots, j_r$ be distinct elements of $\{ 1, \ldots, n \}$.
Then, 
\[
(y_{j_1} \cdots y_{j_r}) = (y_{j_1}) \cap \cdots \cap (y_{j_r}).
\]

\item
If $j \not\in \Supp(B)$, then $y_j$ is a prime element of $R$.
Moreover, if $j \in \Supp(B)$,
then $(y_j) = (x_1, y_j) \cap \cdots \cap (x_l, y_j)$
is a primary decomposition with no embedded primes.
\end{enumerate}
\end{Proposition}

\Proof
(1) This is a consequence of Lemma~\ref{lem:prime:X:a:UFD}.

\medskip
(2)
It is easy to see these facts by the following canonical isomorphisms:
\begin{align*}
R/(x_i,t) & \simeq (A/tA) \lformal X_1, \ldots, \widehat{X_i}, \ldots, X_l, Y \rformal. \\
R/(x_i, y_j) & \simeq A \lformal X_1, \ldots, \widehat{X_i}, \ldots, X_l, 
Y_1, \ldots, \widehat{Y_j}, \ldots, Y_n \rformal \qquad (j \in \Supp(B)).
\end{align*} 

\medskip
(3)
Let us begin with the following claim:

\begin{Claim}
Let $f$ be an element of  $A \lformal X, Y \rformal$, and let
$I \in \NN^n$. If
\[
\rest{f}{X_i= 0} = f(X_1, \ldots, X_{i-1}, 0, X_{i+1}, \ldots, X_l, Y)
\]
is divisible by $t^aY^I$ for
every $i=1, \ldots, l$, then
there are $g, h \in A \lformal X, Y \rformal$ with
$f = t^a Y^I g + X^{\Delta}h$.
\end{Claim}

We set
\[
f = f_1 + X_1 f_2 + X_1 X_2 f_3 + \cdots + X_1 \cdots X_{l-1} f_l + X_1 \cdots X_l h,
\]
where $h \in A \lformal X, Y \rformal$ and
$f_i \in A \lformal X_1, \ldots, X_{i-1}, X_{i+1}, \ldots, X_l, Y \rformal$.
Here we see  $t^a Y^I \mid f_i$ for every $i=1, \ldots, l$ by induction on $i$.
In the case where $i=1$,
$f_1 = f(0, X_2, \ldots, X_l, Y)$. Thus, $t^a Y^I \mid f_1$. In general,
\[
\rest{f}{X_i=0} =
 \rest{f_1}{X_i= 0}  + X_1 \rest{f_2}{X_i= 0} + \cdots + X_1 \cdots X_{i-1} f_i.
 \]
Thus, using hypothesis of induction, we can see that
$t^a Y^I \mid f_i$.
Therefore, we get the claim.

\medskip
Clearly, $(t^a y^{I}) \subseteq (x_1, t^ay^{I}) \cap \cdots \cap (x_l, t^a y^{I})$.
We assume 
\[
\bar{f} \in (x_1, t^ay^{I}) \cap \cdots \cap (x_l, t^a y^{I}),
\]
where
$f \in A \lformal X, Y \rformal$.
Since $\bar{f} \in (x_i, t^ay^{I})$, there are $v, w, z \in A \lformal X_1, \ldots, X_n \rformal$
with $f = X_i v + t^a Y^{I} w + (X^{\Delta} - u t^aY^B)z$.
Thus,
\[
\rest{f}{X_i= 0} = \left( \rest{w}{X_i=0}-\rest{(uz)}{X_i= 0}Y^{B-I} \right) t^a Y^{I}.
\]
Thus, by the above claim, there are $g, h \in A \lformal X_1, \ldots, X_n \rformal$
with $f = t^a Y^{I} g + X^{\Delta} h$.
Therefore, $\bar{f} \in (t^a y^{I})$.

\medskip
(4) 
The decomposition
 \[
 (x_i,  t^a y^{I}) = (x_i, t^a) \cap 
\bigcap_{j \in \Supp(I)}(x_i, y_j^{I(j)})
\]
is obvious because
\[
R/(x_i, t^a y^{I}) \simeq 
A \lformal X_1, \ldots, \widehat{X_i}, \ldots, X_n, Y \rformal/(t^a Y^{I}).
\]
We need to show that
$(x_i, t^a)$ and $(x_i, y_j^{I(j)})$ are primary ideals and
$\sqrt{(x_i, t^a)} = (x_i,t)$ and $\sqrt{(x_i, y_j^{I(j)})} = (x_i, y_j)$.
Let us see the following claim:

\begin{Claim}
Let $C$ be a noetherian ring.
If $C$ satisfies the property:
\[
\text{$xy = 0$ and $y \not\in \sqrt{0}$}
\quad\Longrightarrow\quad
\text{$x = 0$},
\]
then $C \lformal T \rformal$ holds the same property and
$\sqrt{0_{C \lformal T \rformal}} = \sqrt{0_C}\lformal T \rformal$.
\end{Claim}

Since $C$ is noetherian, there is a positive integer $N$
such that $\sqrt{0_C}^N = 0$.
Thus, it is easy to see that $\sqrt{0_{C \lformal T \rformal}} = \sqrt{0_C}\lformal T \rformal$.
We assume that $f g = 0$ and $g \not\in \sqrt{0_{C \lformal T \rformal}}$.
We set $g = \sum_{i} g_i T^i$.
If $a_i \in \sqrt{0_C}$ for all $i$, then $g \in \sqrt{0_{C \lformal T \rformal}}$.
Thus, there is $i$ such that $g_i \not\in \sqrt{0_C}$ and
$g_0, \ldots, g_{i-1} \in \sqrt{0_C}$.
Here we set $g_1 = g_0 + \cdots + g_{i-1}T^{i-1}$. 
Then $g = g_1 + T^i g_2$ and $g_2(0) \not\in \sqrt{0_C}$.
Note that
\[
f(g - g_1)^N = \sum_{j=0}^N(-1)^{N-i} \binom{N}{i} f g^i g_1^{N-i} =  (-1)^N f g_1^N = 0.
\]
Hence, $f g_2^N = 0$. By this observation, we may assume that
$g_0 \not\in \sqrt{0_C}$.
We set $f = \sum_i f_i T^i$.
Then, $fg  = f_0 g_0 + \cdots = 0$. Thus, $f_0 = 0$.
Hence $fg = f_1 g_0 T + \cdots  = 0$. Therefore, $f_1 = 0$.
In the same way, we can see $f_i = 0$ for all $i$.

\medskip
First of all,
\[
R/ (x_i,  t^a) = (A/t^aA) \lformal X_1, \ldots, \widehat{X_{i}}, \ldots, X_l, Y \rformal.
\]
Thus, by the previous claim, $(x_i, t^a)$ is a primary ideal and
\[
\sqrt{0_{R/(x_i, t^a)}} = (tA/t^aA) \lformal X_1, \ldots, \widehat{X_{i}}, \ldots, X_l, Y \rformal,
\]
i.e., $\sqrt{(x_i, t^a)} = (x_i, t )$.
Moreover, if $I(j) \not= 0$, then
\[
R/ (x_i, y_j^{I(j)}) = 
A \lformal X_1, \ldots, \widehat{X_{i}}, \ldots, X_l, Y \rformal / (Y_j^{I(j)}).
\]
Therefore, $(x_i, y_j^{I(j)})$ is a primary and $\sqrt{(x_i, y_j^{I(j)})} =
(x_i, y_j)$.

\medskip
(5)
First, we assume that $A$ is complete.
Let $\pi : A \lformal X, Y \rformal \to R$
be the canonical homomorphism.
Here we consider $\Gamma
= \{ C \in \NN^l \mid \Supp(\Delta) \not\subseteq \Supp(C) \}$.
Then, by Lemma~\ref{lem:flatness:x:1:x:l}, the map
\[
\left\{ \sum_{C \in \Gamma, D \in \NN^n } a_{C, D} X^C Y^D \mid a_{C,D} \in A \right\}
\overset{\pi}{\longrightarrow} R
\]
is bijective. 
In order to see out assertion, it is sufficient to see that
if $y_{j_a}$ divides $y_{j_1} \cdots y_{j_{a-1}} f$, then
$y_{j_a}$ divides $f$.
We set $f = \sum_{C \in \Gamma, D \in \NN^n} a_{C,D} x^C y^D$ and
$y_{j_1} \cdots y_{j_{a-1}} f = y_{j_a} \sum_{C \in \Gamma, D \in \NN^n} b_{C,D} x^Cy^D$.
Then,
\[
\sum_{C\in \Gamma, D \in \NN^n} a_{C,D} x^{C} y^{D + e_{j_1} + \cdots + e_{j_{a-1}}}=
\sum_{C \in \Gamma, D \in \NN^n} b_{C,D} x^{C } y^{D+ e_{j_a}}.
\]
Thus,
\[
\sum_{C\in \Gamma, D \in \NN^n} a_{C,D} X^{C} Y^{D + e_{j_1} + \cdots + e_{j_{a-1}}}=
\sum_{C \in \Gamma, D \in \NN^n} b_{C,D} X^{C } Y^{D+ e_{j_a}}
\]
in $A \lformal X, Y \rformal$.
Therefore, if $D(j_a) = 0$, then $a_{C, D} = 0$.
Hence,
\[
f = \sum_{C\in \Gamma, D \in \NN^n} a_{C, D + e_{j_a}} x^{C} y^{D + e_{j_a}}=
y_{j_a} \sum_{C\in \Gamma, D \in \NN^n} a_{C, D + e_{j_a}} x^{C} y^{D}.
\]
Thus, we get (5) in the case where $A$ is complete.

In general, let $\widehat{A}$ be the completion of $A$.
We set
\[
R' = \widehat{A} \lformal X_1, \ldots, X_l, Y_1, \ldots, Y_n  \rformal/(X^{\Delta} - u t^aY^B),
\]
Then, $R'$ is faithfully flat over $R'$.
By the previous observation,
\[
y_{j_1} \cdots y_{j_r} R' = y_{j_1} R' \cap \cdots \cap y_{j_r}R'.
\]
Therefore, using the faithfully flatness of $R'$ over $R$,
\begin{multline*}
y_{j_1} \cdots y_{j_r} R = (y_{j_1} \cdots y_{j_r} R') \cap R 
= (y_{j_1} R' \cap \cdots \cap y_{j_r}R') \cap R \\
= ((y_{j_1} R') \cap R) \cap \cdots \cap ((y_{j_r}R') \cap R) =
y_{j_1} R \cap \cdots \cap y_{j_r}R.
\end{multline*}

\medskip
(6)
We assume that $j \not\in \Supp(B)$.
Since 
\[
R/(y_j) \simeq A \lformal X, Y_1, \ldots, \widehat{Y_j}, \ldots, Y_n\rformal/
(X^{\Delta} - (\rest{u}{X_j=0}) t^a Y^B),
\]
this is an integral domain by (1).

Next we assume that $j \in \Supp(B)$. Then
\[
R/(y_j)\simeq A \lformal X, Y_1, \ldots, \widehat{Y_{j}}, \ldots, Y_n \rformal/
(X^{\Delta}).
\]
Thus, we get our assertion.
\QED

Finally, we consider the following proposition.

\begin{Proposition}
\label{prop:two:isom:unit}
Let $(A, tA)$ be a discrete valuation ring and $R$ an $A$-algebra.
We assume that we have two isomorphisms
\[
\phi : C = A \lformal X_1, \ldots, X_l, Y_1, \ldots, Y_n  \rformal/(X^{\Delta} - u t^aY^B)
\overset{\sim}{\longrightarrow} R
\]
and
\[
\phi' : C ' = A \lformal X'_1, \ldots, X'_l, Y'_1, \ldots, Y'_n \rformal/({X'}^{\Delta} - u' t^a{Y'}^{B'})
\overset{\sim}{\longrightarrow} R
\]
over $A$,
where $l \geq 2$,
$u$ and $u'$ are units,  $a$ is a positive integer,  $\Delta = (1, \ldots, 1) \in \NN^l$ and
$B, B' \in \NN^n$.
We denote the class of $X_i$ and the class of $Y_j$ in $C$ by $x_i$ and $y_j$ respectively,
and the class of $X'_i$ and the class of $Y'_j$ in $C'$
by $x'_i$ and $y'_j$ respectively.
We assume that there are subsets $\Gamma$ and $\Gamma'$ of $\{ 1, \ldots, n\}$
such that $\Supp(B) \subseteq \Gamma$, $\Supp(B') \subseteq \Gamma'$ and that
$\prod_{j \in \Gamma'}\phi'(y'_j) = v \prod_{j \in \Gamma}\phi(y_j)$
for some $v \in R^{\times}$.
Then, 
after renumbering $X_1, \ldots, X_l$, $Y_1, \ldots, Y_n$,
$X'_1, \ldots, X'_l$ and $Y'_1, \ldots, Y'_n$,
we have $\Gamma = \Gamma'$, $B = B'$ and there are
families $\{ u_i \}_{i=1}^l$ and $\{ v_j \}_{j \in \Gamma}$
of $R^{\times}$ such that
$\phi'(x'_i) = u_i \phi(x_i)$ for all $i=1, \ldots, l$ and
$\phi'(y'_j) = v_j \phi(y_j)$ for all $j \in \Gamma$.
\end{Proposition}

\Proof
By abuse of notation, $\phi(x_i)$, $\phi(y_j)$, $\phi'(x'_i)$ and $\phi'(y'_j)$ are denoted by
$x_i$, $y_j$, $x'_i$ and $y'_j$ respectively.
First, we claim the following:

\begin{Claim}
\label{claim:lem:two:isom:unit:1}
Renumbering $X_1, \ldots, X_l$ and $X'_1, \ldots, X'_l$,
there are $z_1, \ldots, z_l \in R^{\times}$ and
$w_1, \ldots, w_l$ such that
$x'_i = z_i x_i + t^a w_i$ for all $i=1, \ldots, l$.
\end{Claim}

By (3) and (4) of Proposition~\ref{prop:primary:decomp} for $I = 0$,
\[
(t^a) = (x_1, t^a) \cap \cdots \cap (x_l, t^a) = (x'_1,  t^a) \cap \cdots \cap (x'_l, t^a)
\]
are primary decompositions with no embedded primes.
Thus, renumbering $X_1, \ldots, X_l$ and $X'_1, \ldots, X'_l$,
we have $(x_i, t^a)= (x'_i, t^a)$ for
$i = 1, \ldots, l$.
Therefore, there are $z_i,w_i, z'_i, w'_i\in R$ with
\[
\begin{cases}
x_i = x'_i z'_i + t^a w'_i \\
x'_i = x_i z_i + t^a w_i
\end{cases}
\]
Thus,
$x_i = z_i z'_i x_i + t^a (w_iz'_i+ w'_i)$,
which implies that
$(1 - z_i z'_i) x_i \equiv 0 \mod tR$.
If  $z_i z'_i \in (t, x_1, \ldots, x_l, y_1, \ldots, y_n)$,
then $(1 - z_iz'_i) \in R^{\times}$.
Thus, $x_i \equiv 0 \mod tR$.
This is a contradiction because $l \geq 2$. Therefore, $z_iz' _i\in R^{\times}$, i.e.,
$z_i, z'_i \in R^{\times}$. 

\medskip
By using (5) and (6) of Proposition~\ref{prop:primary:decomp},
if we set $I =  (y_{1} \cdots y_{r}) = (y'_{1}\cdots y'_{r'})$,
then
\begin{align*}
I & = \bigcap_{j \in \Gamma \setminus \Supp(B)} (y_j) \cap
\bigcap_{j \in \Supp(B)} (x_1, y_j) \cap
\cdots \cap (x_l, y_j) \\
& = \bigcap_{j \in \Gamma \setminus \Supp(B')} (y'_j) \cap
\bigcap_{j \in \Supp(B')} (x'_1, y'_j) \cap
\cdots \cap (x'_l, y'_j)
\end{align*}
are primary decompositions with no embedded primes.
Here we claim the following;

\begin{Claim}
\label{claim:lem:two:isom:unit:2}
\begin{enumerate}
\renewcommand{\labelenumi}{(\roman{enumi})}
\item
$(y_j) \not= (x'_{i'}, y'_{j'})$ for $j \in \Gamma$,
$i' \in \{1, \ldots, l \}$ and $j' \in \Gamma'$.
Moreover, $(y'_{j'}) \not= (x_i, y_j)$
for $j' \in \Gamma'$,
$i \in \{ 1, \ldots, l \}$ and $j \in \Gamma$.

\item
If $(x_i, y_j) = (x'_{i'}, y'_{j'})$, then $i= i'$,
where $i, i' \in \{ 1, \ldots, l \}$, $j \in \Gamma$ and $j' \in \Gamma'$.

\item
If $(x_1, y_j) = (x'_1, y'_{j'})$ for some $j \in \Gamma$ and $j' \in \Gamma'$, 
then
$(x_i, y_{j}) = (x'_i, y_{j'})$ for all $i=1, \ldots, l$.
\end{enumerate}
\end{Claim}

(i) We assume that $(y_j) = (x'_{i'}, y'_{j'})$.
Note that $x'_{i'} = 0$ in $R/(t, x'_{i'}, y'_{j'})$.
Since $x'_{i'} = z_{i'} x_{i'} + t w_{i'}$, $x_{i'} = 0$ in $R/(t, x'_{i'}, y'_{j'})$.
On the other hand, $x_{i'} \not= 0$ in $R/(t, y_j)$.
This is a contradiction.

(ii) $x_i = 0$ in $R/(t, x_i, y_j)$. Thus, in the same way as (i),
$x'_i = 0$ in $R/(t, x_i, y_j)$. On the other hand,
if $i \not= i'$, then $x'_i \not= 0$ in $R/(t, x'_{i'}, y'_{j'})$.
Thus, $i = i'$.

(iii)
By using (i), (ii) and the above primary decompositions,
we can see $(x_i, y_{j}) = (x'_i, y'_{j''})$ for some $1 \leq j'' \leq r'$.
Thus, $(t, x_1, \ldots, x_l, y_{j}) = (t, x'_1, \ldots, x'_l,y'_{j''})$.
On the other hand, since $(x_1, y_j) = (x'_1, y'_{j'})$,
we have
$(t, x_1, \ldots, x_l, y_{j}) = (t, x'_1, \ldots, x'_l, y'_{j'})$.
Therefore, $j' = j''$.

\medskip
By the above claim, renumbering
$Y_{1}, \ldots, Y_n$ and $Y'_{1}, \ldots, Y'_n$,
we have $\Gamma = \Gamma'$. Moreover,
$(y_j) = (y'_j)$ if $j \not\in \Supp(B)$ and
$(x_1, y_j) = (x'_1, y'_j), \ldots, (x_l, y_j) = (x'_l, y'_j)$
if $j \in \Supp(B)$.
In particular, $(y_j) = (y'_j)$ for all $j \in \Gamma$,
namely, there is $v_j \in R^{\times}$ with $y'_j = v_j y_j$.

Here, using Claim~\ref{claim:lem:two:isom:unit:1},
\begin{align*}
t^a \phi'(u') {y'}^{B'} & = x'_1 \cdots x'_l \\
& = (z_1 x_1 + t^a w_1) \cdots (z_l  x_l + t^a w_l) \\
& = z_1 \cdots z_l x_1 \cdots x_l + 
t^a \sum_i w_i z_1 \cdots \widehat{z_i} \cdots z_l x_1 \cdots \widehat{x_i}
\cdots x_l + t^{2a} h \\
& = t^a z_1 \cdots z_l \phi(u) y^B + 
t^a \sum_i w_i z_1 \cdots \widehat{z_i} \cdots z_l x_1 \cdots \widehat{x_i}
\cdots x_l + t^{2a} h \\
\end{align*}
for some $h \in R$.
Therefore, we get
\[
\phi'(u') {y'}^{B'} =
z_1 \cdots z_l \phi(u) y^B +
 \sum_i w_i z_1 \cdots \widehat{z_i} \cdots z_l x_1 \cdots \widehat{x_i}
\cdots x_l + t^{a} h.
\]
Hence, since $y'_j = v_j y_j$ for $j \in \Gamma$,
\[
\phi'(u') \left( \prod_{j \in \Gamma} v_{j}^{B'(j)} \right) y^{B'} \equiv
z_1 \cdots z_l \phi(u) y^B \quad\mod (t, x_1, \ldots, x_l).
\]
Note that $(A/tA) \lformal Y \rformal 
\overset{\sim}{\longrightarrow} R/ (t, x_1 , \ldots, x_l)$.
Thus, $B = B'$.
Therefore, if we set $J = (t^a y^{B}) = (t^a {y'}^{B'})$, then,
by using (3) and (4) of Proposition~\ref{prop:primary:decomp} for $I = B$,
we have two primary decompositions of $J$:
\[
J  = \bigcap_{i=1}^l \left(( x_i, t^a) \cap \bigcap_{j \in \Supp(B)}
(x_i, y_j^{B(j)})  \right)
= \bigcap_{i=1}^l \left( (x'_i, t^a) \cap 
\bigcap_{j \in \Supp(B)}
(x'_i, {y'}_j^{B(j)})  \right).
\]
First of all, $(x_i, t^a) = (x'_i, t^a)$ for all $i = 1, \ldots, l$.
Moreover, if $(x_i, y_j) = (x'_{i'}, y'_{j'})$
for $i, i' \in \{ 1, \ldots, l\}$ and $j, j' \in \Gamma$, then
$i = i'$ and $j = j'$ because 
\[
\begin{cases}
\text{$x'_i \equiv z_i x_i \mod tR$ for $i= 1, \ldots, l$}, \\
\text{$y'_j = v_j y_j$ for $j = 1, \ldots, r$}, \\
(A/tA) \lformal X_1, \ldots, \widehat{X_i}, \ldots, X_l, Y_1, \ldots, \widehat{Y_j}, \ldots, Y_n \rformal
\simeq R/(t, x_i,y_j), \\
(A/tA) \lformal X'_1, \ldots, \widehat{X'_{i'}}, \ldots, X'_l, Y'_1, \ldots, \widehat{Y'_{j'}}, \ldots, 
Y'_n \rformal
\simeq R/(t, x'_{i'}, y'_{j'}).
\end{cases}
\]
Therefore, by using the uniqueness of primary decomposition,
$(x_i, y_j^{B(j)}) = (x'_{i}, {y'}_{j}^{B(j)})$
for $1 \leq i \leq l$ and $j \in \Supp(B)$.
Hence
\begin{align*}
(x_i) & = (x_i, t^a y^B) =  (x_i, t^a) \cap \bigcap_{j \in \Supp(B)} (x_i, y_{j}^{B(j)}) \\  
& = (x'_i, t^a) \cap \bigcap_{j \in \Supp(B)} (x'_i, {y'}_{j}^{B(j)}) 
= (x'_i, t^a {y'}^B)
= (x'_i).
\end{align*}
Thus, there is $u_i \in R^{\times}$ with $x'_i = u_i x_i$
for $i=1, \ldots, l$.
\QED

\bigskip
Let us start the proof of Theorem~\ref{thm:extension:log:mor}.
First, we claim the following:

\begin{Claim}
\label{claim:thm:extension:log:mor:1}
We may assume that $A/tA$ is algebraically closed and
there is a fine and sharp monoid
$Q$ together with a homomorphism
$\pi_Q : Q \to M_{S, s}$ such that
the induced homomorphism
$Q \to M_{S,\bar{s}} \to \overline{M}_{S,\bar{s}}$ is bijective.
\end{Claim}

Let us denote the special point of $S = \Spec(A)$ by $s$.
By Proposition~\ref{prop:split:log:structure},
there are a local homomorphism
$f : (A, tA) \to (A', n)$ of noetherian local rings and a
fine and sharp monoid $Q$ together with
a homomorphism
$\pi_Q : Q \to M_{\Spec(A'), n} = M_{S, s} \pushout_{A^{\times}} {A'}^{\times}$
such that $f$ is flat and quasi-finite and that the induced homomorphism
$Q \to \overline{M}_{\Spec(A'), \overline{n}}$ is bijective.
Let $P$ be a minimal prime of $A'$ such that $\Spec(A'/P) \to \Spec(A)$ is surjective.
Moreover, let $A_1$ be the localization of the normalization of $A'/P$ at a closed point.
Then, $A_1$ is a discrete valuation ring by \cite[Theorem~11.7]{MatComm}.
Let $t_1$ be the uniformizing parameter of $A_1$.
Moreover, \cite[Theorem~29.1]{MatComm},
there is a ring extension $A_2$ of $A_1$ such that
$(A_2, t_1A_2)$ is a discrete valuation ring and $A_2/t_1A_2$ is algebraically closed.
Therefore, by (2) of Theorem~\ref{thm:descent:log:smooth},
we get our claim.

\medskip
Let $x$ be a closed point of $X$. By the rigidity theorem,
it is sufficient to construct a homomorphism in a Zariski neighborhood of $x$.
By the above claim, there are fine and sharp monoids
$P$ and $P'$ together with
homomorphisms $\pi_{P} : P \to M_{X, \bar{x}}$,
$\pi_{P'} : P' \to M'_{X, \bar{x}}$,
$f : Q \to P$ and $f' : Q \to P'$ with the following properties:
\begin{enumerate}
\renewcommand{\labelenumi}{(\arabic{enumi})}
\item
The induced homomorphisms
$P \to \overline{M}_{X,\bar{x}}$ and $P' \to \overline{M}'_{X,\bar{x}}$
are bijective.

\item The diagrams
\[
\begin{CD}
Q @>{f}>> P \\
@V{\pi_Q}VV @VV{\pi_P}V \\
M_{S, \bar{s}} @>>> M_{X,\bar{x}} 
\end{CD}\qquad
\begin{CD}
Q @>{f'}>> P' \\
@V{\pi_Q}VV @VV{\pi_{P'}}V \\
M_{S, \bar{s}} @>>> M'_{X,\bar{x}} 
\end{CD}
\]
are commutative.

\item
$A \otimes_{A[Q]} A[P] \to \OO_{X, \bar{x}}$ and
$A \otimes_{A[Q]} A[P'] \to \OO_{X,\bar{x}}$ are smooth.
\end{enumerate}

\medskip
First we assume $f$ is smooth at $x$.
Then, $P = f(Q) \times \NN^r$ and $P' = f'(Q) \times \NN^{r'}$ for some
non-negative integers $r$ and $r'$.
Let $Y_1, \ldots, Y_r$ and $Y'_1, \ldots, Y'_{r'}$ be all irreducible elements
of $\NN^r$ and $\NN^{r'}$ respectively.
Since $\OO_{X, \bar{x}}$ is smooth over $A \otimes_{A[Q]} A[P]$ and
$A \otimes_{A[Q]} A[P']$, adding indeterminates
$Y_{r+1}, \ldots, Y_n$ and $Y'_{r'+1}, \ldots, Y'_{n}$, we have
two isomorphisms
\[
\phi : \widehat{A} \lformal Y_1, \ldots, Y_n \rformal \overset{\sim}{\longrightarrow}
\widehat{\OO}_{X, \bar{x}}
\quad\text{and}\quad
\phi' : \widehat{A} \lformal Y'_1, \ldots, Y'_n \rformal \overset{\sim}{\longrightarrow}
\widehat{\OO}_{X, \bar{x}}.
\]
Then,
$\phi'(Y'_1) \cdots \phi'(Y'_{r'}) = v \phi(Y_1) \cdots \phi(Y_r)$
for some $v \in \widehat{\OO}_{X, \bar{x}}^{\times}$
because $\BB_{M_X/M_S} = \BB_{M'_X/M_S}$.
Thus, renumbering $Y_1, \ldots, Y_r$ and $Y'_1, \ldots, Y'_{r'}$,
we have $r = r'$ and there are
$v_1, \ldots, v_r \in  \widehat{\OO}_{X, \bar{x}}^{\times}$
with $\phi'(Y'_i) = v_i \phi(Y_i)$ for all $i =1, \ldots, r$.
Note that by Artin's approximation theorem \cite{Artin}, we may assume that
$v_1, \ldots, v_l \in \OO_{X, \bar{x}}^{\times}$.
Here, we define
\[
H : P \times \OO_{X,\bar{x}}^{\times} \to P' \times \OO_{X,\bar{x}}^{\times}
\]
to be
\[
H((f(q),  Y^{I}), u) = ((f'(q), {Y'}^I), v^I u).
\]
Thus, we get a homomorphism $h : M_{X,\bar{x}} \to M'_{X,\bar{x}}$
such that the following diagram is commutative:
\[
\begin{CD}
P \times \OO_{X,\bar{x}} @>{\sim}>> M_{X,\bar{x}} \\
@V{H}VV @VV{h}V \\
P' \times \OO_{X,\bar{x}} @>{\sim}>> M'_{X,\bar{x}}.
\end{CD}
\]
Therefore, by (1) of Theorem~\ref{thm:descent:log:smooth},
$h : M_{X,\bar{x}} \to M'_{X,\bar{x}}$
descends to a homomorphism $M_X \to M'_X$ around
a Zariski neighborhood of $x$

\medskip
Next we assume that $f$ is not smooth at $x$.
Note that $Q \to P$ and $Q \to P'$ have semistable structures
$(\sigma, q_0, \Delta, B)$ and
$(\sigma', q'_0, \Delta', B')$.
Clearly, by our assumption, $q_0 = q'_0$.
Moreover, we have $l = \#\Supp(\Delta) = \#\Supp(\Delta') \geq 2$
because $l$ is the multiplicity of the special fiber at $x$.
We set
\[
\begin{cases}
\Supp(\Delta) = \{ X_1, \ldots, X_l \}, \\
\sigma \setminus \Supp(\Delta) = \{ Y_1, \ldots, Y_r \}, \\
\Supp(\Delta') = \{ X'_1, \ldots, X'_l \}, \\
\sigma' \setminus \Supp(\Delta') = \{ Y'_1, \ldots, Y'_{r'} \}
\end{cases}
\]
Moreover, we set $\beta(q_0) = t^a u$,
where $\beta : M_{S, \bar{s}} \to A^h$ is the canonical homomorphism,
$u$ is a unit of $A^h$ and $a$ is a positive integer.
Then,
\[
A^h \otimes_{A^h[Q]} A^h[P] = A^h[X, Y]/(X^{\Delta} - ut^a Y^B)
\]
and
\[
A^h \otimes_{A[^hQ]} A^h[P'] = A^h[X', Y']/({X'}^{\Delta} - ut^a {Y'}^{B'}).
\]
Therefore, adding indeterminates $Y_{r+1}, \ldots, Y_n$ and
$Y'_{r'+1}, \ldots, Y'_n$, we have two isomorphisms
\[
\phi : C = \widehat{A} \lformal X_1, \ldots, X_l, Y_1, \ldots, Y_n  \rformal/(X^{\Delta} - u t^aY^B)
\overset{\sim}{\longrightarrow} \widehat{\OO}_{X, x}
\]
and
\[
\phi' : C' = \widehat{A} \lformal X'_1, \ldots, X'_l, Y'_1, \ldots, Y'_n \rformal/({X'}^{\Delta} - u' t^a{Y'}^{B'})
\overset{\sim}{\longrightarrow} \widehat{\OO}_{X, x}
\]
over $\widehat{A}$.
We denote the class of $X_i$ and the class of $Y_j$ in $C$ by $x_i$ and $y_j$ respectively,
and the class of $X'_i$ and the class of $Y'_j$ in $C'$
by $x'_i$ and $y'_j$ respectively.
Since $\BB_{M_X/M_S} = \BB_{M'_X/M_S}$,
we can see that $\phi'(y'_1) \cdots \phi'(y'_{r'}) = v
\phi(y_1) \cdots \phi(y_r)$ for some $v \in  \widehat{\OO}_{X, x}^{\times}$.
Thus, by Proposition~\ref{prop:two:isom:unit},
after renumbering $X_1, \ldots, X_l$, $Y_1, \ldots, Y_{r}$,
$X'_1, \ldots, X'_l$ and $Y'_1, \ldots, Y'_{r'}$,
we have $r = r'$, $B = B'$ and there are
families $\{ u_i \}_{i=1}^l$ and $\{ v_j \}_{j =1}^r$
of $\widehat{\OO}_{X,x}^{\times}$ such that
$\phi'(x'_i) = u_i \phi(x_i)$ for all $i=1, \ldots, l$ and
$\phi'(y'_j) = v_j \phi(y_j)$ for all $j =1, \ldots, r$.
Note that 
\[
\phi(x_1), \ldots, \phi(x_l), \phi(y_1), \ldots, \phi(y_r),
\phi'(x'_1), \ldots, \phi'(x'_l), \phi'(y'_1), \ldots, \phi'(y'_r)
\in \OO_{X,\bar{x}}
\]
Thus, using Artin's approximation theorem \cite{Artin}, we may assume that
\[
u_1, \ldots, u_l, v_1, \ldots, v_r \in \OO_{X,\bar{x}}^{\times}.
\]
On the other hand,
\begin{align*}
v_1^{B(1)} \cdots v_r^{B(r)} u t^a \phi(y)^B & =
ut^a \phi'(y')^{B'} = \phi'(x'_1) \cdots \phi'(x'_l) \\
& = u_1 \cdots u_l \phi(x_1) \cdots \phi(x_l) = u_1 \cdots u_l ut^a \phi(y)^B.
\end{align*}
Therefore, $u_1 \cdots u_l = v_1^{B(1)} \cdots v_r^{B(r)}$ because
$ut^a \phi(y)^B$ is a regular element.
Hence, we can define
\[
H : P \times \OO_{X,\bar{x}}^{\times} \to P' \times \OO_{X,\bar{x}}^{\times}
\]
to be
\[
\begin{cases}
H(X_i, 1) = (X'_i, u_i) & i =1, \ldots, l,\\
H(Y_j, 1) = (Y'_j, v_j) & j =1, \ldots, r, \\
H(f(q), 1) = (f'(q), 1) & q \in Q, \\
H(1, u) = (1, u) & u \in \OO_{X,\bar{x}}^{\times}.
\end{cases}
\]
We need to see that this is well-defined. Indeed,
\begin{align*}
H(X_1 \cdots X_l, 1)&  = (X'_1 \cdots X'_l, u_1 \cdots u_l) \\
& = (f'(q_0) \cdot {Y'}_1^{B(1)} \cdots {Y'}_r^{B(r)}, v_1^{B(1)} \cdots v_r^{B(r)} ) \\
& = H(f(q_0) \cdot  Y_1^{B(1)} \cdots Y_r^{B(r)}, 1).
\end{align*}
Thus, we get a homomorphism $h : M_{X,\bar{x}} \to M'_{X,\bar{x}}$
such that the following diagram is commutative:
\[
\begin{CD}
P \times \OO_{X,\bar{x}} @>{\sim}>> M_{X,\bar{x}} \\
@V{H}VV @VV{h}V \\
P' \times \OO_{X,\bar{x}} @>{\sim}>> M'_{X,\bar{x}}.
\end{CD}
\]
This descends to a homomorphism $M_X \to M'_X$ around
a Zariski neighborhood of $x$ by (1) of Theorem~\ref{thm:descent:log:smooth}.
\QED

\renewcommand{\theTheorem}{\arabic{section}.\arabic{subsection}.\arabic{Theorem}}
\renewcommand{\theClaim}{\arabic{section}.\arabic{subsection}.\arabic{Theorem}.\arabic{Claim}}
\renewcommand{\theequation}{\arabic{section}.\arabic{subsection}.\arabic{Theorem}.\arabic{Claim}}

\section{Semistable schemes of log canonical type}
\subsection{Moderate semistable schemes}
\setcounter{Theorem}{0}

Let $f : X \to S$ be a semistable scheme over a locally noetherian scheme $S$.
For $x \in X$ and $s = f(x)$,
we say $f$ is {\em moderate at $x$} if there are a smooth $\OO_{S,s}$-algebra $B$ and
$X_1, \ldots, X_n, G \in B$ with the following properties:
\begin{enumerate}
\renewcommand{\labelenumi}{(\arabic{enumi})}
\item
$dX_1, \ldots, dX_n$ form a free basis of $\Omega_{B/\OO_{S,s}}$ and
$G \in m_{S,s} B$.

\item
There is an \'{e}tale neighborhood $(U,x')$ of $X$ at $x$ together with
an \'{e}tale morphism $\rho : U \to \Spec(B/(X_1 \cdots X_l - G))$.

\item
$dG \in B dX_{l+1} + \cdots + B dX_n$.
\end{enumerate}
Note that if $f$ is smooth at $x$, then $f$ is moderate at $x$.

First, let us see the following proposition.

\begin{Proposition}
\label{prop:moderate:criterion}
Let $(A, m_A)$ be a noetherian ring such that
$A/m_A$ is algebraically closed.
Let $X \to \Spec(A)$ be a semistable scheme over $A$.
\begin{enumerate}
\renewcommand{\labelenumi}{(\arabic{enumi})}
\item
We assume that $(A, m_A)$ is a discrete valuation ring with
a uniformizing parameter $t$.
For a closed point $x \in X$, if $\OO_{X,x}$ is regular,
then $f$ is moderate at $x$.

\item
If $(A, m_A)$ is complete and $\dim f =1$, then $f$ is moderate at any closed point of $X$.
\end{enumerate}
\end{Proposition}

\Proof
(1) By Proposition~\ref{prop:semistable:alg:closed}, there is
a surjective homomorphism
\[
\phi : \widehat{A} \lformal X_1, \ldots, X_n \rformal \to \widehat{\OO}_{X, x}
\]
with $\Ker(\phi) = (X_1 \cdots X_l - G)$
for some $G \in t\widehat{A} \lformal X_1, \ldots, X_n \rformal$.
Here we set $I = (X_1 \cdots X_l - G)$,
$M = (t, X_1, \ldots, X_n)$ and $G = \sum_{J \in \NN^n} a_J X^J$.
We assume $a_{(0, \ldots, 0)} \in t^2 \widehat{A} \lformal X_1, \ldots, X_n \rformal$.
Then, since $G - a_{(0, \ldots, 0)} \in M^2$, $G \in M^2$.
Thus, $I \subseteq M^2$. Therefore,
\[
m_{X, x}/m_{X, x}^2 \simeq M/(I + M^2) = M/M^2.
\]
This contradicts to the assumption that $\OO_{X, x}$ is regular.
Hence 
\[
a_{(0, \ldots, 0)} \in t\widehat{A} \lformal X_1, \ldots, X_n \rformal \setminus t^2\widehat{A} \lformal X_1, \ldots, X_n \rformal.
\]
Thus, there is $u \in \widehat{A} \lformal X_1, \ldots, X_n \rformal^{\times}$
with $G = t \cdot u$.
Therefore, replacing $X_1$ by $u^{-1}X_1$,
we may take $G$ as $t$. 
Hence, by Proposition~\ref{prop:approximation:etale},
$f$ is moderate at $x$.

\medskip
(2) By Proposition~\ref{prop:semistable:alg:closed}, there is
a surjective homomorphism
\[
\phi : A \lformal X_1, X_2 \rformal \to \widehat{\OO}_{X, x}
\]
with $\Ker(\phi) = (X_1X_2  - G)$
for some $G \in m_A \lformal X_1, X_n \rformal$.
We denote $\phi(X_1)$, $\phi(X_2)$ and $\phi(G)$ by
$x_1$, $x_2$ and $g$.
Here we claim the following:

\begin{Claim}
There are sequences $\{ h_n \}_{n=1}^{\infty}$,
$\{ h'_n \}_{n=1}^{\infty}$ and $\{ a_n \}_{n=1}^{\infty}$
with the following properties:
\begin{enumerate}
\renewcommand{\labelenumi}{(\alph{enumi})}
\item
$h_n, h'_n \in m_A^n \widehat{\OO}_{X, x}$
for all $n \geq 1$.

\item
$a_n \in m_A^n$
for all $n \geq 1$.

\item
$\left(x_1 + \sum_{i=1}^n h_i\right)\left(x_2 + \sum_{i=1}^n h'_i\right)
-\sum_{i=1}^n a_i \in m_A^{n+1} \widehat{\OO}_{X, x}$
for all $n \geq 1$.
\end{enumerate}
\end{Claim}

First $g$ can be written as a form $g = a_1 - x_1 h'_1 - x_2 h_1$ with
$a_1 \in m_A$ and $h_1, h' _1 \in m_A \widehat{\OO}_{X, x}$.
Then,
\[
 (x_1 + h_1)(x_2 + h'_1) - a_1 = h_1 h'_1 \in m_A^2 \widehat{\OO}_{X, x}.
\]
We assume that
$\{ h_1, \ldots, h_n \}$,
$\{ h'_1, \ldots, h'_n \}$ and
$\{ a_1, \ldots, a_n \}$ have been constructed.
Then,
we set
\[
 \left(x_1 + \sum_{i=1}^n h_i\right)\left(x_2 + \sum_{i=1}^n h'_i\right)
-\sum_{i=1}^n a_i = a_{n+1} - x_1 h'_{n+1} - x_2 h_{n+1},
\]
where $h_{n+1}, h'_{n+1} \in m_A^{n+1} \widehat{\OO}_{X, x}$ and
$a_{n+1} \in m_A^{n+1}$.
Then
\begin{multline*}
 \left(x_1 + \sum_{i=1}^{n+1} h_i\right)\left(x_2 + 
\sum_{i=1}^{n+1} h'_i\right)
-\sum_{i=1}^{n+1} a_i =  \\
 h'_{n+1} \sum_{i=1}^{n}h_i + h_{n+1} 
\sum_{i=1}^{n} h'_i + h_{n+1} h'_{n+1} \in m_A^{n+2} \widehat{\OO}_{X, x}.
\end{multline*}
Thus, we get the claim.

Here we set $h = \sum_{i=1}^{\infty} h_i$ and
$h' = \sum_{i=1}^{\infty} h'_i$ in $\widehat{\OO}_{X,x}$,
and $a = \sum_{i=1}^{\infty} a_i$ in $A$.
Then
\[
\begin{cases}
 a = (x_1 + h)(x_2 + h') \\
 \text{$h, h' \in m_A \widehat{\OO}_{X,x}$ and 
$a \in m_A$}.
\end{cases}
\]
Let us take 
$H, H' \in m_A \lformal X_1, X_2 \rformal$ such that
$\phi(H) = h$ and $\phi(H') = h'$.
Then, 
\[
\phi((X_1 + H)(X_2 + H') - a) = 0.
\]
Hence, 
by Lemma~\ref{lem:flatness:x:1:x:l} and
Lemma~\ref{lem:flat:isom},
we can see that
\[
\Ker(\phi) = ((X_1 + H)(X_2 + H')  - a).
\]
Therefore, replacing $X_1$ by $X_1 + H$ and $X_2$ by $X_2 + H'$,
we may assume $G \in m_A$.
Thus, by Proposition~\ref{prop:approximation:etale},
$f$ is moderate at $x$.
\QED

\begin{Remark}
\label{rem:log:smooth:moderate}
Let $f : X \to S$ be a semistable scheme over a locally noetherian scheme $S$.
Let $M_X$ and $M_S$ be fine log structures on $X$ and $S$ respectively.
We assume that $f : X \to S$ extends to an integral and log smooth
morphism $(f, h) : (X, M_X) \to (S, M_S)$.
Let $x \in X$ and $s = f(x)$.
If there are a fine and sharp monoid $Q$ and a homomorphism
$\pi_Q : \to M_{S,\bar{s}}$
such that $Q \to M_{S,\bar{s}} \to \overline{M}_{S,\bar{s}}$
is bijective, then $f$ is moderate at $x$.

Indeed, there is a fine and sharp monoid $P$ together with
homomorphisms $\pi_P : P \to M_{X,\bar{x}}$ and
$Q \to P$ such that $P \to M_{X,\bar{x}} \to \overline{M}_{X,\bar{x}}$
is bijective and the diagram
\[
\begin{CD}
Q @>>> P \\
@V{\pi_Q}VV @VV{\pi_P}V \\
M_{S,\bar{s}} @>{h_{\bar{x}}}>> M_{X,\bar{x}}
\end{CD}
\]
is commutative.
Then,
\[
\OO_{S,\bar{s}} \otimes_{\OO_{S,\bar{s}}[Q]} \OO_{S,\bar{s}}[P] \to
\OO_{X,\bar{x}}
\]
is smooth.
If $f$ is smooth at $x$, then the assertion is obvious, so that we assume that
$f$ is not smooth at $x$.

In the case where $Q \to P$ splits, $P \simeq Q \times N$ and
$N$ is isomorphic to the monoid arising from the monomials of
$\ZZ[U, V]/(U^2 - V^2)$ by the local structure theorem (cf. Theorem~\ref{thm:local:structure:theorem}).
Thus,
\[
\OO_{S,\bar{s}} \otimes_{\OO_{S,\bar{s}}[Q]} \OO_{S,\bar{s}}[P] \simeq
\OO_{S,\bar{s}}[U, V]/(U^2 - V^2) \simeq
\OO_{S,\bar{s}}[X_1, X_2]/(X_1X_2)
\]
because $2$ is invertible in $\OO_{S,\bar{s}}$.
This means that $f$ is moderate at $x$.

In the case where $Q \to P$ does not split,
$P$ has a semistable structure $(\sigma, q_0, \Delta, B)$ over $Q$
by the local structure theorem.
Here $\sigma$ is the set of all irreducible elements of $P$ not coming from $Q$,
$q_0 \in Q$, and $\Delta, B \in \NN^{\sigma}$.
We set $\sigma = \{ X_1, \ldots, X_l, Y_1, \ldots, Y_r  \}$ such that
$\Supp(\Delta) = \{ X_1, \ldots, X_l\}$ and
$\Supp(B) \subseteq \{ Y_{1}, \ldots, Y_{r} \}$.
Moreover, we set $b_j = B(Y_j)$ for all $j$ and $t = \alpha(\pi_Q(q_0))$,
where $\alpha : M_{S,\bar{s}} \to \OO_{S,\bar{s}}$ is the canonical
homomorphism. Then,
\[
\OO_{S,\bar{s}} \otimes_{\OO_{S,\bar{s}}[Q]} \OO_{S,\bar{s}}[P] \simeq
\OO_{S,\bar{s}} [X_1, \ldots, X_l, Y_1, \ldots, Y_r]/
(X_1 \cdots X_l - t Y_1^{b_1} \cdots Y_r^{b_r}).
\]
Therefore, $f$ is moderate at $x$.
\end{Remark}

Next we consider the following proposition.

\begin{Proposition}
\label{prop:moderate:dualizing:sheaf}
Let $(A, tA)$ be a complete discrete valuation ring such that
$A/tA$ is algebraically closed.
Let $f : X \to S = \Spec(A)$ be a semistable scheme over $S$.
We assume that $X \to S$ is moderate at any closed point of $X$.
Let $\mu : Y \to X$ be a generically \'{e}tale morphism over $S$ such that
$\mu$ is \'{e}tale in the generic fiber and $f \circ \mu$ is smooth.
Then, there is the canonical homomorphism
$\mu^*(\omega_{X/S}) \to \omega_{Y/S}$.
\end{Proposition}

\Proof
Let $X_0$ be the locus of points over which $\mu$ is \'{e}tale.
We set $Y_0 = \mu^{-1}(X_0)$.
Then, there is the canonical homomorphism
\[
\alpha_0 : \mu^*(\omega_{X_0/A}) \to \omega_{Y_0/A},
\]
so that we need to show that
$\alpha_0$ extends to
$\alpha : \mu^*(\omega_{X/A}) \to  \omega_{Y/A}$.
Note that the extension is uniquely determined if it exists.

\begin{Claim}
\label{claim:prop:moderate:dualizing:sheaf:1}
Let $\{ \pi_i : Y_i \to Y \}_{i}$ be a family of  \'{e}tale morphisms
with $\codim(Y \setminus \bigcup_i \pi_i(Y_i)) \geq 2$.
If the assertion holds for each $\mu \circ \pi_i : Y_i \to X$ for all $i$,
then so does for $\mu : Y \to X$.
\end{Claim}

We assume that there is the canonical homomorphism
\[
\alpha_i : (\mu \circ \pi_i)^*(\omega_{X/A}) \to \omega_{Y_i/A}.
\]
for each $i$.
Using the uniqueness of the extension and descent theory,
we have $\alpha' : \rest{\mu^*(\omega_{X/A})}{Y'} \to \rest{\omega_{Y/A}}{Y'}$,
where $Y' = \cup_i \pi_i(Y_i)$.
Here $Y$ is normal. Thus, we obtain the extension
$\alpha : \mu^*(\omega_{X/A})\to \omega_{Y/A}$.

\medskip
By the above claim,
we may assume that
the central fiber of $f \circ \mu$ is irreducible.
Let $\gamma$ be the generic point of the central fiber.
We set $\delta = \mu(\gamma)$.
Let us choose a closed point $x \in \overline{\{ \delta \}}$.
Then there are a smooth $A$-algebra $B$ and
$T_1, \ldots, T_n, G \in B$ with the following properties:
\begin{enumerate}
\renewcommand{\labelenumi}{(\arabic{enumi})}
\item
$dT_1, \ldots, dT_n$ form a free basis of $\Omega_{B/A}$ and
$G \in t B$.

\item
There is an \'{e}tale neighborhood $\pi : (U,x') \to X$ at $x$ together with
an \'{e}tale morphism $\rho : U \to \Spec(B/(T_1 \cdots T_l - G))$.

\item
$dG \in B dT_{l+1} + \cdots + B dT_n$.
\end{enumerate}
Since $U \times_{X} Y$ is an \'{e}tale neighborhood of $\gamma$,
by the previous claim again, 
we may assume that $U = X$ and $x' = x$.
Let $s$ be a uniformizing parameter of a discrete valuation ring
$\OO_{Y, \gamma}$.
Then,
we can set ${\mu}^*(t_i) = s^{a_i} u_i$ ($u_i \in \OO_{Y, \gamma}^{\times}$)
for each $i = 1, \ldots, l$, where
$t_i$ is the class of $T_i$ in $B$.
Then, there is a Zariski open set $V$ of $\overline{\{ \gamma \}}$
such that $u_i \in \OO_{Y, y}^{\times}$ for all $y \in V$ and
all $i = 1, \ldots, l$.
Let $c : \bigwedge^{n-1} \Omega^1_{U/A} \to
\omega_{U/A}$ be the canonical homomorphism and
$\omega_0$ a basis of $\omega_{U/A}$ as in
Lemma~\ref{lem:canonical:hom:dualizing}.
Then, by using Lemma~\ref{lem:canonical:hom:dualizing},
\[
c(d t_2 \wedge \cdots \wedge d t_n) = 
t_2 \cdots t_l \cdot \omega_0.
\]
On the other hand,
\[
\mu^*(dt_i) = s^{a_i} d u_i
\]
for all $i = 1, \ldots, l$.
Thus,
\[
u_2 \cdots u_l \cdot \mu^*(d t_2 \wedge \cdots \wedge d t_n)
= \mu^*(t_2 \cdots t_l) \cdot
d u_2  \wedge \cdots \wedge d u_l \wedge d \mu^*(t_{l+1}) \wedge \cdots \wedge
d \mu^*(t_{n}).
\]
Therefore,
\[
\mu^*(\omega_0) =
\frac{d u_1}{u_1}  \wedge \cdots \wedge \frac{d u_l}{u_l}
\wedge d \mu^*(t_{l+1}) \wedge \cdots \wedge d \mu^*(t_{n})
\]
on $V$. Thus, we get our proposition.
\QED

\begin{Corollary}
\label{cor:log:smooth:canonical}
Let $(A, tA)$ be a discrete valuation ring and
$f : X \to S = \Spec(A)$ a semistable scheme over $S$.
Let $M_X$ and $M_S$ be log structures of $X$ and
$S$.
We assume that $f$ extends to a smooth and integral
homomorphism $(X, M_X) \to (S, M_S)$.
Let $\mu : Y \to X$ be a generically \'{e}tale morphism over $S$ such that
$\mu$ is \'{e}tale on the generic fiber and $f \circ \mu : Y \to S$ is smooth.
Then, there is the canonical homomorphism
$\mu^*(\omega_{X/S}) \to \omega_{Y/S}$.
\end{Corollary}

\Proof
By Proposition~\ref{prop:split:log:structure} and
Remark~\ref{rem:log:smooth:moderate},
it is sufficient to see the following lemma.
\QED

\begin{Lemma}
\label{lem:descent:hom:dual}
Let $f : X \to S$ and $g : Y \to S$ be semistable schemes over
a locally noetherian scheme $S$, and $\phi : X \to Y$ be
a morphism over $S$.
Let $Y_0$ be the set of all points $y \in Y$ such that
$g$ is smooth at $y$, and let $X_0 = \phi^{-1}(Y_0)$ and $\phi_0 = \rest{\phi}{X_0}$.
Let $\pi : S' \to S$ be a faithfully flat and quasi-compact morphism of
locally noetherian schemes.
We set $X' = X \times_S S'$, $X'_0 = X_0 \times_S S'$,
$Y' = Y \times_S  S'$, $Y'_0 = Y \times_S S'$
$\phi' = \phi \times_S \operatorname{id}_{S'}$ and
$\phi'_0 = \rest{\phi'}{X'_0}$.
If there is a homomorphism
${\phi'}^*(\omega_{Y'/S'}) \to \omega_{X'/S'}$ as an extension
of the canonical homomorphism ${\phi'_0}^*(\omega_{Y'_0/S'}) \to \omega_{X'_0/S'}$,
then it descents to
${\phi}^*(\omega_{Y/S}) \to \omega_{X/S}$.
\end{Lemma}

\Proof
Note that if we have an extension ${\phi}^*(\omega_{Y/S}) \to \omega_{X/S}$
of the canonical homomorphism
$\phi_0^*(\omega_{Y_0/S}) \to \omega_{X_0/S}$, then
it is uniquely determined.
Thus, our lemma is a consequence of the standard descent theory.
\QED

\subsection{Semistable schemes of log canonical type and automorphisms}
\setcounter{Theorem}{0}
\label{subsec:log:canonical:auto}
Let $(A, tA)$ be a discrete valuation ring.
Let $f : X \to \Spec(A)$ be a semistable scheme over $A$.
We assume that $f$ is generically smooth, that is,
smooth over the generic point of $\Spec(A)$.
Note that $X$ is normal by using Serre's criterion.
Let $\mu : Y \to X$ be a birational morphism of normal schemes over
$\Spec(A)$ such that $f \circ \mu$ is smooth over $\Spec(A)$ and
$\mu$ is an isomorphism on the generic fibers.
Let $\Delta_{\mu}$ be the set of irreducible components of
the special fiber of $f \circ \mu$ such that  $\overline{\mu(\Gamma)}$
has the codimension greater than or equal to $2$ in $X$.
Then, there is an integer $a_{\Gamma}$ for each $\Gamma \in \Delta_{\mu}$
with
\[
\omega_{Y/A} = \mu^*(\omega_{X/A}) + \sum_{\Gamma \in \Delta_{\mu}}
a_{\Gamma} \Gamma.
\]
This $a_{\Gamma}$ is denoted by $d_{\Gamma}(\mu)$.
Let $B$ be an effective Cartier divisor on $X$.
The pair $(X, B)$ is said to be {\em of log canonical type}
if, for any birational morphism $\mu : Y \to X$ of normal schemes over
$\Spec(A)$ such that $f \circ \mu$ is smooth over $\Spec(A)$ and
$\mu$ is an isomorphism on the generic fibers,
$d_{\Gamma}(\mu) \geq \ord_{\Gamma}(\mu^*(B))$ for all
$\Gamma \in \Delta_{\mu}$.

A moderate semistable scheme is of log canonical, namely,
we have the following:

\begin{Proposition}
\label{prop:moderate:log:canonical}
Let $f : X \to \Spec(A)$ be the same as above.
If there is a local flat homomorphism
$(A, tA) \to (A', t'A')$ of discrete valuation rings such that
$X \times_{\Spec(A)} \Spec(A') \to \Spec(A')$ is moderate at any closed points,
then $(X, 0)$ is of log canonical.
\end{Proposition}

\Proof
This is a consequence of Proposition~\ref{prop:moderate:dualizing:sheaf} and
Lemma~\ref{lem:descent:hom:dual}.
\QED

Here we consider the following lemma.

\begin{Lemma}
\label{lem:log:canonical:dualizing:isom}
Let $(A, tA)$ be a discrete valuation ring. Let
$f : X \to \Spec(A)$ and $f' : X' \to \Spec(A)$ be proper and generically smooth semistable 
schemes over $\Spec(A)$, and let $B$ and $B'$ be horizontal effective Cartier divisors
on $X$ and $X'$ respectively. 
Let $\phi : X \dasharrow X'$ be a birational map over $\Spec(A)$ such that
it is an isomorphism on the generic fibers $X_{\eta}$ and $X'_{\eta}$ and that
$\phi_{\eta}^*(B'_{\eta}) = B_{\eta}$.
Let $Y$ be the normalization of the closure of the graph of $\phi$ and
let $\mu : Y \to X$ and $\mu' : Y \to X'$ be the canonical morphisms. We set
$g = f \circ \mu = f' \circ \mu'$ as the following diagram:
\[
\xymatrix{
 & Y \ar[ld]_{\mu} \ar[dd]^>>>>>>>{g} \ar[rd]^{\mu'} \ & \\
 X \ar[rd]_f \ar@{.>}'[r]^>>>>>>{\phi}[rr] & & X' \ar[ld]^{f'} \\
 & \Spec(R) & \\
}
\]
Let $\Delta$ be the set of irreducible components of the special fiber
of $g$. Let us set $\Delta_0$, $\Delta_1$ and $\Delta_2$ as follows:
\[
\begin{cases}
\Delta_1 = \{ \Gamma \in \Delta \mid \codim(\mu(\Gamma)) \geq 2 \}, \\
\Delta_2 = \{ \Gamma \in \Delta \mid \codim(\mu'(\Gamma)) \geq 2 \}, \\
\Delta_0 = \Delta_1 \cap \Delta_2.
\end{cases}
\]
Then, we have the following:
\begin{enumerate}
\renewcommand{\labelenumi}{(\arabic{enumi})}
\item 
If we set $Y' = Y \setminus \bigcup_{\Gamma \in \Delta_0} \Gamma$,
then there is an open set $Y_0$ of $Y'$ such that
$\codim(Y' \setminus Y_0) \geq 2$ and
$g$ is smooth on $Y_0$.

\item
We assume the following:
\[
\begin{cases}
\text{$d_{\Gamma}(\rest{\mu}{Y_0}) \geq \ord_{\Gamma}((\rest{\mu}{Y_0})^*(B))$} & 
\text{$\forall\  \Gamma \in \Delta_1 \setminus \Delta_0$}, \\
\text{$d_{\Gamma}(\rest{\mu'}{Y_0}) \geq \ord_{\Gamma}((\rest{\mu'}{Y_0})^*(B'))$} & 
\text{$\forall\  \Gamma \in \Delta_2 \setminus \Delta_0$}.
\end{cases}
\]
Then, the isomorphism
\[
\phi_{\eta}^* : H^0(X'_{\eta}, n(\omega_{X'_{\eta}/K} + B'_{\eta})) 
\overset{\sim}{\longrightarrow}
H^0(X_{\eta}, n(\omega_{X_{\eta}/K} + B_{\eta}))
\]
on the generic fibers gives rise to
the isomorphism
\[
H^0(X', n(\omega_{X'/A} + B')) \overset{\sim}{\longrightarrow}
H^0(X, n(\omega_{X/A} + B))
\]
for all $n > 0$.
\end{enumerate}
\end{Lemma}

\Proof
(1) Let $\Gamma \in \Delta \setminus \Delta_0$ and
$\gamma$ the generic point of $\Gamma$.
Then, by Zariski's main theorem,
either $\mu$ or $\mu'$ is an isomorphism at $\gamma$.
Therefore, $g$ is smooth at $\gamma$.
Thus, the assertion of (1) is obvious.

\bigskip
(2) For simplicity, $\rest{\mu}{Y_0}$ and $\rest{\mu'}{Y_0}$ are
denoted by $\mu_0$ and $\mu'_0$ respectively.
Let $B_Y$ be the Zariski closure of $(\mu_0)_{\eta}^*(B_{\eta})
= (\mu'_0)_{\eta}^*(B'_{\eta})$.
Then,
\[
\mu_0^*(B) = B_Y + \sum_{\Gamma \in \Delta_1 \setminus \Delta_0} a_{\Gamma} \Gamma
\quad\text{and}\quad
{\mu'}_0^*(B') = B_Y + \sum_{\Gamma \in \Delta_2 \setminus \Delta_0} a'_{\Gamma} \Gamma.
\]
Moreover, we set
\[
\omega_{Y_0/A} = \mu_0^*(\omega_{X/A}) + 
\sum_{\Gamma \in \Delta_1 \setminus \Delta_0} b_{\Gamma} \Gamma
\quad\text{and}\quad
\omega_{Y_0/A} = {\mu'_0}^*(\omega_{X'/A}) + 
\sum_{\Gamma \in \Delta_2 \setminus \Delta_0} b'_{\Gamma} \Gamma.
\]
Therefore,
\[
\begin{cases}
\omega_{Y_0/A} + B_Y = \mu_0^*(\omega_{X/A} + B) +
\sum_{\Gamma \in \Delta_1 \setminus \Delta_0} (a_{\Gamma} - b_{\Gamma}) \Gamma \\
\omega_{Y_0/A} + B_Y = {\mu'_0}^*(\omega_{X'/A} + B') +
\sum_{\Gamma \in \Delta_2 \setminus \Delta_0} (a'_{\Gamma} - b'_{\Gamma}) \Gamma.
\end{cases}
\]
Note that $a_\Gamma - b_\Gamma \geq 0$ for $\Gamma \in \Delta_1 \setminus \Delta_0$
and
$a'_\Gamma - b'_\Gamma \geq 0$
for $\Gamma \in \Delta_2 \setminus \Delta_0$.
Therefore, we have the natural injection
\begin{align*}
H^0(X, n(\omega_{X/A} + B))  & \hookrightarrow 
H^0(Y_0, n(\omega_{Y_0/A} + B_Y)) \\
& \hookrightarrow
H^0\left(Y_0 \setminus \bigcup_{\Gamma \in \Delta_2 \setminus \Delta_0} \Gamma,
\ n(\omega_{Y_0/A} + B_Y)\right) \\
& = H^0\left(Y_0 \setminus \bigcup_{\Gamma \in \Delta_2 \setminus \Delta_0} \Gamma,
\ n{\mu'_0}^*(\omega_{X'/A} + B')\right) 
\end{align*}
Here there are open sets $U$ of 
$Y_0 \setminus \bigcup_{\Gamma \in \Delta_2 \setminus \Delta_0} \Gamma$
and $V$ of $X'$ such that $\mu'_0$ induces the isomorphism
$U \overset{\sim}{\longrightarrow} V$ and
$\codim(X' \setminus V) \geq 2$.
Therefore, we obtain the injection
\[
H^0(X, n(\omega_{X/A} + B)) \hookrightarrow
H^0(V, n(\omega_{X'/A} + B')).
\]
Moreover, since $\codim(X'\setminus V) \geq 2$,
the natural map
\[
H^0(X', n(\omega_{X'/A} + B')) \to H^0(V, n(\omega_{X'/A} + B'))
\]
is an isomorphism.
Hence, by the above observation, we can see that, for every $n > 0$,
there is an injective homomorphism
\[
\alpha : H^0(X, n(\omega_{X/A} + B)) \hookrightarrow
H^0(X', n(\omega_{X'/A} + B'))
\]
such that the following diagram is commutative:
\[
\xymatrix{
H^0(X, n(\omega_{X/A} + B)) \ar@{^{(}->}[r]^{\alpha} \ar@{^{(}->}[d] & 
H^0(X', n(\omega_{X'/A} + B')) \ar@{^{(}->}[d] \\
H^0(X_{\eta}, n(\omega_{X_{\eta}/K} + B_{\eta})) 
\ar[r]_{(\phi_{\eta}^{-1})^*}^{\sim}& 
H^0(X'_{\eta}, n(\omega_{X'_{\eta}/K} + B'_{\eta})).
}
\]
In the same way, we have also an injective homomorphism
\[
\beta : H^0(X', n(\omega_{X'/A} + B')) \hookrightarrow
H^0(X, n(\omega_{X/A} + B))
\]
for all $n > 0$ with the following commutative diagram:
\[
\xymatrix{
H^0(X', n(\omega_{X'/A} + B')) \ar@{^{(}->}[r]^{\beta} \ar@{^{(}->}[d] & 
H^0(X, n(\omega_{X/A} + B)) \ar@{^{(}->}[d] \\
H^0(X'_{\eta}, n(\omega_{X'_{\eta}/K} + B'_{\eta})) 
\ar[r]_{(\phi_{\eta})^*}^{\sim}& 
H^0(X_{\eta}, n(\omega_{X_{\eta}/K} + B_{\eta})).
}
\]
Thus, for $\lambda \in H^0(X, n(\omega_{X/A} + B))$,
$\lambda$ and $\beta(\alpha(\lambda))$ gives rise to the same element
in $H^0(X_{\eta}, n(\omega_{X_{\eta}/K} + B_{\eta}))$.
Therefore, $\lambda = \beta \circ \alpha(\lambda)$,
which means $\beta \circ \alpha = \operatorname{id}$.
In the same way, we can see that $\alpha \circ \beta = \operatorname{id}$.
Hence, we get our assertion.
\QED

\begin{Theorem}
\label{thm:log:canonical:isom}
Let $(A, tA)$ be a discrete valuation ring,
$f : X \to \Spec(A)$ and $f' : X' \to \Spec(A)$ proper and generically smooth semistable 
schemes over $\Spec(A)$, and $B$ and $B'$ horizontal effective Cartier divisors
on $X$ and $X'$ respectively. 
Let $\phi : X \dasharrow X'$ be a birational map over $\Spec(A)$ such that
it is an isomorphism on the generic fibers $X_{\eta}$ and $X'_{\eta}$ and that
$\phi_{\eta}^*(B'_{\eta}) = B_{\eta}$.
We assume that $(X, B)$ and $(X', B')$ are of log canonical type.
Then, we have the following:
\begin{enumerate}
\renewcommand{\labelenumi}{(\arabic{enumi})}
\item
Then, for all $n > 0$,
there is the isomorphism
\[
\phi^* : H^0(X', n(\omega_{X'/A} + B')) \overset{\sim}{\longrightarrow}
H^0(X, n(\omega_{X/A} + B))
\]
such that the following diagram is commutative:
\[
\xymatrix{
H^0(X', n(\omega_{X'/A} + B')) \ar[r]_{\phi^*}^{\sim} \ar@{^{(}->}[d] & 
H^0(X, n(\omega_{X/A} + B)) \ar@{^{(}->}[d] \\
H^0(X'_{\eta}, n(\omega_{X'_{\eta}/K} + B'_{\eta})) 
\ar[r]_{(\phi_{\eta})^*}^{\sim}& 
H^0(X_{\eta}, n(\omega_{X_{\eta}/K} + B_{\eta})).
}
\]

\item
If $\omega_{X/A} + B$ and $\omega_{X'/A} + B'$ are ample
with respect to $f$ and $f'$, then $\phi$ extends to an isomorphism.
\end{enumerate}
\end{Theorem}

\Proof
(1) is an immediate corollary of Lemma~\ref{lem:log:canonical:dualizing:isom}.
Let us consider (2). Since $\omega_{X/A} + B$ and $\omega_{X'/A} + B'$ are ample,
there is a positive integer $n$ such that
$n(\omega_{X/A} + B)$ and $n(\omega_{X'/A} + B')$ are very ample.
Moreover, by (1), we have the following commutative diagram:
\[
\xymatrix{
X \ar@{^{(}->}[r] \ar@{.>}[d]^{\phi}  & 
\Proj\left(\bigoplus_{m \geq 0} \Sym^m(H^0(X, n(\omega_{X/A} + B))^{\vee})\right) 
\ar[d]^{\rotatebox{90}{$\sim$}}_{(\phi^*)^{\vee}} \\
X' \ar@{^{(}->}[r] & 
\Proj\left(\bigoplus_{m \geq 0} \Sym^m(H^0(X', n(\omega_{X'/A} + B'))^{\vee})\right).
}
\]
Thus, $\phi$ must be an isomorphism.
\QED

\renewcommand{\theTheorem}{\arabic{section}.\arabic{Theorem}}
\renewcommand{\theClaim}{\arabic{section}.\arabic{Theorem}.\arabic{Claim}}
\renewcommand{\theequation}{\arabic{section}.\arabic{Theorem}.\arabic{Claim}}

\end{document}